\documentclass[11pt, reqno]{article}
\pdfoutput=1
\usepackage{amssymb,amsmath,amsthm,amsfonts,latexsym,amsbsy}
\usepackage{epsfig,esint}
\usepackage{mathrsfs}
\usepackage{mathtools,textcomp}
\usepackage{graphicx}
\usepackage{epic,eepic,wrapfig,color, ifthen, tocloft} 
\usepackage{url,cite}
\usepackage{geometry}

\usepackage[colorlinks=true, citecolor=blue, filecolor=cyan, linkcolor=red, urlcolor=magenta]{hyperref}
\usepackage{pmboxdraw}
\usepackage{verbatim}
\usepackage[normalem]{ulem}
\usepackage{enumerate}

\usepackage[bottom]{footmisc}

\makeatletter
\def\blfootnote{\gdef\@thefnmark{}\@footnotetext}
\makeatother

\def\1{{\bf 1}}
\def\nn{\nonumber}

 \def\sB {{\mathcal B}} \def\sC {{\mathcal C}}

\def\sE {{\mathcal E}} 
\def\sF {{\mathcal F}}
\def\sG {{\mathcal G}} \def\sH {{\mathcal H}} 
\def\sJ {{\mathcal J}} 
\def\sL {{\mathcal L}}

\def\sN {{\mathcal N}}

\def\sT {{\mathcal T}} \def\sU {{\mathcal U}}

\def\R {{\mathbb R}} \def\Q {{\mathbb Q}}

\def\N{{\mathbb N}} 
\def\E {{\mathbb E}}

\def \diag{{\textrm{diag}}}
\def \diam{{\textrm{diam}}}

\def \bbeta{{\pmb \beta}}

\def \VD {$\mathrm{VD}(\alpha)$}
\def \RVD {${\mathrm {RVD}}(\alpha_0)$}

\def \VDb {$\mathrm{VD}(\alpha')$}
\def \RVDb {${\mathrm {RVD}}(\alpha_0')$}

\def \TJ {${\mathrm{TJ}}(\phi)$}
\def \TJb {${\mathrm{TJ}}(\beta)$}

\def \TE {${\mathrm{TE}}(\phi)$}
\def \TEb {${\mathrm{TE}}(\beta)$}

\def \IVJ {${\mathrm{IJ}}_{2,\gamma}(\phi)$}
\def \IVJb {${\mathrm{IJ}}_{2,\gamma}(\beta)$}

\def \dTJ {${\mathrm{TJ}}_{q}(\phi)$}
\def \dTJb {${\mathrm{TJ}}_{q}(\beta)$}

\def \CS {${\mathrm{CS}}(\phi)$}

\def \S {${\mathrm{SE}}(\phi)$}
\def \Sb {${\mathrm{SE}}(\beta)$}

\def \Nash {${\mathrm{Nash}}_{\nu,b}(\phi)$}
\def \Nashb {${\mathrm{Nash}}_{\nu,b'}(\beta)$}

\def \FK {${\mathrm{FK}}_\nu(\phi)$}
\def \wFK {${\mathrm{WFK}}_\nu(\phi)$}
\def \GFK {${\mathrm{GFK}}_{\nu,b}(\phi)$}

\def \Gl {${\mathrm{LRE}}_{\kappa}(\phi)$}

\def \DUE {${\mathrm{DUE}}(\phi)$}

\def\la{{\langle}}
\def\ra{{\rangle}}

\numberwithin{equation}{section}

\def\qed{{\hfill $\Box$ \bigskip}}

\def\eps{\varepsilon}
\def\wh{\widehat}
\def\wt{\widetilde}
\def\pf{\noindent{\bf Proof. }}

\def\vp{{\varphi}}

\def\p{{\Theta}}

\newcommand{\X}{\mathfrak{X}}

\newcommand{\ifandonlyif}{\impliedby\hspace{-0.3in}\implies}

\DeclareMathOperator*{\esssup}{esup}
\DeclareMathOperator*{\essinf}{einf}

\theoremstyle{plain}
\newtheorem{thm}{Theorem}[section]
\newtheorem{lem}[thm]{Lemma}
\newtheorem{cor}[thm]{Corollary}
\newtheorem{prop}[thm]{Proposition}

\theoremstyle{definition}
\newtheorem{defn}[thm]{Definition}
\newtheorem{remark}[thm]{Remark}
\theoremstyle{remark}
\begin{document}
	
		\title{Stable characterization of diagonal  heat kernel upper bounds  for symmetric Dirichlet forms}	
	\date{}

	\author{{\bf Soobin Cho$^{1}$}\\ \\
		$^1$ Department of Mathematics, University Illinois Urbana-Champaign, \\
		Urbana, IL 61801, USA. Email: {\texttt soobinc@illinois.edu}
	}

	\maketitle
	
	\begin{abstract}We present a stable characterization of on-diagonal upper bounds for heat kernels associated with regular Dirichlet forms on metric measure spaces satisfying the volume doubling property. Our conditions include integral bounds on the jump kernel outside metric balls, a variant of the Faber-Krahn inequality, a cutoff Sobolev inequality, and an integral control of inverse square volumes of balls with respect to the jump kernel. Crucially, we do not assume that the jump kernel has a density, and we show that these assumptions are essentially optimal.

\blfootnote{\!\!\!\!\textit{2020 Mathematics Subject Classification:} 35K08, 31C25, 60J35.

}

	\end{abstract}
	\allowdisplaybreaks
	
	\hypersetup{linkcolor=black}

	\tableofcontents

	\section{Introduction}\label{ss:intro}	
Heat kernels are fundamental in analysis, geometry, and probability, serving as powerful tools for studying stochastic processes and partial (integro-)differential equations. A classical example in the Euclidean setting is the Gauss–Weierstrass kernel on $\R^d$,
	\begin{align*}
		p(t,x,y)= \frac{1}{(4\pi t)^{d/2} } \exp \bigg( - \frac{|x-y|^2}{4t} \bigg), \quad t>0, \; x,y \in \R^d,
	\end{align*}
	which arises from the heat equation  $\partial_t u = \Delta u$. In more general contexts—such as on Riemannian manifolds or fractals—exact formulas of this kind typically do not exist. Consequently, heat kernel estimates have become a major theme in both analysis and probability theory.

		This paper deals with the stable characterization of on-diagonal upper bounds for  heat kernels of symmetric regular Dirichlet forms,  including those with both local and non-local parts, on general metric measure spaces. Such estimates often serve as a foundation for more refined results in spectral theory and stochastic processes. Their connections to functional inequalities and geometric features of the underlying space have been extensively studied, particularly in the settings of Riemannian manifolds and fractals.

On a Riemannian manifold $(M,\mu)$,  Carlen, Kusuoka, and Stroock \cite{CKS87}  showed that
 for the heat kernel of the Laplace-Beltrami operator,  the on-diagonal estimate
	\begin{align}\label{e:on-diagonal-alpha}
		p(t,x,x) \le Ct^{-\alpha/2} \quad \text{for all $t>0$ and $x\in M$}
	\end{align}  is equivalent to the \textit{Nash inequality} \cite{Na58}:
		\begin{align}\label{e:Nash}
		\bigg( \int_{M} |\nabla f|^2 \, d\mu  \bigg) \bigg( \int_{M} |f| \,d\mu \bigg)^{4/\alpha}\ge C	\bigg( \int_{M} f^2 \, d\mu \bigg)^{2+2/\alpha} \quad \text{for all $f \in C_c^\infty(M)$}.
	\end{align} Grigor'yan \cite{Gr94} provided an alternative perspective by showing that  \eqref{e:on-diagonal-alpha} is also equivalent to  the \textit{Faber–Krahn inequality}:
	\begin{align*}
		\lambda_1(U) \ge C \mu(U)^{-2/\alpha},
	\end{align*}
	for all precompact open subsets $U$ of $M$, where $\lambda_1(U)$ denotes the first eigenvalue of the Dirichlet Laplacian on $U$.
	
However, for manifolds with inhomogeneous geometric properties, heat kernel estimates often reflect spatial dependencies, and  \eqref{e:on-diagonal-alpha}-type estimate may fail. Li and Yau \cite{LY86} established that on a complete Riemannian manifold  with non-negative Ricci curvature, the heat kernel satisfies
	\begin{align}\label{e:on-diagonal-manifold}
		\frac{C^{-1}}{V(x, \sqrt t)}\le 	p(t,x,x) \le \frac{C}{V(x, \sqrt t)}, \qquad t>0, \; x \in M,
	\end{align} 
	where $V(x,r)=\mu(B(x,r))$ is the volume of the geodesic ball $B(x,r)$ centered at $x$. This relationship between heat kernel behavior and volume growth reflects the deep interplay between diffusion and geometry.
	 Grigor'yan \cite{Gr94} showed that, 	under  the \textit{volume doubling property}  (VD):
	\begin{align}\label{e:doubling}
		V(x,2r) \le CV(x,r) \quad \text{for all $x\in M$ and $r>0$},
	\end{align}
	 the upper bound in  \eqref{e:on-diagonal-manifold} is equivalent to the  \textit{relative} Faber-Krahn inequality:
	\begin{align*}
		\lambda_1(U) \ge \frac{C}{r^2} \bigg( \frac{V(x,r)}{\mu(U)}\bigg)^\nu,
	\end{align*}
	 for any ball $B(x,r)$  and open set $U\subset B(x,r)$. Similar behavior beyond \eqref{e:on-diagonal-alpha} has been observed on fractals \cite{BK01, HKK02}, where the heat kernel exhibits inhomogeneous decay. Significant progress on heat kernel estimates in metric measure spaces including mainfolds and fractals has been made for symmetric Dirichlet forms, often by using suitable functional inequalities. See, for example,  \cite{BB92, BB04, BBK06, BCS15, CG98, FHK94, GH14, GHL09,  GT12, Ki04, Ki09, Ki12, Sa92, St95,  Va2}.

	Let $(M,d,\mu)$ be a metric measure space in which all metric balls are precomact, and  let $(\sE,\sF)$ be a regular symmetric Dirichlet form on $L^2(M,\mu)$. The associated heat semigroup $(P_t)_{t\ge 0}$ is given by  $P_t = e^{-t\sL}$, where  $\sL$ is the $L^2$-generator of $(\sE,\sF)$.  If  $P_t$ admits an integral kernel $p(t,\cdot,\cdot)$, this kernel is called the \textit{heat kernel} of $(\sE,\sF)$. 
	
	A widely used form of on-diagonal upper bound in  general spaces is
	\begin{align}\label{e:DUE-beta}\tag{DUE$(\beta)$}
		p(t,x,x) \le \frac{C}{V(x,t^{1/\beta})}, \qquad t>0, \; x \in M,
	\end{align}where $\beta>0$ can differ from  $2$. This estimate is connected to the following relative Faber-Krahn inequality with   parameter $\beta$:
	\begin{align}\label{e:Faber-Krahn}\tag{FK$_\nu(\beta)$}
		\lambda_1(U) \ge \frac{C}{r^\beta} \bigg( \frac{V(x,r)}{\mu(U)}\bigg)^\nu,
	\end{align}
for any ball $B(x,r)$  and open set $U\subset B(x,r)$.
 For fractals like the Sierpi\'ski gasket or carpet, and the Vicsek fractal, $\beta$ can exceed $2$, introducing challenges like the lack of suitable cutoff functions in the domain of the Dirichlet form.
Hence, to connect  FK$_\nu(\beta)$ with DUE$(\beta)$  across all $\beta$, an additional assumption is typically required.

	For strongly local Dirichlet forms, Andres and Barlow \cite{AB15} introduced  the \textit{cutoff Sobolev inequality} associated with parameter $\beta$ (CS$(\beta)$), a weaker variant of functional inequalities from  \cite{BB04, BBK06}. Under VD and  the \textit{reverse volume doubling property} (RVD):
	\begin{align}\label{e:reverse-doubling}
		\frac{V(x,R)}{V(x,r)} \ge  c\left( \frac{R}{r} \right)^{\alpha_0} \quad \text{for  all $x \in M$ and $0<r \le R$}
	\end{align}
	for some $c,\alpha_0>0$, they proved that
	 if  CS$(\beta)$ holds, then
	\begin{align}\label{e:characterization-1}
		\text{DUE$(\beta)$} \; \ifandonlyif \; \text{FK$_\nu(\beta)$}.
	\end{align}
	Both CS$(\beta)$ and  FK$_\nu(\beta)$ is stable under rough isometries. A similar result was obtained in \cite{GHL15} under the \textit{generalized capacity condition} (Gcap$(\beta)$), which plays a role similar to CS$(\beta)$.
	
	Parallel questions arise for non-local Dirichlet forms, where one often imposes an upper bound on the jump kernel. A standard assumption is that the jump kernel has a density $J(x,y)$ on $M\times M$ and there exists $C>0$ such that
	\begin{align}\label{e:J-upper}\tag{J$_\le(\beta)$}
		J(x,y) \le \frac{C}{V(x,d(x,y)) d(x,y)^\beta} \quad \text{for all $x,y \in M$}.
	\end{align} 
	 Chen, Kumagai and Wang \cite{CKW-memo} showed that   under VD, RVD and a variant of  CS$(\beta)$,  if J$_\le(\beta)$ holds, then  \eqref{e:characterization-1} holds.   Grigor'yan,  Hu and  Hu \cite{GHH18} obtained a similar result under  Gcap$(\beta)$ for  Alfhors $\alpha$-regular spaces, and \cite{CKW-adv} extends this to more general Dirichlet forms that may include both local and non-local parts.
	
One then naturally asks if the  pointwise condition $J_\le(\beta)$ can be weakened.  A candidate is  the  \textit{tail estimate of the jump kernel}: the jump kernel has the form $J(dx,dy)=J(x,dy)\mu(dx)$ and there exists $C>0$ such that
		\begin{align}\label{e:TJ-upper}\tag{TJ$(\beta)$}
J(x,B(x,r)^c)=	\int_{B(x,r)^c}	J(x,dy) \le \frac{C}{r^\beta} \quad \text{for all $x\in M$ and $r>0$}.
	\end{align} 
Under VD, one sees that  J$_\le(\beta)$ implies TJ$(\beta)$. One may then ask:
	
	\medskip
	
	\noindent {\bf Question 1.}  Assume $(M,d,\mu)$ satisfies VD and RVD. Does the equivalence \eqref{e:characterization-1} hold under TJ$(\beta)$ and CS$(\beta)$ (or Gcap$(\beta)$)?

	\medskip
	
	Recently,  Grigor'yan,  Hu and  Hu \cite{GHH24-annalen} made progress by considering the following $L^q$-tail condition for jump kernels: the jump kernel has a density $J(dx,dy)=J(x,y)\mu(dx)\mu(dy)$, and there exist  $C>0$ and $q\in [1,\infty)$ such that
	\begin{align}\label{e:TJ-q-upper}\tag{TJ$_q(\beta)$}
		\bigg(\int_{B(x,r)^c}	J(x,y)^q \mu(dy) \bigg)^{1/q}  \le \frac{C}{V(x,r)^{1-1/q}r^\beta} \quad \text{for all $x\in M$ and $r>0$}.
	\end{align}
	 They proved that for any $q\ge 2$,
	\begin{align}\label{e:GHH24-annalen}
		\text{VD   $+$  FK$_\nu(\beta)$  $+$   Gcap$(\beta)$  $+$  TJ$_q(\beta)$  $\implies$  DUE$(\beta)$}.
	\end{align}

	In this paper, we consider regular Dirichlet forms, whose jump kernel may not possesses a density, and we show that \textbf{Question 1} has a negative answer.   We  then introduce a new condition that measures the inhomogeneity of space through an integrated jump bound: 
		\begin{align}\label{e:IJ-beta}\tag{IJ$_{2,\gamma}(\beta)$}
		\int_{B(x,2R^{1/\beta})\setminus B(x,R^{1/\beta})} \frac{ J(x,dy)}{\sqrt{V(y,r^{1/\beta})}} \le \frac{C}{R\sqrt{V(x,r^{1/\beta})}} \bigg( \frac{R}{r}\bigg)^\gamma \quad \text{for all $x\in M$ and $r>0$},
	\end{align}
	with some constants $C>0$ and $\gamma \ge 0$. Note that under  uniform comparability of volume of balls (i.e., there exist $V:(0,\infty) \to (0,\infty)$ and $C>1$ such that $C^{-1}V(r) \le V(x,r) \le CV(r)$ for all $x\in M$ and $r>0$), TJ$(\beta)$ implies IJ$_{2,\gamma}(\beta)$  with $\gamma=0$. More generally, under VD, TJ$(\beta)$ implies IJ$_{2,\gamma}(\beta)$ for some finite $\gamma\ge 0$.
	
One of 	our main results in this setting can be stated as follows:
\begin{thm}\label{t:main1-special}
	Assume  that $(M,d,\mu)$ satisfies {\rm VD},  and  $(\sE,\sF)$ satisfies {\rm CS$(\beta)$} and  {\rm TJ$(\beta)$}.  	For any $\nu>0$ and $\gamma \ge 0$ satisfying
	\begin{align}\label{e:main1-special}
		(1-\nu)\gamma< 1+\nu,
	\end{align} 
	we have
	\begin{align*}
		\text{\rm FK$_\nu(\beta)$   $+$  IJ$_{2,\gamma}(\beta)$   $\implies$  DUE$(\beta)$.} 
	\end{align*}
\end{thm} For the optimality of the additional assumption  IJ$_{2,\gamma}(\beta)$  with parameter range \eqref{e:main1-special}, we prove the following:
\begin{thm}\label{t:main1-special-counterexample}
	For any  $\eps>0$, there exists a metric measure space $(M,d,\mu)$ satisfying {\rm VD}, and a regular Dirichlet form $(\sE,\sF)$ on $(M,d,\mu)$ such that {\rm  CS$(\beta)$}, {\rm TJ$(\beta)$}, {\rm FK$_\nu(\beta)$}  and {\rm IJ}$_{2,\gamma}(\beta)$ hold   with constants $\nu>0$ and $\gamma\ge 0$ satisfying $(1-\nu)\gamma<1+\nu + \eps$, and
	yet {\rm DUE$(\beta)$}  does not hold.
\end{thm}
In particular, this gives a negative answer to {\bf Question 1}.

Our final result considers the case where the jump density exists, providing an improvement of \eqref{e:GHH24-annalen}.  Note that \eqref{e:doubling} is equivalent to the existence of  $\alpha>0$ and $C\ge 1$ such that
 	\begin{align}\label{e:doubling-alpha}
 	\frac{V(x,R)}{V(x,r)} \le  C\left( \frac{R}{r} \right)^{\alpha} \quad \text{for  all $x \in M$ and $0<r \le R$},
 \end{align}
 and that \eqref{e:reverse-doubling} is trivially satisfied for $\alpha_0=0$.
\begin{thm}\label{t:main2-special}
	Assume that  $(M,d,\mu)$ satisfies \eqref{e:doubling-alpha} and \eqref{e:reverse-doubling} with $\alpha>0$ and $\alpha_0\ge 0$.  For any $\nu>0$ and $q \in [1,\infty)$ satisfying
	\begin{align}\label{e:main2-special-assumption}
		\frac{(1-\nu)}{2\beta} \left[ \bigg( \frac{2}{q}-1 \bigg)  \alpha-\alpha_0\right]_+<1+\nu,
	\end{align} 
	we have
	\begin{align*}
		\text{\rm FK$_\nu(\beta)$   $+$  CS$(\beta)$  $+$  TJ$_q(\beta)$    $\implies$   DUE$(\beta)$.} 
	\end{align*}
\end{thm}

Compared with \eqref{e:GHH24-annalen}, this result admits a broader range for $q$, since \eqref{e:main2-special-assumption} always holds for $q\ge 2$ and also covers certain $q<2$.

Theorems \ref{t:main1-special}  and \ref{t:main2-special} follow immediately from Theorems \ref{t:main1} and \ref{t:main2}. The proof of Theorem \ref{t:main1-special-counterexample} will be treated in Section \ref{ss:proofs}. For the proofs of Theorems \ref{t:main1-special} and \ref{t:main2-special},  we use truncation methods developed in \cite{GHL14, GHL10, BGK09, CKW-memo}. However, since we will work with rather general scale functions and assume our conditions only hold within a certain localization radius, significant new comparisons and self-improvement arguments for heat semigroups are required. For Theorem \ref{t:main1-special-counterexample}, we explicitly construct a counterexample.

We  remark that, to achieve a stable structure, we will also consider the reverse implication DUE($\beta$)$\implies$FK$_\nu(\beta)$. To address this, we introduce a weaker version of  FK$_\nu(\beta)$ that allows us to proceed without assuming RVD with $\alpha_0>0$. Indeed, without RVD  with $\alpha_0>0$,
the implication DUE($\beta$)$\implies$FK$_\nu(\beta)$ does not hold in general (see \cite[comment 5, p. 9]{GH14}). 
A similar strategy is adopted in \cite{BCS15}.

The paper is organized as follows: In Section \ref{ss:main-results}, we present the detailed framework and our main results.  Section \ref{ss:3} provides  preliminary results. In Section \ref{ss:4}, we establish basic results for Faber-Krahn-type inequalities and explore their connections with \textit{Nash-type} inequalities. These results lead to the existence of (Dirichlet) heat kernels.  In Section \ref{ss:5}, we study the $L^2$-mean value inequality for subharmonic functions,  lower resolvent estimates and survival estimates in our framework. Survival estimates are critical for the study of truncated Dirichlet forms in later sections.  In Section \ref{ss:special}, we consider a specialized setup in which truncation arguments are simplified. This section is the most challenging and provides most of the proofs for our main results in this specialized setting.  In Section \ref{ss:7}, we extend the results established in Section \ref{ss:special} to the general setting by applying a change-of-metric method. This approach draws inspiration from \cite{Ki12, BM18, GHH24+}. Section \ref{ss:proofs} contains the proofs of our main results.  In Section \ref{ss:9}, we present two applications of our main results: stable-like non-local operators with variable orders on $\R^d$ and non-local Dirichlet forms with singular jump kernels. These cases are not covered by previous stability results due to the absence of jump densities. In the Appendix, we recall the existence of regular, pointwise defined heat kernels and provide results for truncated Dirichlet forms. We also establish new comparison inequalities for truncated Dirichlet forms, which are used in our results.
	
	\smallskip

	{\it Notation}:  Throughout the paper,   $c_i$, $i=0,1,2,...$ are fixed in each statement and proof, and their labeling restarts for each proof. 
		For $a,b\in \R$, let $a\wedge b:=\min\{a,b\}$ and  $a\vee b:=\max\{a,b\}$.  For $p \in [1,\infty]$,  $\lVert f \rVert_p$ denotes the $L^p$-norm in $L^p(M)=L^p(M,\mu)$.   For  $B=B(x_0,r)$ and a constant $a>0$, we use $aB$ to denote $B(x_0,ar)$. For a subset $D\subset M$, $\overline D$ denotes the closure of $D$, and $D^c$ denotes the complement of $D$ in $M$. We use $\Q_+$ to denote the set of all positive rational numbers.

	\section{Framework and main results}\label{ss:main-results}
	
	Let $(M,d)$ be a locally compact separable metric space,  and let $\mu$ be a positive Radon measure on $M$ with full support. The triplet $(M,d,\mu)$ is called a \textit{metric measure space}. We assume that all metric balls in $M$ are precompact. Let
	 $(\sE, \sF)$ be a regular symmetric Dirichlet form on $L^2(M):=L^2(M,\mu)$.  By the Beurling-Deny formula,  $(\sE,\sF)$  can be decomposed into the \textit{strongly local part}, the \textit{pure-jump part} and the \textit{killing part}. In this paper, we always assume that $(\sE, \sF)$ has no killing part. Thus,
	 \begin{align*}
	 	\sE(f,g)=\sE^{(L)}(f,g) + \sE^{(J)}(f,g), \quad f,g \in \sF,
	 \end{align*}
 where $\sE^{(L)}$ is the strongly local part of $(\sE,\sF)$ and $\sE^{(J)}$ is the jump part associated with a positive symmetric Radon measure $J$ on $M \times M\setminus \diag$:
	 \begin{align*}
	 	\sE^{(J)}(f,g)&= \int_{M \times M} (f(x)-f(y))(g(x)-g(y))J(dx,dy).
	 \end{align*}
By \cite[Theorem 2.1.3]{FOT},  every $f \in \sF$ admits a quasi-continuous version $\wt f$ on $M$. In this paper, we  always use a quasi-continuous version of $f \in \sF$ without notational distinction.

To state our results precisely, we introduce some definitions. 
Denote by  $B(x,r)$  the open ball centered at $x$ with radius $r$ and  let $V(x,r):=\mu(B(x,r))$.

	\begin{defn}\label{d:VD}
	(i)  We say that $(M,d,\mu)$ satisfies the \textit{volume doubling property}  \VD, if there exist  constants $\alpha>0$ and  $C\ge 1$  such that
			\begin{align*}
		 \frac{V(x,R)}{V(x,r)} \le  C \left( \frac{R}{r} \right)^{\alpha} \quad \text{for  all $x \in M$ and $0<r \le R$}.
			\end{align*}
				(ii)  We say that $(M,d,\mu)$ satisfies the \textit{reverse volume doubling property}  \RVD, if there exist  constants $\alpha_0\ge 0$ and  $C\ge 1$  such that
			\begin{align*}
				\frac{V(x,R)}{V(x,r)} \ge  C^{-1} \left( \frac{R}{r} \right)^{\alpha_0} \quad \text{for  all $x \in M$ and $0<r \le R<\diam(M)$}.
			\end{align*}
	\end{defn}
	
	We emphasize that the constant $\alpha_0$ in \RVD \ can be $0$, and that RVD(0) holds trivially.
	
	\begin{remark}\label{r:volume-doubling}
 (i)	Condition  \VD \ is equivalent to the following: There exist constants $\alpha>0$ and $C\ge 1$ such that
\begin{align}\label{e:VD2}
	\frac{V(y,R)}{V(x,r)} \le C\bigg(\frac{R + d(x,y)}{r}\bigg)^{\alpha}  \quad \text{for all} \;\, x,y \in M \text{ and } 0<r\le R.
\end{align} 
	\noindent (ii) If \RVD \ holds, then for any fixed $k\ge 1$, there exists a constant $C=C(k)\ge 1$ such that 
		\begin{align}\label{e:RVD2}
		\frac{V(x,R)}{V(x,r)} \ge  C^{-1} \left( \frac{R}{r} \right)^{\alpha_0} \quad \text{for  all $x \in M$ and $0<r \le R<k\,\diam(M)$}.
	\end{align}
	\end{remark}

	A function $\phi:M \times (0,\infty)\to (0,\infty)$ is called a \textit{scale function}  if:
 
	\smallskip
	
	$\bullet$  For each $x \in M$, the map $r \mapsto \phi(x,r)$ is strictly increasing and  continuous.

	$\bullet$    There exists a constant $C_1 \ge 1$ such that
	\begin{align}\label{e:phi-comp}
		\phi(y,r)\le C_1\phi(x,r)  \quad \text{for all} \;\, x,y \in M \text{ and } r \ge d(x,y).
	\end{align}
	
	$\bullet$  There exist constants $\beta_1,\beta_2>0$ and  $C_2\ge 1$  such that
	\begin{align}\label{e:phi-scale}
C_2^{-1}\left( \frac{R}{r} \right)^{\beta_1} \le 		\frac{\phi(x,R)}{\phi(x,r)} \le C_2\left( \frac{R}{r} \right)^{\beta_2} \quad \text{for all} \;\, x \in M \text{ and } 0<r\le R.
	\end{align}
	Denote by $\phi^{-1}(x,\cdot)$  the inverse function of $r \mapsto \phi(x,r)$. Throughout this paper,  we always use $\phi$ to denote a scale function.

Let $\sB(M)$ be the $\sigma$-alebra of Borel sets in $M$. A function $J:M \times \sB(M)\to [0,\infty)$ is called a \textit{transition jump kernel} if:

\smallskip

$\bullet$ For each $x \in M$, $J(x,dy)$ is a Borel measure.

$\bullet$  For every $E \in \sB(M)$, the map $x \mapsto J(x,E)$ is Borel measurable.

\begin{defn}
(i)	We say that   \TJ \ holds if there exist a  transition jump kernel $J(x,dy)$ and  $C>0$ so that 
	\begin{align}\label{e:J-density}
		J(dx,dy)=J(x,dy)\mu(dx) \quad \text{in $M\times M$}
	\end{align}
	and
	\begin{equation*}
		J(x, B(x,r)^c) \le \frac{C}{\phi(x,r)} \quad \text{for all $x\in M$ and $r>0$}.
	\end{equation*}
	
\noindent 	(ii) 	For $\gamma\ge 0$,	we say that   \IVJ \  holds if \eqref{e:J-density} is satisfied, and there exists  $C>0$ such that for all $x\in M$ and $0<r\le R$,
\begin{align}\label{e:IVJ}
	\int_{B(x,2\phi^{-1}(x,R))\setminus B(x,\phi^{-1}(x,R))} \frac{ J(x,dy)}{\sqrt{V(y,\phi^{-1}(y,r))}} \le \frac{C}{R\sqrt{V(x,\phi^{-1}(x,r))}} \bigg( \frac{R}{r}\bigg)^\gamma.
\end{align}
\end{defn}

\begin{remark}\label{r:IVJ-general}
(i)	If $(M,d,\mu)$ satisfies \VD \ and \RVD, and $\phi$ satisfies \eqref{e:phi-scale}, then \IVJ \ always  holds with $\gamma= \alpha/(2\beta_1) - \alpha_0/(2\beta_2)$ (see Proposition \ref{p:IVJ-general}). 
	
\noindent (ii) If there exist functions $V:(0,\infty) \to (0,\infty)$ and  $\phi:(0,\infty) \to (0,\infty)$ such that 
\begin{align}\label{e:homoeneous}
	C^{-1} V(r) \le 	V(x,r) \le CV(r)  \quad &\text{for all $x\in M$ and $0<r <\diam(M)$},\nn\\
	C^{-1}\phi(r) \le \phi(x,r) \le C\phi(r)  \quad &\text{for all $x\in M$ and $r>0$},
\end{align}
for some constant $C>1$,	then  \TJ \ implies \IVJ \ with $\gamma=0$. Indeed,   under \eqref{e:homoeneous},  $V(y,\phi^{-1}(y,r))$ and $V(x, \phi^{-1}(x,r))$ are comparable for all $x,y \in M$ and $r>0$. Hence, the implication holds.
\end{remark}

From now on, we fix a constant $
T_0 \in (0, \infty],$ 
which  we refer to as the  \textit{localizing constant}.

For $a>0$, define $\sE_a(u,u)=\sE(u,u)+a \lVert u \rVert_2^2$. For an open set $D\subset M$, let $\sF^D$ be the $\sE_1^{1/2}$-closure of $\sF \cap C_c(D)$ in $\sF$, where $C_c(D)$ is the family of all continuous functions on $D$ with compact supports.  The form $(\sE, \sF^D)$ is called the \textit{part}  of $(\sE,\sF)$ on $D$. It is known that  $(\sE, \sF^D)$ is a regular Dirichlet form on $L^2(D)$. See \cite[Theorem 3.3.9(i)]{CF12}. Define
\begin{align}\label{e:def-lambda1}
	\lambda_1(D)= \inf \left\{ \sE(f,f) : f \in \sF^D \mbox{ with } \Vert f \Vert_{2} =1  \right\}.
\end{align}

\begin{defn}\label{d:FK} 
 	(i)	We say that the \textit{Faber-Krahn inequality} \FK \ holds if there exist constants   $\nu,C>0$ and  $\delta_1\in (0,1)$ such that for any ball $B:=B(x_0,r)$ with $0<r<\phi^{-1}(x_0,\delta_1T_0)$ and  non-empty open set $D\subset B$,
	\begin{equation*}
		\lambda_1(D) \ge \frac{C}{\phi(x_0,r)} \left( \frac{V(x_0,r)}{\mu(D)}\right)^\nu.
	\end{equation*}

 \noindent	(ii)	We say that the \textit{weak Faber-Krahn inequality} \wFK \ holds if there exist constants   $\nu, C,C'>0$ and  $\delta_2\in (0,1)$ such that for any ball $B:=B(x_0,r)$ with $ 0<r<\phi^{-1}(x_0,\delta_2T_0)$ and  non-empty open set $D\subset B$,
	\begin{equation*}
		\lambda_1(D) \ge \frac{C}{\phi(x_0,r)} \bigg[ \left( \frac{V(x_0,r)}{\mu(D)}\right)^\nu - C'\bigg].
	\end{equation*}

	\noindent (iii)	We say that the \textit{generalized Faber-Krahn inequality} \GFK \ holds if there exist constants   $\nu, C, C'>0$, $b\ge 0$ and  $\delta_2\in (0,1)$ such that for any ball $B:=B(x_0,r)$ with $r >0$ and  non-empty open set $D\subset B$,
	\begin{equation}\label{e:GFK}
		\lambda_1(D) \ge \frac{C}{\phi(x_0,r)} \bigg[  \bigg(1\wedge \frac{T_0}{\phi(x_0,r)}\bigg)^b \left( \frac{V(x_0,r)}{\mu(D )}\right)^\nu   -  C' \bigg].
	\end{equation}
\end{defn}

Observe that both   \FK \ and \GFK \   imply \wFK. Moreover, under \RVD \ with $\alpha_0>0$, it can be shown that  \FK \ and \wFK \ are equivalent (see Lemma \ref{l:wFK}).

Define $	\sF'=\left\{ f+a: f \in \sF,  a \in \R \right\}.$ The form $\sE$  extends to $\sF'$ by setting  
\begin{align*}
	\sE(f, g)=\sE(f_0,g_0), \;\; \text{
		where $f=f_0+a$, $g=g_0+b$,  $f_0,g_0 \in \sF$, $a,b\in \R$.} 
\end{align*} This definition is well-defined since  $(\sE,\sF)$ has no killing part. Let $\sF_b:=\sF \cap L^\infty(M)$ and $\sF'_b:=\sF\cap L^\infty(M)$. For  $f\in\sF_b$, there is a unique positive Radon measure  $\Gamma(f,f)$  on $M$ such that 
\begin{align*}
	\int_M g \,d \Gamma(f,f)=\sE(f,fg)-\frac12 \sE(f^2,g) \quad \text{for all $g \in \sF\cap C_c(M)$.}
\end{align*}
The measures $\Gamma(f,f)$  can be uniquely extended to any $f=f_0+a \in \sF'$ as the increasing limit of $\Gamma(((-n) \vee f_0) \wedge n,((-n) \vee f_0) \wedge n))$. The measure $\Gamma(f,f)$  is called the \textit{energy measure} or the \textit{carr\'e du champ} of $f$ for $\sE$. See \cite[Section 3.2]{FOT}. Let $\Gamma^{(L)}(f,f)$  be  the energy measure of $f$ with respect to the strongly local part $\sE^{(L)}$.

For subsets $U,V$ of $M$, we write $U \Subset V$ if $\overline U \subset V$. Let $U$ and $V$ be open subsets of $M$ such that $U \Subset V$. We say that a  measurable function $\vp$ is a \textit{cutoff function} for $U \Subset V$, if $0\le \varphi \le 1$ in $M$, $\varphi = 1$ in $U$ and $\varphi=0$ in $V^c$.

\begin{defn}
	 We say that the \textit{cutoff Sobolev inequality} \CS \ holds if there exist constants $c_1,c_2>0$ such that for any $f \in \sF_b'$, $x_0\in M$ and  $0<R<R+r<R'<\phi^{-1}(x_0,T_0)$,  there exists a cutoff function $\vp \in \sF_b$ for $B(x_0,R) \Subset B(x_0,R+r)$ so that the following holds:
			\begin{align*}
			&\int_{B(x_0,R+r)} f^2 d\Gamma^{(L)}(\vp,\vp)  + \int_{B(x_0,R')\times B(x_0,R')} f(x)^2 (\vp(x)-\vp(y))^2 J(dx,dy) \nn\\
			&\le   c_1\bigg(\int_{B(x_0,R+r)} \varphi^2 d\Gamma^{(L)}(f,f) + \int_{B(x_0,R+r)\times B(x_0,R+r)} \vp(x)^2 (f(x) - f(y))^2 J(dx,dy) \bigg)  \nn\\
			&\quad + \sup_{z \in B(x_0,R')} \frac{c_2}{\phi(z,r)} \int_{B(x_0,R')} f^2 d\mu .
		\end{align*}
\end{defn}

We next introduce   conditions  \S \ and \TE,  referred to as the \textit{survival estimates} and  \textit{tail estimates} of the semigroup $(P_t)_{t\ge 0}$.
\begin{defn}
(i)	We say that    \S \ holds if  there exist constants $\eps_0,a_0\in (0,1)$ such that for any ball $B:=B(x_0,r)$ with $r\in (0,\phi^{-1}(x_0,T_0))$,
	\begin{align*}
		\essinf_{  B(x_0,r/4)}	P^B_t \1_{B} \ge \eps_0 \quad \text{for all $t \in (0, a_0\phi(x,r)]$}.
	\end{align*}
	
\noindent (ii)	We say that    \TE \ holds if  there exists  $C>0$ such that for any ball $B:=B(x_0,r)$ with $r>0$,
	\begin{align}\label{e:TE}
		\esssup_{ B(x_0,r/4)}	P_t \1_{B^c}  \le \frac{Ct}{\phi(x_0,r) \wedge T_0}  \quad \text{for all $t>0$}.
	\end{align}
\end{defn}

Let $D\subset M$ be an open set. A subset $E$ of $D$ is said to be \textit{$D$-regular} if  $\mu(E\cap U_x \cap D)>0$ for any $x \in E$ and any open neighborhood $U_x$ of $x$. An increasing sequence of closed subsets $\{E_n\}_{n\ge 1}$ of $D$ is called a \textit{$\mu$-nest} of $D$ if $\lim_{n\to \infty} \mu(D \setminus E_n)=0$, and is \textit{$D$-regular} if each $E_n$ is $D$-regular. 
For a sequence of closed  subsets $\{E_n\}_{n\ge 1}$ of $M$, let 
$	C(\{E_n\}) :=\left\{ u: u|_{E_n} \text{ is continuous for each $n\ge 1$}\right\}.$ 

Let $(P^D_t)_{t \ge 0}$ be the semigroup  associated with the part $(\sE,\sF^D)$ of $(\sE,\sF)$ on $D$. We say that a non-negative jointly measurable function $p^D$ on $(0,\infty) \times  D \times D$ 
is a \textit{pointwise defined heat kernel} of $(P^D_t)_{t \ge 0}$, if 
there exists a $D$-regular $\mu$-nest $\{E_n\}_{n\ge 1}$  of $D$
such that the following conditions are satisfied for all $t>0$:

\setlength{\leftskip}{3mm}
\smallskip

\noindent 
$\bullet$ For  any  bounded measurable function $f$ on $D$,
$$
P^D_tf(x)= \int_{D} p^D(t,x,y) f(y) \mu(dy), \quad \text{$\mu$-a.e. $x\in D$}.
$$

\noindent  $\bullet$ (Markov) $ p^D(t,x,y)\ge 0$ and $\int_D p^D(t,x,z)\mu(dz)\le 1$ for all  $ x,y \in D$.

\noindent  $\bullet$ (Symmetry) $ p^D(t,x,y)=p^D(t,y,x) $ for all  $ x,y \in D$.

\noindent  $\bullet$ (Semigroup property) For all $s>0$ and $x,y  \in D$,
$$
p^D(t+s,x,y) =\int_{D} p^D(t,x,z) p^D(s,z,y) \mu(dz).
$$

\noindent  $\bullet$ (Integral regularity)  For each $f \in L^1(D) \cap L^2(D)$,
\begin{align*}
	\int_D p^D(t,\cdot, z) f(z) \mu(dz) \in C(\{E_n\}) \quad \text{and} \quad 	\int_D p^D(t,z, \cdot) f(z) \mu(dz) \in C(\{E_n\}) .
\end{align*}

\noindent  $\bullet$ (Pointwise regularity) For each $x \in D$,
\begin{align*}
p^D(t,x,\cdot) \in C(\{E_n\}) \cap L^1(D)\cap L^2(D) \quad \text{and} \quad p^D(t,\cdot,x) \in C(\{E_n\}) \cap L^1(D)\cap L^2(D).
\end{align*}

Moreover, 
\begin{align}\label{e:heatkernel-domain}
	p^D(t,x,y)=0 \quad \text{if $x \in D \setminus \cup \cup_{n=1}^\infty E_n$ or $y \in  D \setminus \cup \cup_{n=1}^\infty E_n$.} 
\end{align}

\setlength{\leftskip}{0mm}

We denote the heat kernel of $(P_t)_{t\ge 0}$, if it exists, by $p(t,x,y)$.
  
  Throughout this paper, a heat kernel always refers to a pointwise defined heat kernel.
\begin{defn}\label{d:DUE} 
	We say that \textit{on-diagonal upper estimate} \DUE \ holds if  $p(t,x,y)$  exists on $(0,\infty) \times M \times M$ and for each $k\ge 1$, there is $C=C(k)>0$ such that 
	\begin{align*}
	 p(t,x,x) \le  \frac{C }{  V(x,\phi^{-1}(x,t))} \quad \text{for all $x\in M$ and  $t \in (0, kT_0)$}.
	\end{align*}
\end{defn}

By the semigroup property and the Cauchy-Schwarz inequality, one sees that \DUE \ is equivalent to the following condition: For any $k\ge 1$, there is  $C=C(k)>0$ such that 
	\begin{equation}\label{e:DUE-2}
		p(t,x,y) \le  \frac{C }{ \sqrt{V(x,\phi^{-1}(x,t))V(y,\phi^{-1}(y,t))}} \quad \text{for all $x,y \in M$ and $t \in (0, kT_0)$}.
	\end{equation}

A  Dirichlet form $(\sE,\sF)$ is  \textit{conservative} if  $P_t\1_M = \1_M$ for all $t>0$.

The following is the first main result of this paper.

\begin{thm}\label{t:main1}
	Assume  that $(M,d,\mu)$ satisfies \VD,  and  $(\sE,\sF)$ satisfies \TJ \ with $\phi$  satisfying \eqref{e:phi-comp} and \eqref{e:phi-scale}. Then the following statements hold.

	\medskip
	
\noindent 	(i) \DUE \ $\implies$ {\rm GFK}$_{\nu,1}(\phi)$  $\implies$ \wFK  \, for $\nu=\beta_1/\alpha$; 

	\S \ $\implies$ \CS;

	\DUE \ $+$ \S \; $\ifandonlyif$  \ \DUE \ $+$  \TE \ and  $(\sE,\sF)$ is conservative.

\noindent (ii)
For any $\nu>0$ and $\gamma \ge 0$ satisfying
\begin{align}\label{e:main1-assumption}
	(1-\nu)\gamma< 1+\nu,
\end{align} 
we have
\begin{align}\label{e:main1-result}
	\text{\wFK  \ $+$ \CS \ $+$ \IVJ  \  $\implies$   \DUE \ $+$ \S.} 
\end{align}

	\noindent If  $(M,d,\mu)$ additionally satisfies \RVD \ with $\alpha_0>0$, then the above statements hold with \FK \ replacing \wFK.
\end{thm}

If $\gamma\le 1$,  then \eqref{e:main1-assumption} always holds for any $\nu>0$.  Using this fact, we  obtain the following corollary from Theorem \ref{t:main1}. 
\begin{cor}\label{c:main1}
		Assume  that $(M,d,\mu)$ satisfies \VD \ and \RVD \ with $\alpha_0\ge 0$,  and $\phi$  satisfies \eqref{e:phi-comp} and \eqref{e:phi-scale}.  Suppose also that one of the following holds:

		\vspace{-0.1in}

		\begin{enumerate}[(1)]

			\setlength{\leftskip}{3mm}
			
			\rm	\item $(\sE,\sF)$ is strongly local. 
			
			\vspace{-0.1in}

			\item $(\sE,\sF)$ satisfies \TJ \ and  $\alpha/\beta_1 \le \alpha_0/\beta_2 + 2$.

			\vspace{-0.1in}

			\item $(\sE,\sF)$ satisfies \TJ,  and that there exist   functions $V:(0,\infty) \to (0,\infty)$ and  $\phi:(0,\infty) \to (0,\infty)$ such that for some  $C>1$,
			\begin{align*}
			C^{-1} V(r) \le 	V(x,r) \le CV(r)  \quad &\text{for all $x\in M$ and $0<r <\diam(M)$},\\
			C^{-1}\phi(r) \le \phi(x,r) \le C\phi(r)  \quad &\text{for all $x\in M$ and $r>0$}.
			\end{align*}
		\end{enumerate}

		\vspace{-0.1in}

	\noindent Then the following statements are all equivalent:

	\vspace{-0.1in}

	\begin{enumerate}[(a)]

		\setlength{\leftskip}{3mm}

		\rm	\item \wFK  \ $+$ \CS \, for some $\nu>0$.

	\vspace{-0.1in}

	\item	\DUE \ $+$ \S.

	\vspace{-0.1in}

	\item   \wFK  \ $+$ \CS \, for $\nu=\beta_1/\alpha$.
	\end{enumerate}
\end{cor}

We now demonstrate that the parameter range in \eqref{e:main1-assumption} is essential for the implication \eqref{e:main1-result} and, in general, cannot be improved.

\begin{thm}\label{t:main1-counterexample}
For any  $\eps>0$, there exists a metric measure space $(M,d,\mu)$ satisfying \VD, a scale function $\phi$  satisfying \eqref{e:phi-comp} and \eqref{e:phi-scale}, and a regular Dirichlet form $(\sE,\sF)$ on $(M,d,\mu)$ such that  \TJ, \wFK, \CS \ and \IVJ \ hold   with constants $\nu>0$ and $\gamma\ge 0$ satisfying $(1-\nu)\gamma<1+\nu + \eps,$
yet \DUE \ fails.
\end{thm}

Our last result addresses the case where the  jump kernel has a density. In this context, the range of $\nu$ in \eqref{e:main1-assumption}  can be improved, specifically for the case  $\gamma= \alpha/(2\beta_1) - \alpha_0/(2\beta_1)$. To achieve this,
we introduce an additional  condition on the tail of the jump density, as proposed in \cite{GHH24-annalen}.

\begin{defn}
	Let $q \in [1,\infty)$.		We say that condition  \dTJ \ holds if there exist a   measurable function $J(x,y)$ on $M\times M \setminus \diag$ and a constant $C>0$ so that 
	\begin{align}\label{e:J-kernel-density}
		J(dx,dy)=J(x,y)\mu(dx)\mu(dy) \quad \text{in $M\times M\setminus \diag$}
	\end{align}
	and
	\begin{equation*}
		\bigg( \int_{B(x,r)^c} J(x,y)^{q} \mu(dy)\bigg)^{1/q}\le \frac{C}{V(x,r)^{(q-1)/q}\phi(x,r)} \quad \text{for all $x\in M$ and $r>0$}.
	\end{equation*}
\end{defn}

\begin{remark}\label{r:dTJ}
	By \cite[Proposition 3.1]{GHH24-annalen}, if \VD \ holds, then  \dTJ \ for some $q \in [1,\infty)$ implies \TJ. Moreover, under \VD, if \eqref{e:J-kernel-density} holds with $J(x,y)$ satisfying
	\begin{align*}
		J(x,y) \le \frac{C}{V(x,d(x,y)) \phi(x,d(x,y))} \quad \text{for all $x,y\in M$},
	\end{align*} then \dTJ \ holds for every $q \in [1,\infty)$.
\end{remark}

\begin{thm}\label{t:main2}
Suppose that  $(M,d,\mu)$ satisfies \VD \ and \RVD \ with $\alpha_0\ge 0$,   and  $\phi$ satisfies \eqref{e:phi-comp} and \eqref{e:phi-scale}. For any $\nu>0$ and $q \in [1,\infty)$ satisfying
	\begin{align}\label{e:main2-assumption}
			\frac{(1-\nu)}{2} \left[ \bigg( \frac{2}{q}-1 \bigg) \frac{ \alpha}{\beta_1}-\frac{\alpha_0}{\beta_2}\right]_+<1+\nu,
	\end{align} 
	we have
	\begin{align}\label{e:main2-result}
		\text{\wFK  \ $+$ \CS \ $+$ \dTJ  \  $\implies$   \DUE \ $+$ \S.} 
	\end{align}
	In particular, for any $\nu>0$,
	\begin{align}\label{e:main2-result-2}
		\text{\wFK  \ $+$ \CS \ $+$ {\rm TJ}$_{2}(\phi)$    $\implies$   \DUE \ $+$ \S.} 
	\end{align}
\end{thm}

\begin{remark}
\cite{GHH24-annalen} proves a similar result to  \eqref{e:main2-result-2} under a spatially independent localizing constant $\overline R$ instead of $\phi^{-1}(x_0, \delta T_0)$. The two frameworks  coincide if $\overline R=\diam(M)$ or if $\sup_{x \in M} \phi(x, 1) \le c \inf_{x\in M} \phi(x,1)$ for some $c\ge 1$. More precisely, with modifications to the localizing constant in the definitions of  conditions, 
\cite[Corollary 2.14]{GHH24-annalen} shows that  for any $\nu>0$ and  $q \ge 2$,
\begin{align*}
		\text{\VD \ $+$ \FK  \ $+$ \CS \ $+$ \dTJ  \  $\implies$   \DUE.} 
\end{align*}
Indeed,  \cite{GHH24-annalen} uses the \textit{generalized capacity condition} (see condition (Gcap) therein), while in this paper \CS \ is used. By following the arguments in \cite[Theorem 2.11]{GHH24} and applying Lemma \ref{l:capacity},  one can show that these two conditions are equivalent under \VD, \wFK \ and \TJ. Note that our result has an advantage in the range of $q$;  \eqref{e:main2-assumption} covers some $q$  strictly smaller than $2$.
\end{remark}

\section{Preliminaries}\label{ss:3}

Recall that $\phi$ is a scale function satisfying \eqref{e:phi-comp} and \eqref{e:phi-scale}. We begin with a basic property of $\phi$.

\begin{lem}\label{l:phi-basic}
	There exists $C>0$ such that 
	\begin{align}\label{e:phi-scale2}
		\frac{\phi(y,R)}{\phi(x,r)} \le C\bigg(\frac{R + d(x,y)}{r}\bigg)^{\beta_2}  \quad \text{for all} \;\, x,y\in M \text{ and } 0<r \le R.
	\end{align}
\end{lem}
\pf By the monotonicity of $\phi$,  \eqref{e:phi-comp} and \eqref{e:phi-scale}, we obtain for all $x,y \in M$ and $0<r\le R$,
\begin{align*}
	\frac{\phi(y,R)}{\phi(x,r)} \le \frac{\phi(y,R + d(x,y))}{\phi(x,r)} \le  \frac{c_1\phi(x,R + d(x,y))}{\phi(x,r)}\le c_2\bigg(\frac{R + d(x,y)}{r}\bigg)^{\beta_2}.
\end{align*}
\qed

The following lemma is standard (cf.  \cite[Remark 1.7]{CKW-memo}).
\begin{lem}\label{l:beta<2}
	Suppose that $\sE$ is of pure-jump type, namely, $\sE^{(L)}=0$, and $\phi$  satisfies \eqref{e:phi-scale} with $\beta_2<2$.   Then   \TJ \ implies \CS.
\end{lem}
\pf  Let $x_0 \in M$ and $0<R<R+r<R'$. Choose a cutoff function $\vp(x)=h(d(x_0,x))$ for $B(x_0,R)\Subset B(x_0,R+r)$, where $h$ is a $C^1$ function on $[0,\infty)$ such that $0 \le h\le1$, $h(s)=1$ for all $s\le R$, $h(s)=0$ for all $s\ge R+r$ and $|h'(s)|\le 2/r$ for all $s\ge 0$. By \TJ \ and \eqref{e:phi-scale},  for all $x\in M$,
\begin{align}\label{e:beta<2}
		&\int_M (\vp(x)-\vp(y))^2 J(x,dy)  \le  \int_{M} \left[ \bigg( \frac{2d(x,y)}{r}\bigg)^2  \wedge 1 \right] J(x,dy) \nn\\
		&\le   \sum_{n=0}^\infty  2^{2-2n} \int_{B(x,2^{-n}r) \setminus B(x,2^{-n+1}r)}  J(x,dy)+ \int_{B(x,r)^c} J(x,dy) \nn\\
		&\le c_1 \sum_{n=0}^\infty   \frac{2^{2-2n}}{\phi(x, 2^{-n+1}r)}+ \frac{c_1}{\phi(x,r)} \le \frac{c_1}{\phi(x,r)} \bigg(1 +  c_2\sum_{n=0}^\infty 2^{(\beta_2-2)n} \bigg) = \frac{c_3}{\phi(x,r)}. 
\end{align}
Thus, for any $f\in \sF'_b$ and $R'>R+r$, we get
\begin{align*}
	& \int_{B(x_0,R')\times B(x_0,R')} f(x)^2 (\vp(x)-\vp(y))^2 J(dx,dy)\\
	& =\int_{B(x_0,R')} f(x)^2 \int_{B(x_0,R')} (\vp(x)-\vp(y))^2 J(x,dy) \, \mu(dx) \le  \sup_{z \in B(x_0,R')} \frac{c_3}{\phi(z,r)} \int_{B(x_0,R')} f^2 d\mu .
\end{align*}\qed

\begin{prop}\label{p:IVJ-general}
	Suppose that $(M,d,\mu)$ satisfies  \VD \ and \RVD \ with $\alpha_0\ge 0$.
	 If  \TJ \ holds, then  \IVJ \ holds with $\gamma=\alpha/(2\beta_1)-\alpha_0/(2\beta_2)$.
\end{prop}
\pf  By \VD \ and \eqref{e:phi-scale}, we have for all $y\in M$ and $0<r\le R$,
\begin{align}\label{e:V-phi-1}
	\frac{V(y,\phi^{-1}(y,R))}{V(y,\phi^{-1}(y,r))} \le c_1\bigg( \frac{\phi^{-1}(y,R)}{\phi^{-1}(y,r)}\bigg)^{\alpha} \le c_2\bigg( \frac{R}{r}\bigg)^{\alpha/\beta_1}.
\end{align}
Further,  by  \eqref{e:RVD2} and \eqref{e:phi-scale},  for all $y\in M$ and $0<r\le R$ with $\phi^{-1}(y,R) < 2\,\diam(M)$,
\begin{align}\label{e:V-phi-2}
	\frac{V(y,\phi^{-1}(y,R))}{V(y,\phi^{-1}(y,r))} \ge c_3\bigg( \frac{\phi^{-1}(y,R)}{\phi^{-1}(y,r)}\bigg)^{\alpha_0} \ge c_4\bigg( \frac{R}{r}\bigg)^{\alpha_0/\beta_2}.
\end{align}

Let $x\in M$ and $0<r\le R$. If $\phi^{-1}(x,R)\ge 2 \,\diam(M)$, then \eqref{e:IVJ} clearly holds. Assume $\phi^{-1}(x,R)< 2 \,\diam(M)$.   For all $y \in B(x,2\phi^{-1}(x,R))$, by \eqref{e:phi-scale2}, we have $	\phi(y, \phi^{-1}(y,R))=	R= \phi(x, \phi^{-1}(x,R)) \ge c_5 \phi(y, \phi^{-1}(x,R)).$
By \eqref{e:phi-scale} and the monotonicity of $\phi$, it follows that $\phi^{-1}(y,R) \ge c_6\phi^{-1}(x,R)$.  Using this, \eqref{e:V-phi-1}, \eqref{e:VD2} and \eqref{e:V-phi-2}, we  obtain for all $y \in B(x,2\phi^{-1}(x,R))$,
\begin{align*}
	\frac{1}{\sqrt{V(y,\phi^{-1}(y,r))}}  &\le \frac{c_2^{1/2}(R/r)^{\alpha/(2\beta_1)}}{\sqrt{V(y,\phi^{-1}(y,R))}} \le \frac{c_7(R/r)^{\alpha/(2\beta_1)} (1+ 3\phi^{-1}(x,R)/\phi^{-1}(y,R))^{\alpha/2}}{\sqrt{V(x,\phi^{-1}(x,R) + \phi^{-1}(y,R))}}  \\
	&  \le \frac{c_8(R/r)^{\alpha/(2\beta_1)}}{\sqrt{V(x,\phi^{-1}(x,R))}}  \le \frac{c_9(R/r)^{\alpha/(2\beta_1)- \alpha_0/(2\beta_2)}}{\sqrt{V(x,\phi^{-1}(x,r))}}.
\end{align*}
Thus, by using \TJ, we arrive at
\begin{align*}
&	\int_{B(x,2\phi^{-1}(x,R))\setminus B(x,\phi^{-1}(x,R))} \frac{ J(x,dy)}{\sqrt{V(y,\phi^{-1}(y,r))}}\\
 &\le  \frac{c_9(R/r)^{\alpha/(2\beta_1)- \alpha_0/(2\beta_2)}}{\sqrt{V(x,\phi^{-1}(x,r))}} J(x,B(x,\phi^{-1}(x,R))^c )\le \frac{c_{10}(R/r)^{\alpha/(2\beta_1)- \alpha_0/(2\beta_2)}}{R\sqrt{V(x,\phi^{-1}(x,r))}} .
\end{align*}\qed

\begin{lem}\label{l:IVJ}
	Suppose that \VD \ and \IVJ \ hold. For any $a>1$, there exists  $C=C(a)>0$ such that for all $x\in M$ and $0<r\le aR$,
	\begin{align*}
		\int_{B(x,a \phi^{-1}(x,R))\setminus B(x,\phi^{-1}(x,R))} \frac{ J(x,dy)}{\sqrt{V(y,\phi^{-1}(y,r))}} \le \frac{C}{R\sqrt{V(x,\phi^{-1}(x,r))}} \bigg( \frac{R}{r}\bigg)^{\gamma}.
	\end{align*}
\end{lem}
\pf Let $n_0=n_0(a)\ge 1$ be such that $2^{n_0-1}\le a<2^{n_0}$. For all $x\in M$ and $0<r\le aR$, using \VD \ and \eqref{e:phi-scale} in the first inequality below,  \IVJ \ in the second, and \VD \ and  \eqref{e:phi-scale} in the third, we obtain 
\begin{align*}
	&\int_{B(x,a \phi^{-1}(x,R))\setminus B(x,\phi^{-1}(x,R))} \frac{ J(x,dy)}{\sqrt{V(y,\phi^{-1}(y,r))}}\\
	&\le  \sum_{n=1}^{n_0}\int_{B(x,2^n \phi^{-1}(x,R))\setminus B(x,2^{n-1}\phi^{-1}(x,R))} \frac{ J(x,dy)}{\sqrt{V(y,\phi^{-1}(y,r/a))}} \\
	&\le \frac{c_1}{(r/a)^\gamma}\sum_{n=1}^{n_0}  \frac{\phi(x, 2^{n-1} \phi^{-1}(x,R))^\gamma}{ \phi(x, 2^n \phi^{-1}(x,R))     \sqrt{V(x,\phi^{-1}(x, r/a))}}  \le  \frac{c_2 n_0 2^{(n_0-1)\beta_2 \gamma} a^{\alpha/(2\beta_1) +\gamma} }{R   \sqrt{V(x,\phi^{-1}(x, r))}} \bigg( \frac{R}{r}\bigg)^\gamma.
\end{align*}
\qed

By a standard covering argument, we establish that the constant $1/4$ in \S \  can be replaced by any constant $\delta\in (0,1)$. 

\begin{lem}\label{l:S-self-improvement}
	Suppose that 	there exist constants $\delta,\eps_0,c_0\in (0,1)$ such that for any ball $B:=B(x_0,r)$ with $r\in (0,\phi^{-1}(x_0,T_0))$,
	\begin{align}\label{e:S-self-improvement}
		\essinf_{ \delta B}	P^B_t \1_{B} \ge \eps_0 \quad \text{for all $t \in (0, c_0\phi(x,r)]$}.
	\end{align}
	Then \S \ holds.
\end{lem}
\pf   Let $B:=B(x_0,r)$  with $r\in (0,\phi^{-1}(x_0,T_0))$. By    \eqref{e:phi-scale} and \eqref{e:phi-comp}, there exists a constant $\eta\in (0,1/2)$ independent of $x_0$ and $r$ such that 
\begin{align*}
	\sup_{z\in \overline{4^{-1}B} } \phi(z_i, \eta r) \le  c_1 \eta^{\beta_1}\sup_{z\in \overline{4^{-1} B} } \phi(z_i, r)\le  \phi(x_0,r)<T_0.
\end{align*}
Let  $\{B(z_i,\delta\eta r)\}_{i\ge 1}$ be a covering of $4^{-1}B$, where $z_i \in \overline{4^{-1}B}$ for all $i\ge 1$.  Note that $B(z_i, \eta r) \subset   B$ for all $i\ge 1$.   Applying \eqref{e:S-self-improvement}, we obtain for all $t\in  (0, c_0 \inf_{i\ge 1} \phi(z_i,\eta r) ]$,
\begin{align*}
	\essinf_{4^{-1}B} P_{t}^B\1_B & \ge \inf_{i\ge 1}  	\essinf_{B(z_i,\delta\eta r)} P_{t}^B\1_B \ge \inf_{i\ge 1}  	\essinf_{B(z_i,\delta \eta r)} P_{t}^{B(z_i,\eta r)}\1_{B(z_i,\eta r)} \ge \eps_0. \end{align*}
By \eqref{e:phi-scale2}, $ \inf_{i\ge 1} \phi(z_i,\eta r)  \ge   c_2\phi(x_0,r)$ for some $c_2\in (0,1)$ independent of $x_0$ and $r$. We conclude that  \S \ holds with $a_0=c_0c_2$ and $\eps_0$. \qed

\begin{lem}\label{l:TE->S}
	Suppose that $(\sE,\sF)$ is conservative. Then \TE \ implies \S.
\end{lem}
\pf Let $B:=B(x_0,r)$  with $r\in (0,\phi^{-1}(x_0,T_0))$. By \cite[Proposition 4.1]{GHL10},  for all $t>0$,
\begin{align}\label{e:TE->S-1}
	\essinf_{12^{-1} B}  P^{B}_{t} \1_{B}\ge 	\essinf_{12^{-1} B}  P^{B}_{t} \1_{3^{-1} B} &\ge 	\essinf_{12^{-1} B} P_{t} \1_{3^{-1} B} -  \sup_{0<s\le t} \esssup_{\overline{(2/3) B}^c} P_{s} \1_{3^{-1} B}.
\end{align}
Since $(\sE,\sF)$ is conservative,   by \TE \ and the monotonicity of $\phi$, we obtain for all $t>0$, 
\begin{align}\label{e:TE->S-2}
	\essinf_{y \in 12^{-1} B}  P_{t} \1_{3^{-1} B}(y)  &= 1- 
	\esssup_{y \in 12^{-1} B}  P_{t} \1_{(3^{-1}B)^c}(y) \ge 1 - \frac{c_1t}{\phi(x_0,r)}.
\end{align}
Let $\{B(z_i, r/12)\}_{i\ge 1}$ be a covering of $\overline{(2/3) B}^c$ where $z_i \in ( (2/3)B)^c$. Note that  $d(x_0,z_i)/8 \ge r/12$ and $3^{-1}B \cap B(z_i, d(x_0,z_i)/2)= \emptyset$ for all  $i\ge 1$.   By \TE \ and \eqref{e:phi-scale2}, we get that for all  $s>0$,
\begin{align}\label{e:TE->S-3}
	&\esssup_{\overline{(2/3) B}^c} P_{s} \1_{3^{-1} B} \le \sup_{i\ge 1} \esssup_{ B(z_i,r/12)} P_{s} \1_{\frac13 B} \le \sup_{i\ge 1} \esssup_{ B(z_i,d(x_0,z_i)/8)} P_{s} \1_{B(z_i, d(x_0,z_i)/2)^c}\nn\\
	&\le \sup_{i\ge 1} \frac{c_2s}{\phi(z_i, d(x_0,z_i)/2) \wedge T_0} \le \sup_{i\ge 1} \frac{c_3s}{\phi(x_0, 3d(x_0,z_i)/2) \wedge T_0} \le   \frac{c_3s}{\phi(x_0, r)}.
\end{align}
Set $\eps_0:=1/(1+c_1+c_3)$. Combining \eqref{e:TE->S-1}, \eqref{e:TE->S-2} and \eqref{e:TE->S-3}, we obtain for all $t\in (0, \eps_0\phi(x_0,r)]$,
\begin{align*}
	\essinf_{12^{-1} B}  P^{B}_{t} \1_{B} \ge 1 - \frac{(c_1+c_3)t}{\phi(x_0,r)} \ge \eps_0.
\end{align*} 
By Lemma \ref{l:S-self-improvement}, we conclude that \S \ holds.  \qed

\begin{lem}\label{l:S->CS}
If \TJ \ holds, then 	\S \ implies \CS.
\end{lem}
\pf   The result follows from the proofs of \cite[Lemmas 6.2 and 13.5]{GHH24}.  Although \cite{GHH24} uses a spatially independent  localizing constant $\overline R \in (0,\infty]$ (see conditions (S) and (ABB) therein) and we use $\phi^{-1}(x_0, T_0)$,  the proofs remain valid.  \qed

\section{Faber-Krahn type inequalities}\label{ss:4}

\begin{lem}\label{l:wFK}
	If \RVD \ holds with $\alpha_0>0$, then \wFK \ and \FK \ are equivalent.
\end{lem}
\pf Clearly, \FK \  implies \wFK. Thus, it suffices to prove the other implication. 

Suppose that \wFK \ holds with  $\delta_2\in (0,1)$  and let $\delta_1\in (0,\delta_2)$ be a constant whose  value  to be determined later. 
Let $B:=B(x_0,r)$  with $r\in (0,\phi^{-1}(x_0,\delta_1T_0)) $ and   $D\subset B$ be  a non-empty open set. Applying \wFK \ with $B$ replaced by  $(\delta_2/\delta_1)B$ in the first inequality below, and using \eqref{e:phi-scale} and \RVD \ in the second, we obtain 
\begin{align}\label{e:wFK-1}
	\lambda_1(D) &\ge \frac{c_1}{\phi(x_0,(\delta_2/\delta_1)r)} \bigg[ \left( \frac{V(x_0,(\delta_2/\delta_1)r)}{\mu(D)}\right)^\nu - c_2\bigg]\nn\\
	&	 \ge \frac{c_3}{(\delta_2/\delta_1)^{\beta_2}\phi(x_0,r)} \bigg[ \left( \frac{c_4 (\delta_2/\delta_1)^{\alpha_0}V(x_0,r)}{\mu(D)}\right)^\nu - 1\bigg],
\end{align}
for some $c_3>0$ and $c_4 \in (0,1)$. Set $\delta_1:= (2^{-1/\nu}c_4)^{1/\alpha_0}\delta_2$. It follows from \eqref{e:wFK-1} that
\begin{align*}
	\lambda_1(D)  \ge \frac{c_5}{\phi(x_0,r)} \bigg[ 2\left( \frac{V(x_0,r)}{\mu(D)}\right)^\nu - 1\bigg]\ge \frac{c_5}{\phi(x_0,r)}  \left( \frac{V(x_0,r)}{\mu(D)}\right)^\nu .
\end{align*}This proves the lemma. \qed

The proof of the next lemma is motivated by the approach in \cite[p. 553]{GH14}, where the implication \DUE \ $\Rightarrow$ \FK \ was established under \RVD \ with $\alpha_0>0$, for the case where $T_0=\infty$ and $\phi(x,r)=r^\beta$ for some $\beta>0$. Recall that the implication \DUE \ $\Rightarrow$ \FK \ does not generally hold  without \RVD \ with $\alpha_0>0$. 

\begin{lem}\label{l:DUE-GFK}
	If \VD \ holds, then	\DUE \ implies \GFK \ with  $\nu=\beta_1/\alpha$ and $b=1$.
\end{lem}
\pf Set  $\nu:=\beta_1/\alpha$. Let  $B:=B(x_0,r)$  and    $D\subset B$ be a non-empty open set.     For a constant $c_0\ge 1$ to be chosen later,  define
$$
T=c_0\bigg( \frac{\mu(D)}{V(x_0,r)} \bigg)^{\nu} \phi(x_0, r).
$$
We consider the following two cases separately:

\smallskip

Case 1: Assume  $T <T_0$. By  \eqref{e:DUE-2}, there exists $c_1>0$ such that  for all $x,y \in B$,
\begin{align}\label{e:FK-claim-1}
	p(T,x,y) \le \frac{c_1}{\sqrt{V(x,\phi^{-1}(x,T))V(y,\phi^{-1}(y,T))}}.
\end{align}
We distinguish two subcases:

(i) Suppose $2c_1\mu(D) \ge \inf_{z\in B}V(z,r)$. Then by \eqref{e:VD2}, 
\begin{align*}
	\bigg(1\wedge \frac{T_0}{\phi(x_0,r)}\bigg) \left( \frac{V(x_0,r)}{\mu(D )}\right)^\nu   \le  \left( \frac{V(x_0,r)}{\mu(D )}\right)^\nu  \le   \left( \frac{2c_1V(x_0,r)}{\inf_{z\in B}V(z,r) }\right)^\nu   \le  c_2.
\end{align*}
By taking $C'$ larger than $c_2$, we get \eqref{e:GFK} with $b=1$.

(ii) Suppose $2c_1\mu(D)<\inf_{z\in B}V(z,r)$. We first prove that, by choosing $c_0$ large enough, 
\begin{align}\label{e:FK-claim}
	\sup_{x,y \in D} p(T,x,y) \le  \frac{1}{2\mu(D)}.
\end{align}
To prove \eqref{e:FK-claim}, by \eqref{e:FK-claim-1}, it suffices to show that $	V(x,\phi^{-1}(x,T)) \ge 2c_1 \mu(D)$ for all $x\in D$. 
Let $x \in D$.   If $T\ge \phi(x,r)$, then  $	V(x,\phi^{-1}(x,T)) \ge V(x,r) >2c_1\mu(D)$.  If $T< \phi(x,r)$, then using  \VD \ in the first inequality below, \eqref{e:phi-scale} in the second, and \eqref{e:VD2} and  \eqref{e:phi-scale2} in the last, we obtain
\begin{align*}
	V(x,\phi^{-1}(x,T)) &\ge c_3V(x,r)\bigg( \frac{\phi^{-1}(x,T)}{r}\bigg)^\alpha  \ge c_4V(x,r)\bigg( \frac{T}{\phi(x,r)}\bigg)^{1/\nu} \\
	&= \frac{c_4V(x,r) \mu(D)}{ V(x_0,r)}\bigg( \frac{c_0 \phi(x_0,r)}{\phi(x,r)}\bigg)^{1/\nu} \ge c_5c_0^{1/\nu} \mu(D).
\end{align*} By choosing $c_0:= (2c_1/c_5)^{1/\nu}+1$, we conclude that  \eqref{e:FK-claim} holds.

Fix $\eps>0$. Let  $f \in \sF^D$ be such that $\lVert f \rVert_2=1$ and $\sE(f,f)<\lambda_1(D)+\eps$. By \cite[Lemma 1.3.4(i)]{FOT},  we have
\begin{align}\label{e:FK-2}
	\sE(f,f) = \sup_{t>0} t^{-1}\langle f-P^{ D}_t f , f \rangle   &\ge T^{-1}\langle f-P^{ D}_T f , f \rangle  = T^{-1}( 1- \langle P^D_T f, f\rangle ).
\end{align}
Using \eqref{e:FK-claim} and  the Cauchy-Schwarz inequality, we get
\begin{align*}
	\langle P^D_T f, f\rangle  \le \int_{D\times D} p(T,x,y) |f(x)| |f(y)|  \mu(dx)\mu(dy) \le \frac{1}{2\mu(D)} \bigg(\int_D |f| d\mu \bigg)^2 \le \frac12.
\end{align*}
Combining this with \eqref{e:FK-2},  we arrive at $\lambda_1(D) + \eps > \sE(f,f) \ge  1/(2T)$. Since $\eps$ is arbitrary, this proves  \eqref{e:GFK}.

Case 2:  Assume $T \ge T_0$. In this case,  we have
\begin{align*}
	\bigg(1\wedge \frac{T_0}{\phi(x_0,r)}\bigg) \left( \frac{V(x_0,r)}{\mu(D )}\right)^\nu \le \bigg(1\wedge c_0\bigg( \frac{\mu(D)}{V(x_0,r)} \bigg)^{\nu}\bigg) \left( \frac{V(x_0,r)}{\mu(D )}\right)^\nu  \le c_0.
\end{align*}
Thus, by taking $C'$ larger than $c_0$, we get \eqref{e:GFK} with $b=1$.

\smallskip

The proof is complete.  \qed

\begin{lem}\label{l:GFK}
	Suppose that  \VD, \TJ, \wFK \ and \CS \ hold. Then \GFK \ holds with $b=(1+\nu)\alpha/\beta_1$.
\end{lem}
\pf Let $x_0 \in M$,  $r>0$ and $D\subset B(x_0,r)$ be a non-empty open set. If $\phi(x_0,r)<\delta_2T_0$, then \eqref{e:GFK} follows from \wFK.  

Assume $\phi(x_0,r) \ge  \delta_2T_0$.   Set $r_0:=4^{-1}\phi^{-1}(x_0,\delta_2T_0) \in (0, r/4]$. Fix  $\eps>0$ and let $f \in \sF^D$ be such that $\lVert f \rVert_2=1$ and $\sE(f,f) \le \lambda_1(D)+\eps$.   By \eqref{e:VD2}, there exists  $c_1>0$ such that $V(y,r_0) \ge c_1V(z,r_0)$ for all $y,z \in M$ with $d(y,z) \le r_0$. We claim that there exists $z \in B(x_0,r)$ such that
\begin{align}\label{e:GFK-claim}
	\frac{1}{V(z,r_0)}	\int_{B(z,r_0)} f^2 d\mu \ge  \frac{c_1}{V(x_0,r+r_0)}.
\end{align}
Indeed, if \eqref{e:GFK-claim} does not hold for all $z \in B(x_0,r)$, then  
\begin{align*}
	1 &= \int_{B(x_0,r)} f(y)^2 \mu(dy)  = \int_{B(x_0,r)} \frac{f(y)^2}{V(y,r_0)} \int_{B(y,r_0)} \mu(dz)\mu(dy)\\
	&\le \int_{B(x_0,r+r_0)} \frac{1}{c_1V(z,r_0)} \int_{B(z,r_0)} f(y)^2 \mu(dy) \mu(dz) < \frac{1}{V(x_0,r+r_0)} \int_{B(x_0,r+r_0)}  \mu(dz) =1,
\end{align*}
which is a contradiction. 

Choose $z\in B(x_0,r)$ such that \eqref{e:GFK-claim} holds. Denote $U_a:=B(z,ar_0)$ for $a>0$. By \CS,  there exists a cutoff function $\vp \in \sF_b$ for $U_1\Subset U_2$ such that
\begin{align}\label{e:GFK-1}
	&\int_{U_2} f^2 d\Gamma^{(L)}(\vp,\vp)  + \int_{U_3\times U_3} f(x)^2 (\vp(x)-\vp(y))^2 J(dx,dy) \nn\\
	&\le   c_2\bigg(\int_{U_2} \vp^2 d\Gamma^{(L)}(f,f) + \int_{U_2 \times U_2} \vp(x)^2 (f(x) - f(y))^2 J(dx,dy) \bigg) + \sup_{w \in U_3} \frac{c_3}{\phi(w,r_0)}\nn\\
	&\le c_2 \sE(f,f)  + \sup_{w \in U_3} \frac{c_3}{\phi(w,r_0)}.
\end{align}
Using the  local property of $\sE^{(L)}$, the  Leibniz rule and the Cauchy-Schwarz inequality, we get
\begin{align}\label{e:GFK-2}
	&\sE^{(L)}(f\vp,f\vp)	 = \int_{U_2} f^2d\Gamma^{(L)}(\vp,\vp)  + 2\int_{U_2} f\vp\, d\Gamma^{(L)}(f,\vp)  + \int_{U_2} \vp^2d\Gamma^{(L)}(f,f) \nn\\
	&\le 2  \int_{U_2} f^2d\Gamma^{(L)}(\vp,\vp)  + 2 \int_{U_2} d\Gamma^{(L)}(f,f) \le 2  \int_{U_2} f^2d\Gamma^{(L)}(\vp,\vp)  + 2 \sE^{(L)}(f,f) .
\end{align}
Further, using  $\vp=0$ in $U_2^c$ and \eqref{e:J-density} in the equality below, the Cauchy-Schwarz inequality in the first inequality, and $\vp^2\le 1$, \TJ \ and $\lVert f \rVert_2=1$ in the second inequality, we obtain
\begin{align}\label{e:GFK-3}
	&\sE^{(J)}(f\vp,f\vp)\nn\\
	&= \int_{U_3\times U_3} (f(x)\vp(x)-f(y)\vp(y))^2 J(dx,dy)  + 2 \int_{U_2} f(x)^2\vp(x)^2  \int_{U_3^c}J(x,dy) \mu(dx) \nn\\
	&\le  2 \int_{U_3\times U_3} f(x)^2 (\vp(x) -\vp(y))^2 J(dx,dy) +2 \int_{U_3\times U_3} \vp(y)^2 (f(x)-f(y))^2 J(dx,dy)\nn\\
	&\quad +  2 \int_{U_2} f(x)^2  J(x,B(x,r_0)^c) \mu(dx)\nn\\
	&\le  2 \int_{U_3\times U_3} f(x)^2 (\vp(x) -\vp(y))^2 J(dx,dy) +2 \sE^{(J)}(f,f)+  \sup_{w\in U_2} \frac{c_4}{\phi(w,r_0)}.
\end{align}
Combining \eqref{e:GFK-1}, \eqref{e:GFK-2} and \eqref{e:GFK-3},   we arrive at
\begin{align}\label{e:GFK-5}
	&	\sE(f\vp, f\vp) = \sE^{(L)}(f\vp,f\vp) + \sE^{(J)}(f\vp,f\vp)\nn\\
	&\le 2  \int_{U_2} f^2d\Gamma^{(L)}(\vp,\vp)   +  2 \int_{U_3\times U_3} f(x)^2 (\vp(x) -\vp(y))^2 J(dx,dy) +2 \sE(f,f)+  \sup_{w\in U_2} \frac{c_4}{\phi(w,r_0)} \nn\\
	&\le (2c_2+2) \sE(f,f) +   \sup_{w\in U_3} \frac{2c_3+c_4}{\phi(w,r_0)} \le (2c_2+2) \sE(f,f) +  \frac{c_5}{\phi(z,r_0)},
\end{align}
where we used \eqref{e:phi-scale2} in the last inequality.

On the other hand,  since $f\vp \in \sF^{D\cap U_2}$, using \wFK \ and \eqref{e:phi-scale}  in the first inequality below,   \eqref{e:GFK-claim} in the third and  \eqref{e:VD2} in the fourth,  we obtain
\begin{align}\label{e:GFK-6}
	\sE(f\vp, f\vp) &\ge \frac{c_6}{\phi(z, r_0)} \bigg[ \left( \frac{V(z,2r_0)}{\mu(D \cap U_2)}\right)^\nu -c_7\bigg] \int_{U_2} f^2 \vp^2d\mu  \nn\\
	&\ge \frac{c_6}{\phi(z,r_0)}  \left( \frac{V(z,2r_0)}{\mu(D )}\right)^\nu \int_{U_1} f^2d\mu - \frac{c_8}{\phi(z,r_0)} \int_{U_2} f^2 d\mu \nn\\
	& \ge \frac{c_1c_6V(z,r_0)}{\phi(z,r_0)V(x_0,r+r_0)}  \left( \frac{V(z,2r_0)}{\mu(D )}\right)^\nu  - \frac{c_8}{\phi(z,r_0)} \nn\\
	& \ge \frac{c_9}{\phi(z,r_0)}  \bigg(\frac{r_0}{3r}\bigg)^{(1+\nu)\alpha} \left( \frac{V(x_0,r)}{\mu(D )}\right)^\nu  - \frac{c_8}{\phi(z,r_0)} .
\end{align}
Using  \eqref{e:GFK-5}  in the second inequality below, \eqref{e:GFK-6} in the third,  and \eqref{e:phi-comp} with $z\in B(x_0,r)$ and \eqref{e:phi-scale} in the fourth,  we arrive at
\begin{align*}
	\lambda_1(D) + \eps\ge \sE(f,f)& \ge \frac{1}{2c_2+2}  \bigg(\sE(f\vp,f\vp) -   \frac{c_{5}}{\phi(z,r_0)} \bigg)_+\nn\\
	&\ge \frac{1}{\phi(z,r_0)} \bigg[ c_{10} \bigg(\frac{r_0}{r}\bigg)^{(1+\nu)\alpha} \left( \frac{V(x_0,r)}{\mu(D )}\right)^\nu   -  c_{11} \bigg]_+\nn\\
	&\ge \frac{c_{12}}{\phi(x_0,r_0)} \bigg[   c_{13} \bigg(\frac{\delta_2T_0}{\phi(x_0,r)}\bigg)^{(1+\nu)\alpha/\beta_1} \left( \frac{V(x_0,r)}{\mu(D )}\right)^\nu   -  c_{11} \bigg]_+\nn\\
	&\ge \frac{c_{12} c_{13}}{\phi(x_0,r)} \bigg[  \bigg(\frac{\delta_2T_0}{\phi(x_0,r)}\bigg)^{(1+\nu)\alpha/\beta_1} \left( \frac{V(x_0,r)}{\mu(D )}\right)^\nu   -  \frac{c_{11}}{c_{13}} \bigg].
\end{align*}
Since $\eps>0$ is arbitrary, this completes the proof. \qed

We introduce condition \Nash, a variant of the Nash inequality.
\begin{defn}
	We say that  \Nash \ holds if there exist constants $\nu,C>0$ and $b\ge 0$ such that for any ball $B:=B(x_0,r)$ with   $r>0$ and $f \in \sF^{B}$,
	\begin{align}\label{e:GFK-Nash}
		\lVert f \rVert_{2}^{2+2\nu} \le   \frac{C \phi(x_0,r)}{   V(x_0,r)^\nu}  \bigg( 1\wedge \frac{T_0}{\phi(x_0,r)} \bigg)^{-b} \left[ \sE(f,f) +  \frac{\lVert f \rVert_2^2}{\phi(x_0,r)}\right] \lVert f \rVert_1^{2\nu}.
	\end{align}
\end{defn}

\begin{lem}\label{l:Nash}
	For any $\nu>0$ and $b\ge 0$, we have
	\begin{align*}
		\text{\GFK \  $\;\; \Leftrightarrow \;\;$ \Nash}.
	\end{align*}
\end{lem}
\pf $(\Rightarrow)$ For this implication, we adopt the arguments of  \cite[Lemma 5.4]{GH14}.  Fix  $B:=B(x_0,r)$. Since  $ \sF \cap C_c(B)$ is dense in $\sF^{B}$ with respect to $\sE_1$-norm and $\sE(|f|,|f|) \le \sE(f,f)$ for all $f \in \sF$,  by linearity, 
it suffices to prove \eqref{e:GFK-Nash} for $0\le f \in \sF \cap C_c(B)$ with $\lVert f \rVert_1=0$.

Let $0\le f \in \sF \cap C_c(B)$ with  $\lVert f \rVert_1=1$.  For $a>0$, define $D_a:=\{y \in B:f(y)>a\}$.  Note that $(f-a)_+ \in \sF \cap C_c(D_a) \subset \sF^{D_a}$ and
\begin{align*}
	\int_{D_a} (f-a)_+^2 \,d\mu = \int_{B} (f-a)_+^2 \,d\mu \ge \int_{B} (f^2-2af)\,d\mu = \lVert f \rVert_2^2 - 2a.
\end{align*}
Further, by the Markov's inequality, $
\mu(D_a) \le 1/a$. Thus, using the Markovian property and \GFK, we obtain for all $a>0$,
\begin{align*}
	\sE(f,f)  &\ge \sE((f-a)_+,(f-a)_+) \nn\\	
	&\ge  \frac{c_1}{\phi(x_0,r)} \bigg[  \bigg( 1\wedge \frac{T_0}{\phi(x_0,r)} \bigg)^{b} \left( \frac{V(x_0,r)}{\mu(D_a)}\right)^\nu  -c_2 \bigg] \int_{D_a} (f-a)_+^2 \,d\mu \nn\\
	&\ge  \frac{c_1 ( \lVert f \rVert_2^2 - 2a )}{\phi(x_0,r)} \bigg[   \bigg( 1\wedge \frac{T_0}{\phi(x_0,r)} \bigg)^{b}   \left(aV(x_0,r)\right)^\nu    -c_2 \bigg] .
\end{align*}
By taking $a= \lVert f\rVert_{2}^2/4$, it follows that
\begin{align*}
	\sE(f,f) \ge  \frac{c_1  \lVert f \rVert_2^2 }{2\phi(x_0,r)} \bigg[   \bigg( 1\wedge \frac{T_0}{\phi(x_0,r)} \bigg)^{b}   \left(\frac{V(x_0,r)\lVert f\rVert_{2}^2}{ 4}\right)^\nu    -c_2 \bigg].
\end{align*}
This proves \eqref{e:GFK-Nash}.

$(\Leftarrow)$ 
Let $B:=B(x_0,r)$ and $D$ be a non-empty open subset of $B$. For any $f \in \sF^D$ with $\lVert f \rVert_2=1$, using \Nash \ and the Cauchy-Schwarz inequality, we obtain
\begin{align*}
	1=	\lVert f \rVert_{2}^{2+2\nu} &\le   \frac{c_1 \phi(x_0,r)}{   V(x_0,r)^\nu}  \bigg( 1\wedge \frac{T_0}{\phi(x_0,r)} \bigg)^{-b} \left[ \sE(f,f) +  \frac{1}{\phi(x_0,r)}\right] \bigg( \int_D |f| d\mu \bigg) ^{2\nu}\nn\\
	&\le \frac{c_1 \phi(x_0,r) \mu(D)^\nu}{   V(x_0,r)^\nu}  \bigg( 1\wedge \frac{T_0}{\phi(x_0,r)} \bigg)^{-b}  \left[ \sE(f,f) +  \frac{1}{\phi(x_0,r)}\right],
\end{align*}
implying that
\begin{align*}
	\sE(f,f) \ge \frac{1}{c_1 \phi(x_0,r)} \bigg[  \bigg( 1\wedge \frac{T_0}{\phi(x_0,r)} \bigg)^{b}  \bigg(\frac{V(x_0,r)}{\mu(D)}\bigg)^\nu  - c_1\bigg].
\end{align*}
By the definition \eqref{e:def-lambda1} of $\lambda_1(D)$, this proves \GFK. \qed

\section{$L^2$-mean value inequality and resolvent estimates}\label{ss:5}

In this section, we establish the  $L^2$-mean value inequality for subharmonic functions  and the lower resolvent estimates \Gl \  for  $\kappa > 0$ (see Definition \ref{d:Gl} for the precise definition).  The proofs mainly  follow the arguments in  \cite{GHH24}, where  the $L^2$-mean value inequality and a lower estimate of mean exit time from a ball (which correspond to \Gl \ for $\kappa=0$) were proven using \FK \ in place of \wFK.  
Note that, for \Gl,  $\kappa>0$ is  essential to  our approach, as \Gl \ for $\kappa>0$ alone implies \S  \ (see Lemma \ref{l:G->S}).

We begin with the notions of subharmonic and superharmonic functions.

\begin{defn}
	Let $D\subset M$ be an open set. A function  $u \in \sF'$ is said to be \textit{subharmonic} (resp. \textit{superharmonic})  in $D$ if
	\begin{align*}
		\sE(u,\vp) \le 0 \;\; (\text{resp. $\sE(u,\vp)\ge 0$}) \quad \text{ for all $0\le \vp \in \sF^D$.}
	\end{align*}
\end{defn}

\begin{lem}\label{l:MVI}
	Suppose that \VD, \TJ, \wFK \ and \CS \ hold. Let $u\in \sF'_b$ be subharmonic in $B:=B(x_0,R)$ with $R\in (0,\phi^{-1}(x_0,\delta_2T_0))$. Then there exists $C>0$ independent of $u$ and $B$ such that the following holds for all $0<b_1<b_2$ and  $r \in (0,R)$:
	\begin{align*}
		a_2 \le  \frac{C\phi(x_0,R)}{(b_2-b_1)^{2\nu} V(x_0,R)^\nu }   \bigg( \sup_{z \in B} \frac{1}{\phi(z,R-r)} + \frac{\sT }{b_2-b_1} \bigg)\,a_1^{1+\nu},
	\end{align*}
	where $\sT:=\esssup_{x \in B(x_0, (r+R)/2)} \int_{B^c} u_+(y)J(x,dy)$,
	\begin{align*}
		a_1:=\int_{B} (u-b_1)_+^2 d\mu \quad \text{and} \quad a_2:=\int_{B(x_0,r)} (u-b_2)_+^2 d\mu.
	\end{align*}
\end{lem}
\pf   We follow the proof of \cite[Lemma 10.2]{GHH24}, with $\overline R$, $r_1$ and $r_2$ replaced by $\phi^{-1}(x_0,\delta_2T_0)$, $R$ and $r$, respectively.  Note  that the condition corresponding to \CS \ is referred to as condition (ABB) in \cite{GHH24}.  By following the proof of  \cite[Lemma 8.2]{GHH24}, we see that  condition (EP) in \cite{GHH24} is satisfied with  $\phi^{-1}(x_0,\delta_2T_0)$ in place of $\overline R$.

In the proof of \cite[Lemma 10.2]{GHH24},
condition \FK \ is only used to get \cite[(10.8)]{GHH24}. With \wFK, instead of  \cite[(10.8)]{GHH24}, we obtain that for any $\eps>0$,
\begin{align*}
	a_2 \le c_1 \phi(x_0,R) \bigg[ \bigg( \frac{V(x_0,R)}{\mu(E)+\eps}  \bigg)^\nu -c_2  \bigg]_+^{-1} \sE(\vp v, \vp v),
\end{align*}
where $\vp$ stands for the cutoff function $\phi$ in the proof of \cite[Lemma 10.2]{GHH24}.  Following the rest of  the proof of  \cite[Lemma 10.2]{GHH24}, we  deduce that
\begin{equation}\label{e:(10.4)}
	a_2\le c_3 \phi(x_0,R) \bigg[ \bigg( \frac{(b_2-b_1)^2V(x_0,R)}{a_1}  \bigg)^\nu -c_2  \bigg]_+^{-1} \bigg( \sup_{z \in B} \frac{c_4}{\phi(z,(R-r)/2)} + \frac{3\sT}{b_2-b_1} \bigg)a_1.
\end{equation}
If $(b_2-b_1)^2 V(x_0,R)  \ge (2c_2)^{1/\nu} a_1$, then by \eqref{e:(10.4)} and \eqref{e:phi-scale}, we obtain
\begin{align*}
	a_2\le  \frac{2^\nu c_3 \phi(x_0,R)  } {(b_2-b_1)^{2\nu }V(x_0,R)^\nu }  \bigg( \sup_{z \in B} \frac{c_5}{\phi(z,R-r)} + \frac{3\sT}{b_2-b_1} \bigg)a_1^{1+\nu}.
\end{align*}
If $(b_2-b_1)^2 V(x_0,R)  < (2c_2)^{1/\nu}a_1$, then since $a_2\le a_1$, we get
\begin{align*}
&	\frac{\phi(x_0,R)}{(b_2-b_1)^{2\nu} V(x_0,R)^\nu }   \bigg( \sup_{z \in B} \frac{1}{\phi(z,R-r)} + \frac{\sT }{b_2-b_1} \bigg)a_1^{1+\nu}\\
& \ge \frac{\phi(x_0,R)a_1^{1+\nu}}{(b_2-b_1)^{2\nu} V(x_0,R)^\nu \phi(x_0, R-r)} \ge   \frac{a_1^{\nu}a_2}{(b_2-b_1)^{2\nu} V(x_0,R)^\nu } \ge  \frac{a_2}{2c_2}.
\end{align*}The proof is complete. \qed

\begin{prop}\label{p:MVI}
{\rm \bf ($L^2$-mean value inequality)}	
 Suppose that \VD, \TJ, \wFK \ and \CS \ hold. Let $B:=B(x_0,R)$ with $R\in (0,\phi^{-1}(x_0,\delta_2T_0))$, and   $u\in \sF'_b$ be non-negative and subharmonic in $B$.  There exists $C>0$ independent of $B$ and $u$ such that for any $\eps>0$,
	\begin{align*}
		\esssup_{2^{-1}B} u \le C(1+\eps^{-1/(2\nu)}) \bigg( \frac{1}{V(x_0,R)} \int_{B} u^2 d\mu\bigg)^{1/2}+ \eps \lVert u_+ \rVert_{L^\infty( (2^{-1} B)^c)}.
	\end{align*}
\end{prop}
\pf  Note that condition \FK \ is used in the proof of \cite[Theorem 10.1]{GHH24} only through \cite[Lemma  10.2]{GHH24}. By using Lemma \ref{l:MVI}, one can follow the proofs of \cite[Theorem 10.1 and Corollary 10.3]{GHH24}  line by line to obtain the  result. \qed

Using Proposition \ref{p:MVI} and adapting  the arguments from  \cite[Corollary 10.3 and Lemma 11.2]{GHH24}, we obtain  the following lemma of growth for superharmonic functions. The details are omitted.
\begin{lem}\label{l:LG}
Suppose that \VD, \TJ, \wFK \ and \CS \ hold.  Then there exist constants $\delta,\eta\in (0,1)$ such that, for any ball $B:=B(x_0,r)$ with $r\in (0,\phi^{-1}(x_0,\delta_2T_0))$, and any $u \in \sF'_b$ that is superharmonic in $B$ and is non-negative in $M$, the following holds: If 
\begin{align*}
	\mu(B \cap \{ u<a\}) \le \delta \mu(B)
\end{align*}
for some $a>0$, then $\essinf_{2^{-1}B} u \ge \eta a$.
\end{lem}

For  $\lambda > 0$,  the \textit{$\lambda$-resolvent}  $G^D_\lambda$ of $(\sE,\sF^D)$ is defined by
\begin{align*}
	G^D_\lambda f:=\int_0^\infty e^{-\lambda t} P^D_t fdt
\end{align*}
for   $f\in L^p(D)$, where $p\in [1,\infty]$.  For any $f \in L^2(D)$, we have  (see  \cite[(1.3.7)]{FOT}):
\begin{align}\label{e:resolvent}
	\sE_\lambda (G^D_\lambda f ,\vp) = \la f,\vp \ra  \quad \text{for all $\vp \in \sF^D$}.
\end{align}

\begin{lem}\label{l:resolvent-superharmonic}
	Let $D\subset M$ be an open set and  $\lambda> 0$. Then $G^D_\lambda \1_D+a$ is superharmonic in $D$ for any $a\in \R$.
\end{lem} 
\pf  Note that  $\lambda G_\lambda^D \1_D \le \lambda \1_D \int_0^\infty e^{-\lambda t} dt =\1_D$. 
By \eqref{e:resolvent}, we obtain for any $0\le \vp \in \sF^D$,
\begin{align*}
	\sE(G^D_\lambda \1_D+a,\vp) = 	\sE(G^D_\lambda \1_D,\vp) = \la \1_D, \vp \ra - \lambda\la  G^D_\lambda \1_D, \vp\ra  = \la \1_D - \lambda G^D_\lambda \1_D, \vp \ra \ge 0.
\end{align*}
\qed 

We  introduce   condition  \Gl, referred to as the \textit{lower resolvent estimates}.
\begin{defn}\label{d:Gl}
	We say that  \Gl \ holds  if there  exist constants $\kappa> 0$ and  $c_1>0$  such that for any ball $B:=B(x_0,r)$ with $r\in (0,\phi^{-1}(x_0,T_0))$,
	\begin{equation*}
		\essinf_{ B(x_0,r/4)} G^{B}_{\kappa/\phi(x_0,r)} \1_{B} \ge c_1\phi(x_0,r).
	\end{equation*}
\end{defn}

\begin{lem}\label{l:G->S}
	For every $\kappa> 0$, \Gl \ implies \S.	
\end{lem}
\pf   Let $B:=B(x_0,r)$ with $r\in (0,\phi^{-1}(x_0,T_0))$ and   $\lambda:=\kappa/\phi(x_0,r)$.  By the semigroup property, we have for all $t>0$ and $\mu$-a.e. $y \in 4^{-1}B$,
\begin{align}\label{e:G->S}
	G^{B}_\lambda \1_{B}(y)&=  \int_0^t e^{-\lambda s}P^{B}_s \1_{B}(y) ds  + e^{-\lambda t} P^{B}_t G^{B}_\lambda \1_{B}(y)  \le   t + P_t^{B}\1_{B} (y) \esssup_{ B} G^{B}_\lambda \1_{B}. \end{align}
By \Gl, there exists $c_1>0$ such that  $ \essinf_{4^{-1}B}G^{B}_\lambda \1_{B} \ge c_1^{-1}\phi(x_0,r)$. Further, by the  $L^\infty$-contraction property of the semigroup $(P^B_t)_{t\ge 0}$, we have
$\esssup_B G_\lambda^{B} \1_B \le   \int_0^\infty e^{-\lambda t} dt =\kappa^{-1} \phi(x_0,r)$. Thus, we deduce from  \eqref{e:G->S} that   for all $t\le \phi(x_0,r)/(2c_1)$,
\begin{align*}
	\essinf_{4^{-1} {B}} P_{t}^{B}\1_{B}  \ge \frac{ \essinf_{ 4^{-1} B}G^{B}_\lambda \1_{B} -t}{\esssup_{B} G^{B}_\lambda \1_{B} }  \ge  \frac{ c_1^{-1}\phi(x_0,r) -t}{\kappa^{-1}\phi(x_0,r) }  \ge \frac{\kappa}{2c_1}.
\end{align*}  \qed

\begin{lem}\label{l:capacity}
	Suppose that \TJ \ and \CS \ hold. Then there exists $C>0$ such that for any ball $B:=B(x_0,r)$ with $r\in (0,\phi^{-1}(x_0,T_0))$, there exists a cutoff function $\vp$ for $2^{-1}B\Subset B$ so that
	\begin{align*}
	\sE(\vp,\vp) \le \frac{CV(x_0,r)}{\phi(x_0,r)}.
	\end{align*}
\end{lem}
\pf Fix  $f \in \sF$  such that $f=1$ in $B$. Then $d\Gamma^{(L)}(f,f)=0$ in $B$ and $f(x)-f(y)=0$ for all $x,y \in B$. By \CS, there exists   a cutoff function $\vp$ for $\frac12 B\Subset \frac34 B$ such that
\begin{align}\label{e:capacity-1}
	&\int_{\frac34 B} d\Gamma^{(L)}(\vp,\vp)  + \int_{B\times B}  (\vp(x)-\vp(y))^2 J(dx,dy) \nn\\
	&= \int_{\frac34 B} f^2 d\Gamma^{(L)}(\vp,\vp)  + \int_{B\times B} f(x)^2 (\vp(x)-\vp(y))^2 J(dx,dy)\nn\\
	&\le   \sup_{z\in B}\frac{c_1}{\phi(z,r/4)} \int_B f^2 d\mu  =  \sup_{z\in B}\frac{c_1V(x_0,r)}{\phi(z,r/4)}.
\end{align}
On the other hand, using \TJb,  we obtain
\begin{align*}
	&	\int_{B^c\times M}  (\vp(x)-\vp(y))^2 J(dx,dy)=	\int_{\frac34B} \vp(y)^2 \int_{B^c}   J(y,dx)\,\mu(dy)  \\
	&\le V(x_0, 3r/4)\sup_{z\in B}  J(z,B(z, r/4)^c)  \le \sup_{z\in B}\frac{c_2V(x_0,r)}{\phi(z, r/4)}.
\end{align*}
Combining this with \eqref{e:capacity-1} and using \eqref{e:phi-scale2}, we conclude that
\begin{align*}
	 \sE(\vp,\vp)	& \le \sup_{z\in B}\frac{(c_1+2c_2)V(x_0,r)}{\phi(z, r/4)} \le \frac{c_3V(x_0,r)}{\phi(x_0,r)}.
\end{align*} \qed

We now establish \Gl. The proof of the next proposition follows from  \cite[Lemma 12.4]{GHH24}, with some modifications.
\begin{prop}\label{p:G}
	Suppose that \VD, \TJ, \wFK \ and \CS \ hold. Then \Gl \   holds for all $\kappa> 0$.
\end{prop}
\pf  Let $\kappa > 0$, $B:=B(x_0,r)$ be a ball with $r\in (0,\phi^{-1}(x_0,T_0))$ and $u:=G^B_{\kappa/\phi(x_0,r)} \1_B$. We consider the following two cases separately.

\smallskip

Case 1: $r<\phi^{-1}(x_0,\delta_2T_0)$. Fix $\eps>0$ and let $u_\eps:=u+\eps$.  Let  $\delta,\eta\in (0,1)$ be the constants in Lemma \ref{l:LG}. Choose $a_0>0$ such that 
$$	a_0 \int_{2^{-1}B} \frac{d\mu}{u_\eps} d\mu= \delta V(x_0,r/2).$$
By Markov's inequality,  we have $	\mu( 2^{-1} B \cap \{ u_\eps<a_0\} )  \le a_0 \int_{2^{-1}B} u_\eps^{-1} d\mu =   \delta V(x_0,r/2)$. Note that, by Lemma \ref{l:resolvent-superharmonic}, $u_\eps$ is superharmonic in $B$. Applying Lemma \ref{l:LG}, we obtain
\begin{align}\label{e:Gl-2}
	\essinf_{4^{-1} B} u_\eps \ge \eta a_0 = \eta \delta V(x_0,r/2)\bigg(\int_{2^{-1}B} \frac{d\mu}{u_\eps} \bigg)^{-1} .
\end{align}
By Lemma \ref{l:capacity}, there exists a cutoff function $\vp$ for $2^{-1}B \Subset B$ such that  $\sE(\vp,\vp)\le \frac{c_1V(x_0,r)}{\phi(x_0,r)}$. Hence, using \cite[Proposition 12.3]{GHH24} in the first inequality below, we get 
\begin{align}\label{e:Gl-3}
	\sE(u_\eps, \frac{\vp^2}{u_\eps})  = \sE(u, \frac{\vp^2}{u_\eps}) \le 	3\sE(\vp,\vp) \le \frac{3c_1V(x_0,r)}{\phi(x_0,r)}.
\end{align}Further, by \eqref{e:resolvent}, since $\vp$ is a cutoff function for $2^{-1}B\Subset B$, we have
\begin{align}\label{e:Gl-4}
	\sE(u_\eps, \frac{\vp^2}{u_\eps}) &= \la \1_B , \frac{\vp^2}{u_\eps} \ra - \frac{\kappa}{\phi(x_0,r)} \la u,  \frac{\vp^2}{u_\eps} \ra \ge \int_{2^{-1}B} \frac{1}{u_\eps} d\mu -\frac{\kappa}{\phi(x_0,r)} \int_B \frac{u}{u_\eps} d\mu.
\end{align}
Combining \eqref{e:Gl-2}, \eqref{e:Gl-3} and \eqref{e:Gl-4}, and using \VD, we obtain
\begin{align*}
	\essinf_{4^{-1}B} u_\eps \ge c_2 \eta \delta V(x_0,r) \bigg(\frac{3c_1V(x_0,r)}{\phi(x_0,r)} + \frac{\kappa}{\phi(x_0,r)} \int_B \frac{u}{u_\eps} d\mu \bigg)^{-1}.
\end{align*}
Taking the limit as $\eps \to 0$, it follows that 
\begin{align*}
	\essinf_{4^{-1} B} u \ge  c_2 \eta \delta V(x_0,r) \bigg(\frac{3c_1V(x_0,r)}{\phi(x_0,r)} + \frac{\kappa}{\phi(x_0,r)} \int_B d\mu \bigg)^{-1}=\frac{c_2 \eta \delta \phi(x_0,r)}{3c_1+\kappa},
\end{align*}
proving that \Gl \ holds in this case.

Case 2: $r\ge \phi^{-1}(x_0,\delta_2T_0)$. Set $r_0:=2^{-1}\phi^{-1}(x_0,\delta_2T_0)$ and let $\{B(z_i,r_0/4)\}_{i\ge 1}$ be a covering of $4^{-1}B$,  where $z_i \in \overline{4^{-1}B}$ for all $i \ge 1$. Note that $B(z_i, r_0)\subset B(z_i,r/2)\subset B$ for all $i\ge 1$. Moreover, 
by \eqref{e:phi-scale2}, there exists $c_3\ge 1$ such that 
\begin{align}\label{e:Gl-5}
c_3	\phi(x_0,r) \ge  \phi(z_i,r) \ge \phi(z_i,r_0) \ge c_3^{-1} \phi(x_0,r) \quad \text{for all $i\ge 1$}.
\end{align}
Using the result from Case 1 with $\kappa$ replaced by $c_3\kappa$ and  \eqref{e:Gl-5}, we obtain
\begin{align*}
		\essinf_{4^{-1}B} u  &\ge  \inf_{i\ge 1}\essinf_{B(z_i,r_0/4)} G^{B(z_i,r_0)}_{c_3\kappa/\phi(z_i,r_0)} \1_{B(z_i,r_0)} \ge c_4 \inf_{i\ge 1} \phi(z_i, r_0) \ge c_3^{-1}c_4 \phi(x_0,r).
\end{align*}The proof is complete. \qed

\section{The case of special scale function}\label{ss:special}

 In this section, we establish on-diagonal upper estimates and tail estimates of the semigroup for  the special case where the scale function $\phi$ is given by $\phi(x,r)=r^\beta$ for some $\beta>0$, and the localizing constant is independent of the spatial variable $x$. To obtain these results, we  consider  truncated Dirichlet forms and their relationships with the original Dirichlet form  $(\sE,\sF)$, which are presented in
  Subsection \ref{ss:truncation}. The assumption $\phi(x,r)=r^\beta$  simplifies the arguments for the truncated Dirichlet forms. 
 We will see in Section \ref{ss:7} that, by applying a metric transformation, this special case leads to our main results for a general situation.
 
 \smallskip

Let $\beta>0$ and $R_0\in (0,\infty]$ be fixed.  

\begin{defn}
	We say that  \TJb \ (resp. \dTJb \ or \IVJb) holds if \TJ \ (resp. \dTJ \ or \IVJ)  holds with  $\phi(x,r)=r^\beta$.
\end{defn}

\begin{defn}
	We say that  \Nashb \ holds if there exist constants $\nu,C>0$ and $b'\ge 0$ such that for any ball $B:=B(x_0,r)$ with   $r>0$ and $f \in \sF^{B}$,
	\begin{align*}
			\lVert f \rVert_{2}^{2+2\nu} \le   \frac{C r^\beta}{   V(x_0,r)^\nu}  \bigg( 1\wedge \frac{R_0}{r} \bigg)^{-b'} \left[ \sE(f,f) +  \frac{\lVert f \rVert_2^2}{r^\beta}\right] \lVert f \rVert_1^{2\nu}.
	\end{align*}
\end{defn}

\begin{defn}
(i)	We say that    \Sb \ holds if  there exist constants $\eps_0,a_0\in (0,1)$ such that for any ball $B:=B(x_0,r)$ with $r\in (0,R_0)$,
	\begin{align*}
		\essinf_{ \frac14 B}	P^B_t \1_{B} \ge \eps_0 \quad \text{for all $t \in (0, a_0 r^\beta]$}.
	\end{align*}
	
	\noindent (ii) 	We say that    \TEb \ holds if  there exists  $C>0$ such that for any ball $B:=B(x_0,r)$ with $r>0$,
	\begin{align}\label{e:TEb}
		\esssup_{ \frac14 B}	P_t \1_{B^c}  \le \frac{Ct}{(r\wedge R_0)^\beta} \quad \text{for all $t>0$}.
	\end{align}
\end{defn}
$$
\textbf{In the remainder of this section, we assume that \eqref{e:J-density} is satisfied.}
$$

Define  for $\rho>0$ and  $x\in M$,
\begin{align*}
	J_1^{(\rho)}(x,dy) :=  \1_{\{d(x,y) < \rho\}} J(x,dy) \quad \text{ and } \quad J_2^{(\rho)}(x,dy) = J(x,dy)- J_1^{(\rho)}(x,dy).
\end{align*} 
The \textit{$\rho$-truncated  form} $(\sE^{(\rho)},\sF)$ of $(\sE,\sF)$ on $L^2(M)$ is defined by
\begin{align*}
	\sE^{(\rho)}(f,g) = \sE^{(L)}(f,g) + \int_{M\times M} (f(x)-f(y))(g(x)-g(y))  \,J_1^{(\rho)}(x,dy)\mu(dx).
\end{align*}
The form $(\sE^{(\rho)},\sF)$  is  a symmetric  Dirichlet form on $L^2(M)$. While   $(\sE^{(\rho)},\sF)$ may not be   regular in general,  it is regular under \TJb.

\begin{lem}\label{l:comparison-L2}
	Suppose that \TJb \ holds. There exists $C>0$  such that  for any  $\rho>0$,
	\begin{align*}
		\sE(f,f) -	\sE^{(\rho)}(f,f)\le   C\rho^{-\beta}\lVert f \rVert_2^2 \quad \text{for all $f \in \sF$}.
	\end{align*}
	Consequently, the norms  $\sE_1^{1/2}$ and $(\sE^{(\rho)}_1)^{1/2}$ are equivalent on $\sF$, and  $(\sE^{(\rho)}, \sF)$ is  regular.
\end{lem}
\pf  The result follows immediately from  Lemma \ref{l:truncated-general}. \qed 

A countable union of exceptional sets is  an exceptional set. By  Lemma \ref{l:comparison-L2}, whenever \TJb \ is assumed,  we will, without loss of generality, assume that any $f \in \sF$ is  $\sE^{(\rho)}$-quasi-continuous for all $\rho\in \Q_+$.

\subsection{Tail estimates for truncated  Dirichlet form}\label{ss:Tail-truncated}
 For an open subset $D$ of $M$, denote by  $(P^{(\rho),D}_t)_{t\ge 0}$ the semigroup associated with the part of $(\sE^{(\rho)},\sF)$ on $D$. Let $P^{(\rho)}_t:=P^{(\rho),M}_t$.

The following lemma  is motivated by \cite[Theorem 3.1]{GHL14}.

\begin{lem}\label{l:truncated-off-diagonal}
	Let $\rho>0$.	Suppose that $(\sE^{(\rho)},\sF)$ is regular and  there exist   constants $\lambda,T>0$ such that 
	for  any  $z\in M$,
	\begin{align}\label{e:truncated-off-diagonal-ass}
	\essinf_{B(z,\rho/12)}  P^{(\rho),B(z,\rho)}_{t} \1_{B(z,\rho)}   \ge 1-\lambda \quad \text{for all $t \in (0,T]$}.
	\end{align}
	Then for any $x_0 \in M$ and   $r\ge 3\rho$,  we have
	\begin{align}\label{e:truncated-off-diagonal}
		\esssup_{ B(x_0,r/4)}  P^{(\rho)}_{t} \1_{B(x_0,r)^c}\le \lambda^{r/(3\rho)-1} \quad \text{for all $t\in (0, T]$}.
	\end{align}
\end{lem}
\pf If $\lambda\ge 1$, then \eqref{e:truncated-off-diagonal} clearly holds. Assume $\lambda<1$.  Define for $k\ge 0$,
\begin{align*}
	u_k(s,x,y):= P^{(\rho)}_{s} \1_{B(x,9k\rho/4)^c}(y), \quad s>0, \, x,y \in M.
\end{align*}
Let  $N:=\max\{n\in \N: n\rho \le r/3\}$ and $\{B(z_i, \rho/12)\}_{i\ge 1}$ be a covering of $B(x_0,r/4)$, where $z_i \in \overline{B(x_0,r/4)}$ for all $i\ge 1$. For all $t>0$, we have
\begin{align*}
	\esssup_{y \in B(x_0,r/4)}  P^{(\rho)}_{t} \1_{B(x_0,r)^c}(y) &\le \sup_{i\ge 1}  \esssup_{y \in B(z_i,\rho/12)}   P^{(\rho)}_t \1_{B(z_i,3r/4)^c}(y)\le  \sup_{i\ge 1}  \esssup_{y \in B(z_i,\rho/12)}  u_N(t,z_i,y).
\end{align*}
Applying Proposition \ref{p:comparison-entrance} and using \eqref{e:truncated-off-diagonal-ass}, we get that for all $i\ge 1$ and  $t\in (0, T]$,
\begin{align*}
	\esssup_{y \in B(z_i,\rho/12)}	u_N(t,z_i,y) &	 \le  	\esssup_{y \in B(z_i,\rho/12)} \big( 1-  P^{(\rho),B(z_i,\rho)}_t \1_{B(z_i,\rho)}(y) \big) \sup_{0<s\le t} \esssup_{w \in B(z_i,9\rho/4 )}  u_N(s,z_i,w)  \nn\\
	&\le \lambda \sup_{0<s \le T} \esssup_{y \in B(z_i,9\rho/4)}  u_N(s,z_i,y) .
\end{align*}
Fix $i\ge 1$. Let $\{B(z_{ij}, \rho/12)\}_{j\ge 1}$ be a  covering of $B(z_i,9\rho/4)$, where $z_{ij} \in \overline{B(z_i,9\rho/4)}$ for all $j\ge 1$. Applying  Proposition \ref{p:comparison-entrance}  and using \eqref{e:truncated-off-diagonal-ass}, we obtain
\begin{align*}
	& \sup_{0<s\le T} \esssup_{y \in B(z_i,9\rho/4 )}  u_N(s,z_i,y)   \le \sup_{0<s< T} \sup_{j \ge 1}
	\esssup_{y\in B(z_{ij},\rho/12)} u_{N-1}(s,z_{ij},y) \\
	& \le \sup_{0<s\le  T} \sup_{j \ge 1}
	\esssup_{y\in B(z_{ij},\rho/12)} \big( 1-  P^{(\rho),B(z_{ij},\rho)}_s \1_{B(z_i,\rho)}(y) \big) \sup_{0<s'\le s} \esssup_{y' \in B(z_{ij},9\rho/4 )}  u_{N-1}(s',z_{ij},y')\\
	& \le   \lambda  \sup_{0<s\le  T} \sup_{j\ge 1} \esssup_{y\in B(z_{ij},9\rho /4)}  u_{N-1}(s,z_{ij},y).
\end{align*}
By iterating this procedure, we can deduce the existence of sequences $\{z_{i_1i_2 \cdots i_k}\}_{i_1,\cdots, i_k\ge 1}$ for $1\le k\le N$ that satisfy the following  properties  for all $2\le k\le N$ and $i_1,\cdots, i_{k-1} \ge 1$:

\setlength{\leftskip}{3mm}

\smallskip

\noindent (1)   $B(z_{i_1\cdots i_{k-1}},9\rho/4)\subset \cup_{i_k\ge 1}B(z_{i_1\cdots i_k},\rho/12)$ and  $z_{i_1 \cdots i_k} \in \overline{B(z_{i_1 \cdots i_{k-1}},9\rho/4)}$ for all $i_k \ge 1$.  

\noindent (2) It holds that
\begin{align*}
	&	\sup_{0<s \le  T}  \esssup_{y \in B(z_{i_1\cdots i_{k-1}},9\rho/4)} u_{N-k+1}(s,z_{i_1\cdots i_{k-1}},y) \le\lambda\sup_{0<s\le  T}  \sup_{i_k\ge 1} \esssup_{y\in B(z_{i_1\cdots i_{k}},9\rho/4)} u_{N-k}(s,z_{i_1\cdots i_{k}},y) .
\end{align*}

\setlength{\leftskip}{0mm}

\noindent It follows that for all $t \in (0,T]$,
\begin{align*}
	&	\esssup_{y \in B(x_0,r/4)}  P^{(\rho)}_{t} \1_{B(x_0,r)^c}(y) \le \lambda \sup_{0<s< T} \sup_{i_1\ge 1}  \esssup_{y \in B(z_{i_1},9\rho/4 )}  u_N(s,z_{i_1},y)\\
	&\le \cdots \le \lambda^N \sup_{0<s < T} \sup_{i_1, \cdots, i_{N}\ge 1}   \esssup_{y \in B(z_{i_1\cdots i_{N}},9\rho/4) } u_{1}(s,z_{i_1\cdots i_{N}},y) \le \lambda^N.
\end{align*}
Since $\lambda<1$ and $N+1>r/(3\rho)$,  this completes the proof.  \qed

\begin{lem}\label{l:Markov}
	Let $\rho>0$. Suppose that $(\sE^{(\rho)},\sF)$ is  conservative and regular.
	If  there exist  constants $k,\lambda,T>0$ and a  function $f:[k\rho,\infty) \to (0,\infty)$ such that
	 for any $r\ge k\rho$ and $t\in (0,T]$,
	\begin{align}\label{e:Markov-assumption}
	\sup_{z\in M}	\esssup_{ B(z,r/4)} P^{(\rho)}_{t} \1_{B(z,r)^c} \le f(r),
	\end{align}
	then for any $x_0\in M$, $r\ge 3k\rho$ and $t\in (0,T]$,
		\begin{align*}
		\essinf_{ B(x_0,r/12)} P^{(\rho), B(x_0,r)}_{t} \1_{B(x_0,r)} \ge  1-2f(r/3).
	\end{align*}
\end{lem} 
\pf Let $B:=B(x_0,r)$  with $r\ge 3k\rho$. By Proposition \ref{p:comparison-entrance}, we obtain for all $t \in (0,T]$,
\begin{equation}\label{e:Markov-1}
		\essinf_{12^{-1}B}  P^{(\rho), B}_{t} \1_{B}\ge 	\essinf_{12^{-1}B} P^{(\rho), B}_{t} \1_{3^{-1} B}  \ge 	\essinf_{12^{-1}B} P^{(\rho)}_{t} \1_{3^{-1} B} -  \sup_{0<s\le t} \esssup_{ B(x_0,r+2\rho) \setminus \overline{(2/3) B} } P^{(\rho)}_{s} \1_{3^{-1} B}.
\end{equation}
Since $(\sE^{(\rho)},\sF)$ is conservative,  by \eqref{e:Markov-assumption}, we have for all $t\in (0,T]$, 
\begin{align}\label{e:Markov-2}
	\essinf_{12^{-1}B} P^{(\rho)}_{t} \1_{3^{-1} B}  &= 1- 
	\esssup_{12^{-1}B} P^{(\rho)}_{t} \1_{(3^{-1}B)^c}  \ge 1 - f(r/3).
\end{align}
Let $\{B(z_i, r/12)\}_{i\ge 1}$ be a covering of $B(x_0,r+2\rho) \setminus \overline{(2/3)B} $, where $z_i \in \overline{B(x_0,r+2\rho)} \setminus (2/3)B $ for all $i\ge 1$.   By \eqref{e:Markov-assumption}, we obtain for all $t\in (0,T]$,
\begin{align}\label{e:Markov-3}
	\sup_{0<s\le t} \esssup_{ B(x_0,r+2\rho) \setminus \overline{(2/3) B} } P^{(\rho)}_{s} \1_{3^{-1} B} &\le 	\sup_{0<s\le t} \sup_{i\ge 1} \esssup_{ B(z_i,r/12)} P^{(\rho)}_{s} \1_{B(z_i,R/3)^c}\le f(r/3).
\end{align}
Combining \eqref{e:Markov-1}, \eqref{e:Markov-2} and \eqref{e:Markov-3}, we arrive at the result. \qed

\begin{lem}\label{l:mE-truncated-consequence}
	Suppose that  \TJb \ and \Sb \ hold. There exists  a constant $a_1\in (0,1)$  such that 
	for  any  $\rho>0$, ball $B:=B(x_0,r)$ with $r \in (0,R_0)$ and $t\in (0,  a_1(\rho \wedge r)^\beta]$,
	\begin{align*}
		\essinf_{ 4^{-1} B}  P^{(\rho),B}_{t} \1_{B}\ge \frac{\eps_0}{2},
	\end{align*}
	where $\eps_0  \in (0,1)$ is the constant in \Sb.
\end{lem}
\pf  Using  \eqref{e:comparison-tail-2}, \TJb \ and \Sb, we obtain for all $t\in (0,a_0r^\beta]$,
\begin{align*}
		\essinf_{ 4^{-1} B} P^{(\rho),B}_{t} \1_{B} \ge 	\essinf_{ 4^{-1}B} P^B_{t} \1_{B} - 2t \sup_{z\in M} J_2^{(\rho)}(z,M)  \ge  \essinf_{ 4^{-1} B} P^{B}_{t} \1_{B}  - \frac{c_1t}{\rho^\beta} \ge \eps_0 - \frac{c_1t}{\rho^\beta},
\end{align*}
where $a_0$ is the  constant in \Sb. By taking $a_1:=a_0 \wedge (\eps_0/(2c_1))$,   we get the  result.\qed

As a consequence of Lemma \ref{l:mE-truncated-consequence}, we get the following conservativeness result.

\begin{lem}\label{l:conservative}
	Suppose that \Sb \ holds. Then $(\sE,\sF)$  is conservative. Moreover, if \TJb \ further holds, then $(\sE^{(\rho)},\sF)$ is conservative for all $\rho>0$.
\end{lem}
\pf  By Lemma \ref{l:mE-truncated-consequence}, the results follow from \cite[Lemma 4.6]{GHH17}. \qed

Recall that if \TJb \ holds, then  for any $\rho>0$, $(\sE^{(\rho)},\sF)$ is a regular Dirichlet on $L^2(M)$. 
Using Lemma \ref{l:mE-truncated-consequence}, we  apply Lemma \ref{l:truncated-off-diagonal} with  $\lambda= 1-\eps_0/2$ and   $T=a_1\rho^\beta$ to obtain the following result.
  \begin{lem}\label{l:truncated-off-diagonal-1}
  	Suppose that  \TJb \  and \Sb \ hold.	Then  for any $\rho \in (0,R_0)$, $x_0\in M$  and  $r\ge 3\rho$, 
  	\begin{align*}
  		\esssup_{ B(x_0,r/4)}  P^{(\rho)}_{t} \1_{B(x_0,r)^c}\le  \big(1-\frac{\eps_0}{2}\big)^{r/(3\rho)-1} \quad \text{for all $t\in (0, a_1\rho^\beta]$},
  	\end{align*}
where  	$a_1\in (0,1)$ is the constant in Lemma \ref{l:mE-truncated-consequence}.
  \end{lem}

  By Lemmas \ref{l:Markov} and   \ref{l:conservative},  the next lemma follows  from Lemma \ref{l:truncated-off-diagonal-1}.
  
   \begin{lem}\label{l:truncated-off-diagonal-2}
  	Suppose that  \TJb \  and \Sb \ hold.	Then  for any $\rho \in (0,R_0)$, $x_0\in M$  and  $r\ge 9\rho$,
  	\begin{align*}
  		\essinf_{ B(x_0,r/12)} P^{(\rho), B(x_0,r)}_{t} \1_{B(x_0,r)} \ge  1-2\big(1-\frac{\eps_0}{2}\big)^{r/(9\rho)-1} \quad \text{for all $t\in (0, a_1\rho^\beta]$},
  	\end{align*}
  	where  	$a_1\in (0,1)$ is the constant in Lemma \ref{l:mE-truncated-consequence}.
  \end{lem}

 	 \begin{prop}\label{p:truncated-off-diagonal-3-t}
 		Suppose that  \TJb \  and \Sb \ hold.	There exists  $C\ge 1$  such that  for any $\rho \in \Q_+ \cap (0,R_0)$, $x_0\in M$  and  $r\ge 9\rho$,
 		\begin{align}\label{e:truncated-off-diagonal-3-t}
 			\esssup_{B(x_0,r/4)} P^{(\rho)}_{t} \1_{B(x_0,r)^c}\le  \bigg( \frac{Ct}{\rho^\beta}\bigg)^{r/(5\rho)} \quad \text{for all $t>0$}.
 		\end{align}
 	\end{prop}
 	\pf Fix $\rho \in \Q_+ \cap (0,R_0)$ and set $\eps:=10^{-\beta}a_1$, 
 	where  	$a_1\in (0,1)$ is the constant in Lemma \ref{l:mE-truncated-consequence}. If $t\ge \eps\rho^\beta$, then by taking $C$ larger than $1/\eps$, we get \eqref{e:truncated-off-diagonal-3-t}.

 	Suppose $t<\eps\rho^\beta$. Let $\rho_t\in \Q_+$ be such that  
 	$$(t/a_1)^{1/(10\beta)}(\rho/10)^{9/10} < \rho_t \le (t/a_1)^{1/(10\beta)}(\rho/9)^{9/10}.$$ Then  we have $\rho \ge 9\rho_t$ and $t<a_1\rho_t^\beta$. 	For any $z \in M$ and $s \in (0,t]$, using \eqref{e:comparison-tail-1} and \TJb \ in the first inequality below,  Lemma \ref{l:truncated-off-diagonal-2} in the second, and the fact that $\sup_{a>0}a^{9\beta}e^{-a} <\infty$  in the third, we obtain
 	\begin{align*}
 	&\essinf_{B(z,\rho/12)}	P_s^{(\rho), B(z,\rho)}\1_{B(z,\rho)} \ge \essinf_{ B(z,\rho/12)}	P_s^{(\rho_t), B(z,\rho)}\1_{B(z,\rho)}  - \frac{c_1s}{\rho_t^\beta} \\
 	&\ge  1-2e^{-c_2\rho/\rho_t} - \frac{c_1t}{\rho_t^\beta} \ge 1- c_3\bigg(\frac{\rho_t}{c_2\rho}\bigg)^{9\beta} - \frac{c_1t}{\rho_t^\beta}  \ge 1- \frac{c_4 t^{9/10}}{\rho^{9\beta/10}}.
 	\end{align*}
 	 By Lemma \ref{l:truncated-off-diagonal}, it follows that for any $x_0 \in M$ and  $r\ge 9\rho$,
 		\begin{align*}
 		\esssup_{B(x_0,r/4)} P^{(\rho)}_{t} \1_{B(x_0,r)^c}\le \bigg(  \frac{c_4 t^{9/10}}{\rho^{9\beta/10}}\bigg)^{r/(3\rho)-1}  \le \bigg(  \frac{c_4 t^{9/10}}{\rho^{9\beta/10}}\bigg)^{2r/(9\rho)}  =  \bigg(  \frac{c_4^{10/9}t}{\rho^\beta}\bigg)^{r/(5\rho)}. 
 	\end{align*}
 	The proof is complete. \qed

Using Lemmas \ref{l:Markov} and \ref{l:conservative}, we deduce the next result  from Proposition \ref{p:truncated-off-diagonal-3-t}.
 \begin{prop}\label{p:truncated-off-diagonal-4-t}
 	Suppose that  \TJb \  and \Sb \ hold.	There exists   $C\ge 1$  such that  for any $\rho \in \Q_+\cap (0,R_0)$, $x_0\in M$  and  $r\ge 27\rho$, 
 	\begin{align*}
 		\essinf_{B(x_0,r/12)} P^{(\rho), B(x_0,r)}_{t} \1_{B(x_0,r)}  \ge  1-2 \bigg( \frac{Ct}{\rho^\beta}\bigg)^{r/(15\rho)} \quad \text{for all $t>0$}.
 	\end{align*}
 \end{prop}

 \subsection{On-diagonal upper estimate  for truncated Dirichlet form}
 
 In this subsection,  we   establish an on-upper estimate for the heat kernel of  $(P^{(\rho)}_t)_{t\ge 0}$.  The main result of this subsection is Theorem \ref{t:NDU-truncated-on}.

 \begin{lem}\label{l:Nash-truncated-consequence}
 	 	 (i) Suppose that \Nashb  \ holds. Then  there exists $C>0$ such that for any  ball $B:=B(x_0,r)$ with $r>0$,
 	\begin{align*}
 		\lVert P^{B}_t  \rVert_{L^1(B) \to L^\infty(B)} \le \frac{Cr^{\beta/\nu} }{  V(x_0,r) (t \wedge r^\beta)^{1/\nu}}  \bigg( 1\wedge \frac{R_0}{r} \bigg)^{-b'/\nu}  \quad \text{for all $t>0$}.
 	\end{align*}

 	\noindent (ii) 
 	Suppose that \TJb \ and \Nashb  \ hold.  Then there exists  $C>0$ such that for any $\rho>0$ and any ball $B:=B(x_0,r)$ with $r>0$,
 	\begin{align}\label{e:Nash-truncated-consequence-1}
 		\lVert  P^{(\rho),B}_t  \rVert_{L^1(B) \to L^\infty(B)} \le \frac{Cr^{\beta/\nu} }{  V(x_0,r) (t\wedge \rho^\beta \wedge r^\beta )^{1/\nu}}  \bigg( 1\wedge \frac{R_0}{r} \bigg)^{-b'/\nu}   \quad \text{for all $t>0$}.
 	\end{align}
 \end{lem}
 \pf  Since the proofs are similar, we only give the proof for (ii).

  Let $B:=B(x_0,r)$.  By  the semigroup property and duality,  for all $t>0$,
 \begin{equation}\label{e:Nash-UC-1}
 	\lVert P^{(\rho),B}_t  \rVert_{L^1(B) \to L^\infty(B)}   \le 	\lVert P^{(\rho),B}_{t/2}  \rVert_{L^1(B) \to L^2(B)}	\lVert P^{(\rho),B}_{t/2}  \rVert_{L^2(B) \to L^\infty(B)} = 	\lVert P^{(\rho),B}_{t/2}  \rVert_{L^1(B) \to L^2(B)}^2.
 \end{equation} 
 Let $\sF^{(\rho),B}$ denote the $(\sE^{(\rho)}_1)^{1/2}$-closure of $\sF \cap C_c(B)$ in $\sF$. By Lemma \ref{l:comparison-L2},  $\sF^{(\rho),B}=\sF^B$. Further, 
 by \Nashb, \TJ \ and Proposition \ref{p:comparison-tail}, we obtain  $	\lVert f \rVert_{2}^{2+2\nu} \le   a \left(  \sE^{(\rho)}(f,f) +  \delta\lVert f \rVert_2^2\right) \lVert f \rVert_1^{2\nu}$ for all $f\in \sF^{(\rho),B}=\sF^B$, where 
 $$
 a:=\frac{c_1 r^\beta}{   V(x_0,r)^\nu}  \bigg( 1\wedge \frac{R_0}{r} \bigg)^{-b'} \quad \text{and} \quad \delta:=\rho^{-\beta} + r^{-\beta}.
 $$
 By the  arguments in \cite[Theorem 2.1 and line 8 on page 252]{CKS87}, this implies that for all $t>0$,
 \begin{align*} \lVert P^{(\rho),B}_{t/2} \rVert_{L^1(B) \to L^2(B)}^2 \le  \frac{c_1^{1/\nu} r^{\beta/\nu}}{  \nu^{1/\nu} V(x_0,r) t^{1/\nu}}  \bigg( 1\wedge \frac{R_0}{r} \bigg)^{-b'/\nu} \exp \left( \frac{t}{\rho^\beta} + \frac{t}{r^\beta}\right).
 \end{align*}
 Combining this with \eqref{e:Nash-UC-1}, we deduce that \eqref{e:Nash-truncated-consequence-1} holds for all $t\in (0,(\rho \wedge r)^\beta]$. 
 
 For  $t>(\rho \wedge r)^\beta$, since \eqref{e:Nash-truncated-consequence-1} holds for $t=(\rho \wedge r)^\beta$, using the semigroup property and the $L^\infty$-contraction property, we obtain
 \begin{align*}
 	\lVert P^B_t  \rVert_{L^1(B) \to L^\infty(B)}  &\le 	\lVert P^B_{r^\beta}  \rVert_{L^1(B) \to L^\infty(B)}  	\lVert P^B_{t-r^\beta}  \rVert_{L^\infty(B) \to L^\infty(B)} \le \frac{c_2 }{  V(x_0,r) }  \bigg( 1\wedge \frac{R_0}{r} \bigg)^{-b'/\nu}.
 \end{align*}
  \qed

 	To obtain an on-diagonal upper estimate for the heat kernel of  $(P^{(\rho)}_t)_{t\ge 0}$, we first establish the  uniform ultracontractivity of $(P^{(\rho),B}_t)_{t\ge 0}$ for balls $B$ with arbitrarily large radii. A similar approach was used  in \cite[Lemma 5.6]{GH14} for local Dirichlet forms, and in \cite[Proposition 4.23]{CKW-memo} for non-local Dirichlet forms. However, in our setting, the underlying assumptions  hold only within a localizing constant. This limitation necessitates a significantly different proof, which relies on a self-improvement type argument to address the constraints of the localized framework.

 	We will use the following elementary lemma. The proof is straightforward and therefore omitted.
 	\begin{lem}\label{l:elementary}
 		Let $p_0,q\ge 0$. For any $a,b\in (0,1)$, there exists $C=C(a,b)$ such that the following holds: Define $p_{n+1}:= q+ a(q p_n)^{1/2} + b p_n$ for $n \ge 0$. Then $\lim_{n\to \infty} p_n=Cq$.
 	\end{lem}

  \begin{prop}\label{p:NDU-pre-1}
  	Suppose that \VDb, \TJb, \Sb \ and \Nashb \ hold. There exist  $K\ge 27$ and  $C>0$  such that for any $\rho \in \Q_+ \cap (0,R_0/K)$ and  ball $B:=B(x_0,R)$ with  $R \ge 3K\rho$, 
  	\begin{align*}
  		\lVert P_t^{(\rho),B}\rVert_{L^1(B(x_0, \rho))\to L^\infty(B)} \le  \frac{C\rho^{\beta/\nu}  }{  V(x_0,\rho) t^{1/\nu}} \quad \text{for all $t\in (0,\rho^\beta/K]$}.
  	\end{align*}
  \end{prop}
  \pf   Let   $K\ge 27\vee (15/\nu)$ be a constant  to be chosen later, and let $\rho \in \Q_+ \cap (0,R_0/K)$ and  $B:=B(x_0,R)$ with  $ R \ge 3K\rho$. 
  By Lemma \ref{l:Nash-truncated-consequence}(ii) and  Proposition \ref{p:existence-heat-kernel},  $(P^{(\rho),B}_t)_{t\ge 0}$ has a heat kernel $p^{(\rho),B}(t,x,y)$ on $(0,\infty)\times B \times B$ such that
   \begin{align}\label{e:NDU-pre-1-1}
   	\sup_{x,y\in B}	 p^{(\rho),B}(t,x,y)  \le  \frac{c_1R^{\beta/\nu} }{  V(x_0,R) (t\wedge \rho^\beta)^{1/\nu}}  \bigg( 1\wedge \frac{R_0}{R} \bigg)^{-b'/\nu} \quad \text{for all $t>0$}.
   \end{align}
   Define  for $t>0$, 
   $$	F(t):= \frac{ K^{\beta/\nu}\rho^{\beta/\nu}}{V(x_0,\rho)t^{1/\nu}} \quad \text{ and } \quad \p_0(t):= c_1 F(t) \bigg( \frac{R}{K \rho}\bigg)^{(\beta+b')/\nu}.$$
  By   \eqref{e:NDU-pre-1-1}, since $K\rho<R_0 \wedge R$, we have  for all  $t\in (0, \rho^\beta/K]$, \begin{align}\label{e:NDU-pre-1-2}
   	\esssup_{x\in B, \, y\in U} p^{(\rho),B}(t,x,y) \le  \frac{c_1K ^{\beta/\nu}  \rho^{\beta/\nu} }{  V(x_0,\rho) t^{1/\nu}}  \bigg(\frac{R}{K\rho}\bigg)^{\beta/\nu}\bigg( 1\wedge \frac{K\rho}{R} \bigg)^{-b'/\nu} = \p_0(t).
   \end{align}
   
   Set
   $U:=B(x_0,\rho)$, $V:=B(x_0, 2K \rho)$ and $W:=B(x_0,(2K+3)\rho)$. 
   Let $f \in L^1(U)\cap L^\infty(U)$ with $\lVert f \rVert_1=1$.
   By Corollary  \ref{c:comparison-entrance-2}, we have for all  $t \in (0,\rho^\beta/K]$  and  $\mu$-a.e. $x \in B$,
 	\begin{align}\label{e:NDU-pre-1-3}
 |P^{(\rho),B}_t f(x)| \le 	P^{(\rho),B}_t |f|(x)\le  P^{(\rho),W}_t |f|(x) + 	\sup_{0<s\le t} \lVert P^{(\rho),B}_{s} |f|\rVert_{L^\infty(V_{2\rho}\setminus \overline V)}.
 \end{align}
  By Lemma \ref{l:Nash-truncated-consequence}(ii) and \VDb,  since $(2K+3)\rho< 3K\rho<3R_0$, we have for all $t \in (0,\rho^\beta/K]$,
  \begin{align}\label{e:NDU-pre-1-4}
  \esssup_{B}	P^{(\rho),W}_t |f| &\le \frac{c_2(3K \rho)^{\beta/\nu} }{  V(x_0,K \rho) t^{1/\nu}}    \le 3^{\beta/\nu}c_2 F(t).
  \end{align}
  
  Let  $\{B(z_i,K\rho/12)\}_{i\ge 1}$ be a  covering of $V_{2\rho}\setminus \overline V$ where $z_i \in \overline{V_{2\rho}}\setminus  V$ for all $i\ge 1$. Fix $i \ge 1$ and $t \in (0,\rho^\beta/K]$. By the semigroup property and Proposition \ref{p:comparison-entrance}, we obtain for all $s\in (0,t]$ and $\mu$-a.e. $y \in B(z_i, K \rho/12)$,
  \begin{align}\label{e:NDU-pre-1-5}
  	& P^{(\rho),B}_{s} |f| (y)  = P^{(\rho),B}_{s/2}  (P^{(\rho),B}_{s/2} |f| )(y)  \nn\\
  	&\le  P^{(\rho),B(z_i, K \rho )}_{s/2}  (P^{(\rho),B}_{s/2} |f| )(y)   +  \Big( 1 -P^{(\rho), B(z_i, K \rho )}_{s/2} \1_{B(z_i,K\rho)}(y) \Big) \sup_{0<a\le s/2} \lVert P_{s/2+a}^{(\rho),B} |f| \rVert_{L^\infty(B)}.
  \end{align}
  By Proposition  \ref{p:truncated-off-diagonal-4-t}, we have
\begin{align*}
	&	  1 - P^{(\rho), B(z_i, K \rho )}_{s/2}\1_{B(z_i,K \rho)}(y) \le  2\bigg( \frac{c_3s}{\rho^\beta}\bigg)^{K/15}  \text{\quad for all  $s\in (0,t]$ and $\mu$-a.e. $y \in B(z_i, K \rho/12)$.} 
\end{align*}
Combining this with \eqref{e:NDU-pre-1-2},  since $K \ge 15/\nu$ so that  $s\mapsto s^{K/15}\p_0(s)$ is non-decreasing, we obtain for all $s\in (0,t]$ and $\mu$-a.e. $y \in B(z_i, K \rho/12)$,
\begin{align}\label{e:NDU-pre-1-6}
	&	 \Big( 1 - P^{(\rho), B(z_i, K \rho )}_{s/2}\1_{B(z_i,K \rho)}(y) \Big) \sup_{0<a\le s/2} \lVert P_{s/2+a}^{(\rho),B} |f| \rVert_{L^\infty(B)} \nn\\
	&\le  2\bigg( \frac{c_3s}{\rho^\beta}\bigg)^{K/15}  \sup_{0<a\le s/2}  \p_0(s/2+a) \lVert f \rVert_1  = 2^{1+1/\nu}\bigg( \frac{c_3s}{\rho^\beta}\bigg)^{K/15}   \p_0(s) \nn\\
	&\le 2^{1+1/\nu}\bigg( \frac{c_3t}{\rho^\beta}\bigg)^{K/15}   \p_0(t) \le 2^{1+1/\nu}(c_3/K)^{K/15}   \p_0(t).
\end{align}
On the other hand, by Lemma \ref{l:Nash-truncated-consequence}(ii) and  Proposition \ref{p:existence-heat-kernel},  $( P^{(\rho),B(z_i,K\rho)}_s)_{s\ge 0}$ has a heat kernel $p^{(\rho),B(z_i,K \rho)}(s,x,y)$ such that for all $s>0$,
\begin{align*}
	\sup_{x,y\in B(z_i,K \rho)}	p^{(\rho),B(z_i,K\rho)}(s,x,y)  &\le  \frac{c_4(K\rho)^{\beta/\nu} }{  V(z_i,K\rho) (s\wedge \rho^\beta)^{1/\nu}} .
\end{align*}
Since $d(x_0,z_i) \le (2K+2)\rho<3K \rho$ and $t\le \rho^\beta/K$, by \eqref{e:VD2}, it follows that
\begin{align}\label{e:NDU-pre-1-7}
	\sup_{x,y\in B(z_i,K \rho)} p^{(\rho),B(z_i,K\rho)}(s,x,y)  
	&\le   c_5 F(s) \quad \text{ for all $s \in (0,t]$}.
\end{align}
For all $s\in (0,t]$ and $\mu$-a.e. $y \in B(z_i, K \rho/12)$, using \eqref{e:NDU-pre-1-7} in the first inequality  below,  the Cauchy-Schwarz inequality, $B(z_i,K\rho) \subset B$ and symmetry in the second,  the semigroup property  and $\lVert f \rVert_1=1$ in the third, and  \eqref{e:NDU-pre-1-2} in the fourth, we obtain
\begin{align}\label{e:NDU-pre-1-8}
	&	P^{(\rho),B(z_i, K \rho )}_{s/2}  (P^{(\rho),B}_{s/2} |f| )(y)    = \int_U f(w) \int_{B(z_i,K \rho)} p^{(\rho),B(z_i,K \rho)}(s/2, y, z) p^{(\rho), B} (s/2, z, w)    \mu(dz) \mu(dw) \nn\\
	&\le \left( c_5 F(s/2) \right)^{1/2} \int_U f(w) \int_{B(z_i,K \rho)} p^{(\rho),B(z_i,K\rho)}(s/2, y, z)^{1/2} \,  p^{(\rho), B} (s/2, z, w)   \mu(dz) \mu(dw)\nn\\
		&\le   \left( c_5 F(s/2) \right)^{1/2} \int_U f(w) \bigg[ \int_{B(z_i,K \rho)} p^{(\rho),B}(s/2, y, z)  p^{(\rho), B} (s/2, z, w) \mu(dz) \bigg]^{1/2} \nn\\
	&\qquad \qquad  \qquad \qquad \qquad \quad    \times   \bigg[ \int_{B(z_i,K \rho)} p^{(\rho), B} (s/2, w,z) \mu(dz) \bigg]^{1/2} \mu(dw)\nn\\
		&\le   \left( c_5 F(s/2) \right)^{1/2}  \Big[ \esssup_{w \in U} p^{(\rho),B}(s, y, w)\, P^{(\rho),B}_{s/2} \1_{B(z_i,K \rho)}(w)\Big]^{1/2} \nn\\
			&\le  \left(2^{1/\nu} c_5  F(s) \p_0(s) \right)^{1/2}   \Big[ \esssup_{w \in U} P^{(\rho),B}_{s/2} \1_{B(z_i,K \rho)}(w)\Big]^{1/2}.
\end{align}
Since $B(x_0,K\rho) \cap B(z_i,K \rho) = \emptyset$, by  Proposition \ref{p:truncated-off-diagonal-3-t}, we have
\begin{align*}
\esssup_{w \in U} P^{(\rho),B}_{s/2} \1_{B(z_i,K \rho)}(w)
	&\le  \esssup_{w \in U} P^{(\rho),B}_{s/2} \1_{B(x_0,K \rho)^c}(w) \le  \bigg(\frac{c_6s}{ \rho^\beta}\bigg)^{K/5}.
\end{align*}
Thus,  from  \eqref{e:NDU-pre-1-8}, since   $s\mapsto F(s)^{1/2}\p_0(s)^{1/2}s^{K/10}$ is increasing and $t \le \rho^\beta/K$, we deduce that  for all $s\in (0,t]$ and $\mu$-a.e. $y \in B(z_i, K \rho/12)$,
\begin{align}\label{e:NDU-pre-1-9}
	&P^{(\rho),B(z_i, K \rho )}_{s/2}  (P^{(\rho),B}_{s/2} |f| )(y) \le  \left(2^{1/\nu} c_5 F(s) \p_0(s) \right)^{1/2}  \bigg(\frac{c_6s}{ \rho^\beta}\bigg)^{K/10} \nn\\
	&\le \left(2^{1/\nu} c_5 F(t) \p_0(t) \right)^{1/2} \bigg(\frac{c_6t}{\rho^\beta}\bigg)^{K/10}   \le   \frac{2^{1/(2\nu)} c_5^{1/2}F(t)^{1/2}\p_0(t)^{1/2} }{(K /c_6)^{K/10}}  .
\end{align}
Combining  \eqref{e:NDU-pre-1-5},  \eqref{e:NDU-pre-1-6} and \eqref{e:NDU-pre-1-9},  we deduce that for all $t \in (0,\rho^\beta/K]$,
\begin{align}\label{e:NDU-pre-1-10}
	&	\sup_{0<s\le t} \lVert  P^{(\rho),B}_{s} |f|\rVert_{L^\infty(V_{2\rho}\setminus \overline V)} \le 	\sup_{i\ge 1} \sup_{0<s\le t} \lVert  P^{(\rho),B}_{s} |f|\rVert_{L^\infty(B(z_i,K\rho/12))} \nn\\
	&\le    \frac{2^{1/(2\nu)} c_5^{1/2} F(t)^{1/2}\p_0(t)^{1/2} }{(K/c_6)^{K/10}}   +2^{1+1/\nu}(c_3/K)^{K/15}   \p_0(t).
\end{align}

 We   choose $K\ge 27 \vee (15/\nu)$ such that 
 $$\frac{2^{1/\nu} c_5  }{3^{\beta/\nu}c_2(K/c_6)^{K/5}}  \vee  2^{1+1/\nu}(c_3/K)^{K/15}  \le \frac14.$$ By \eqref{e:NDU-pre-1-3}, \eqref{e:NDU-pre-1-4}  and \eqref{e:NDU-pre-1-10},   since $L^1(B)\cap L^\infty(B)$ is dense in $L^1(B)$, we obtain
 \begin{align}\label{e:NDU-pre-1-11}
 \lVert P_t^{(\rho),B}\rVert_{L^1(U) \to L^\infty(B)} = 	\esssup_{x\in B, \, y\in U} p^{(\rho),B}(t,x,y) &\le   \p_1(t) \quad \text{for all $t\in (0, \rho^\beta/K]$} 
 \end{align}
 where
 \begin{align*}
 	\p_1(t):= 3^{\beta/\nu} c_2F(t) +  \frac12\left( 3^{\beta/\nu} c_2 F(t) \p_0(t) \right)^{1/2}  + \frac14  \p_0(t).
 \end{align*}
 Define for $n \ge 1$,
 \begin{align*}
 \p_{n+1}(t):= 3^{\beta/\nu} c_2  F(t) +  \frac12\left( 3^{\beta/\nu}c_2 F(t) \p_n(t) \right)^{1/2}  + \frac14  \p_n(t).
 \end{align*}
   By repeating the arguments for  \eqref{e:NDU-pre-1-11}, with \eqref{e:NDU-pre-1-2} replaced by \eqref{e:NDU-pre-1-11}, we obtain
 \begin{align*}
\lVert  P_t^{(\rho),B}\rVert_{L^1(U) \to L^\infty(B)} = 	\esssup_{x\in B, \, y\in U} p^{(\rho),B}(t,x,y) &\le   \p_2(t) \quad \text{for all $t\in (0, \rho^\beta/K]$} .
 \end{align*} 
By iterating this procedure and applying Lemma \ref{l:elementary}, we conclude that for all $t \in (0,\rho^\beta/K]$,
\begin{align*}
	\lVert P_t^{(\rho),B}\rVert_{L^1(U) \to L^\infty(B)} \le  \limsup_{n\to \infty}  \p_{n}(t) = c_7 F(t).
\end{align*}
The proof is complete. \qed

\begin{thm}\label{t:NDU-truncated-on}
	Suppose that \VDb, \TJb,  \Sb \  and  \Nashb \ hold. For any $\rho \in \Q_+ \cap (0,R_0/K)$, where $K \ge 27$ is the constant in Proposition \ref{p:NDU-pre-1}, the semigroup  $(P^{(\rho)}_t)_{t\ge 0}$ has a heat kernel $p^{(\rho)}(t,x,y)$ on $(0,\infty)\times M \times M$. Moreover, there 
	exists   $C>0$ independent of $\rho$ such that 
	\begin{align*}
		p^{(\rho)}(t,x,y) \le   \frac{C\rho^{\beta/\nu}}{V(x,\rho) (t\wedge \rho^\beta)^{1/\nu}} \quad \text{for all $t>0$ and $x,y \in M$.} 
	\end{align*}
\end{thm}
\pf Let $\{z_i\}_{i\ge 1}\subset M$ be a countable family such that $M= \cup_{i\ge 1}B(z_i,\rho)$. By  duality,  for all $i\ge 1$ and $t>0$, we have $	\lVert P^{(\rho)}_t\rVert_{L^1(M)\to L^\infty(B(z_i,\rho))}= 	\lVert P^{(\rho)}_t\rVert_{L^1(B(z_i,\rho))\to L^\infty(M)}.$ Thus, by Proposition \ref{p:NDU-pre-1},  we get that for all $i\ge 1$ and  $t\in (0,\rho^\beta/K]$,
\begin{align*}
	&	\lVert P^{(\rho)}_t\rVert_{L^1(M)\to L^\infty(B(z_i,\rho))}= 	\lVert P^{(\rho)}_t\rVert_{L^1(B(z_i,\rho))\to L^\infty(M)} \\
	&= \lim_{R\to \infty} 	\lVert P^{(\rho),B(z_i,R)}_t\rVert_{L^1(B(z_i,\rho))\to L^\infty(M)} \le \frac{c_1\rho^{\beta/\nu}}{V(z_i,\rho)t^{1/\nu}}.
\end{align*}
Hence, by Proposition \ref{p:existence-heat-kernel},  $(P^{(\rho)}_t)_{t\ge0}$  has a heat kernel $p^{(\rho)}(t,x,y)$ on $(0,\infty) \times M\times M$ such that
\begin{align}\label{e:NDU-truncated-on-1}
	\sup_{x\in B(z_i,\rho),\, y \in M}	p^{(\rho)}(t,x,y)  \le \frac{c_1\rho^{\beta/\nu}}{V(z_i,\rho)t^{1/\nu}} \quad \text{ for all $i\ge 1$ and   $t\in (0,\rho^\beta/K]$.} 
\end{align}

For $x\in M$, let $i(x)$ be such that $x\in B(z_{i(x)},\rho)$.
For  all $t\in (0,\rho^\beta/K]$ and  $x,y \in M$, by \eqref{e:NDU-truncated-on-1} and \eqref{e:VD2}, we obtain
\begin{align}\label{e:NDU-truncated-on-2}
	p^{(\rho)}(t,x,y)\le \frac{c_1\rho^{\beta/\nu}}{V(z_{i(x)},\rho)t^{1/\nu}} \le \frac{c_2\rho^{\beta/\nu}}{V(x,\rho)t^{1/\nu}}.
\end{align}
Moreover, for all    $t>\rho^\beta/K$ and  $x,y \in M$, by the semigroup property, symmetry and  \eqref{e:NDU-truncated-on-2}, we get
\begin{align*}
		p^{(\rho)}(t,x,y) &= \int_{M}p^{(\rho)}(\rho^\beta/K, x,z)p^{(\rho)}(t-\rho^\beta/K, y,z) \mu(dz) \nn\\
		&\le \frac{c_2K^{1/\nu}}{V(x,\rho)}\int_{M}p^{(\rho)}(t-\rho^\beta/K, y,z) \mu(dz)\le  \frac{c_2K^{1/\nu}}{V(x,\rho)} .
\end{align*}
This completes the proof.  \qed

\subsection{Some off-diagonal estimates}
In this subsection, we establish \TE \  and an   off-diagonal upper estimate for the heat kernels of the semigroups associated with truncated Dirichlet forms.

\begin{lem}\label{l:TEb}
	Suppose that  \TJb \  and \TEb \ hold. For any $k\ge 1$, there exists  $C=C(k)>0$ such that for all $\rho \in \Q_+\cap (0,R_0)$, $t>0$ and ball $B:=B(x_0,r)$ with $r\in (0,k\rho]$,
	\begin{align*}
		\esssup_{ 4^{-1} B}	P^{(\rho)}_t \1_{B^c}  \le \frac{Ct}{r^\beta}.
	\end{align*}
\end{lem}
\pf By  \eqref{e:comparison-tail-2}, \TJb \ and \TEb, since $r\le k\rho<kR_0$, we obtain
\begin{align*}
	\esssup_{ 4^{-1} B}	P^{(\rho)}_t \1_{B^c}&\le  
	\esssup_{ 4^{-1} B}	P_t \1_{B^c}  + \frac{c_1t}{\rho^\beta}\le \frac{c_2 t}{(r\wedge R_0)^\beta} + \frac{c_1t}{\rho^\beta} \le   \frac{c_3t}{r^\beta}.
\end{align*}
\qed

\begin{prop}\label{p:TEb}
	Suppose that  \TJb \  and \Sb \ hold. Then \TEb \ holds. Moreover,  for any $\eps>0$, there exist  $\eta_0>0$, independent of $\eps$, and $C=C(\eps)>0$ such that for any ball  $B:=B(x_0,r)$,
	\begin{align}\label{e:TEb-general}
		\esssup_{ 4^{-1} B}	P_t \1_{B^c}  \le \frac{Ct}{r^\beta} \bigg( 1 + \frac{r}{R_0}\bigg)^{\eps} \quad \text{ for all  $t\in (0, \eta_0 R_0^\beta)$}.
	\end{align}
\end{prop}
\pf  Let $B:=B(x_0,r)$ 
and $\rho \in \Q_+$ be such that $(r\wedge R_0)/10\le \rho\le (r\wedge R_0)/9$.  Using \eqref{e:comparison-tail-2}, \TJb \ and Proposition \ref{p:truncated-off-diagonal-3-t}, we obtain for all $t>0$,
\begin{align}\label{e:TEb-1}
	\esssup_{4^{-1}B}  P_{t} \1_{B^c}&\le   	\esssup_{4^{-1}B} P^{(\rho)}_{t} \1_{B^c}+ \frac{c_1t}{(r\wedge R_0)^\beta}\le  \bigg( \frac{c_2t}{(r\wedge R_0)^\beta}\bigg)^{r/(5\rho)} + \frac{ c_1t}{(r\wedge R_0)^\beta}.
\end{align}
If $c_2t \ge (r\wedge R_0)^\beta$, then by taking $C$ larger than $ c_2$, we get \eqref{e:TEb}. If  $c_2t < (r\wedge R_0)^\beta$, then  since $r/(5\rho)\ge 1$, \eqref{e:TEb} follows from
\eqref{e:TEb-1}. Thus, \TEb \ holds.

We next prove \eqref{e:TEb-general}. Let $\eps >0$. We assume, without loss of generality, $\eps<\beta$. If  $r<R_0$, then \eqref{e:TEb-general} follows from \TEb.  Assume $r\ge R_0$. Let $\rho \in \Q_+$ be such that $r^{1-\eps/\beta}R_0^{\eps/\beta}/10\le \rho\le r^{1-\eps/\beta}R_0^{\eps/\beta}/9$.  By \eqref{e:comparison-tail-2},  \TJb \ and \TEb,  we have for all  $z\in M$, $R\ge \rho/3$ and $t>0$,
\begin{align*}
	\esssup_{ B(z,R/4)} P^{(\rho)}_{t} \1_{B(z,R)^c} &\le 
	\esssup_{ B(z,R/4)} P_{t} \1_{B(z,R)^c}  + \frac{c_3t}{\rho^\beta} \le 
	\frac{c_4t}{((\rho/3)\wedge R_0)^\beta}  + \frac{c_5t}{r^{\beta-\eps}R_0^\eps} \le \frac{c_6t}{R_0^\beta}.
\end{align*}
Thus, by  Lemmas \ref{l:Markov} and \ref{l:conservative}, we get that  $	\essinf_{B(z,\rho/12)}  P^{(\rho),B(z,\rho)}_{t} \1_{B(z,\rho)}  \ge 1 - 2c_6R_0^{-\beta}t$  for all $z\in M$ and $t>0$. By Lemma \ref{l:truncated-off-diagonal}, it follows that for all $z\in M$, $R\ge 3\rho$ and $t>0$,
\begin{align}\label{e:TEb-general-1}
	\esssup_{ B(z,R/4)}  P^{(\rho)}_{t} \1_{B(z,R)^c}\le \left(\frac{2c_6t}{R_0^\beta}\right)^{R/(3\rho)-1} .
\end{align}
Note that the constant $c_6$ is independent of $\eps$. Set $\eta_0:=1/(2ec_6)$. For all $t\in (0,\eta_0R_0^\beta)$, using  \eqref{e:comparison-tail-2} and  \TJb  \ in the first inequality below, \eqref{e:TEb-general-1} in the second, $r\ge 9\rho$ in the third, and the fact that  $\sup_{a\ge 1} a^{(\beta-\eps)\beta/\eps}e^{-a}<\infty$ in the fifth, we  obtain 
\begin{align*}
	\esssup_{ 4^{-1} B}	P_t \1_{B^c} &\le \esssup_{ 4^{-1} B}	P^{(\rho)}_t \1_{B^c}  + \frac{c_7t}{\rho^\beta} \le \left(\frac{t}{e\eta_0R_0^\beta}\right)^{r/(3\rho)-1}   + \frac{c_8t}{r^{\beta-\eps}R_0^\eps}\\
	& \le \frac{t}{e\eta_0R_0^\beta} \left(\frac{t}{e\eta_0R_0^\beta}\right)^{r/(9\rho)}   + \frac{c_8t}{r^{\beta-\eps}R_0^\eps} \le \frac{te^{-r/(9\rho)} }{e\eta_0R_0^\beta}  + \frac{c_8t}{r^{\beta-\eps}R_0^\eps}\\
	&\le \frac{c_9t}{\eta_0R_0^\beta} \bigg( \frac{9\rho}{r}\bigg)^{(\beta-\eps)\beta/\eps}  + \frac{c_8t}{r^{\beta-\eps}R_0^\eps} =\frac{c_{10}t}{r^{\beta-\eps}R_0^\eps}.
\end{align*}  \qed

\begin{prop}\label{p:NDU-truncated-off-1}
	Suppose that \VDb, \TJb,  \Sb \  and  \Nashb \ hold. Then there exists $K_0\ge 1$ such that for any $\rho \in \Q_+ \cap (0,R_0/K)$, where $K \ge 27$ is the constant in Proposition \ref{p:NDU-pre-1}, the following holds:
For any $k\ge 20$, there 
	exists a constant  $C=C(k)>0$ independent of $\rho$ such that for all $t\in (0,\rho^\beta/K_0]$ and $x,y \in M$ with $d(x,y) \ge k\rho$,
	\begin{align*}
		p^{(\rho)}(t,x,y) \le   \frac{C\rho^{\beta/\nu-\alpha'}t^{k/20-1/\nu}}{V(x,\rho)d(x,y)^{\beta k/20 -\alpha'}}.
	\end{align*}
\end{prop}
\pf  Let $K_0\ge 1$ be a constant  whose  value to be determined later and  $t\in (0,\rho^\beta/K_0]$.   For any $x,y \in M$ with $d(x,y)\ge k\rho$, using the semigroup property in the first inequality below, Theorem \ref{t:NDU-truncated-on} in the second, and \eqref{e:VD2} and Proposition \ref{p:esssup1} in the third, we get
\begin{align*}
	&	p^{(\rho)}(t,x,y) \le \bigg(\int_{B(x,d(x,y)/2)^c} + \int_{B(y,d(x,y)/2)^c}\bigg) p^{(\rho)}(t/2,x,z) p^{(\rho)}(t/2, y,z)\mu(dz)  \nn\\
	&\le \frac{c_1\rho^{\beta/\nu}}{V(y,\rho)t^{1/\nu}} \int_{B(x,d(x,y)/2)^c} p^{(\rho)}(t/2,x,z)\mu(dz) +  \frac{c_1\rho^{\beta/\nu}}{V(x,\rho)t^{1/\nu}} \int_{B(y,d(x,y)/2)^c} p^{(\rho)}(t/2,y,z)\mu(dz)\nn\\
	&\le \frac{c_2\rho^{\beta/\nu}}{V(x,\rho)t^{1/\nu}}\bigg(\frac{d(x,y)}{\rho}\bigg)^{\alpha'} \bigg( \esssup_{B(x,d(x,y)/8)} P^{(\rho)}_{t/2} \1_{B(x,d(x,y)/2)^c} +  \esssup_{B(y,d(x,y)/8)} P^{(\rho)}_{t/2} \1_{B(y,d(x,y)/2)^c} \bigg) .
\end{align*}
Thus, by Proposition \ref{p:truncated-off-diagonal-3-t},  we deduce that for any $x,y \in M$ with $d(x,y)\ge k\rho$,
\begin{align}\label{e:NDU-truncated-2}
	p^{(\rho)}(t,x,y) &\le \frac{c_3\rho^{\beta/\nu-\alpha'}d(x,y)^{\alpha'}}{V(x,\rho)t^{1/\nu}} \bigg( \frac{c_4t}{2\rho^\beta}\bigg)^{d(x,y)/(10\rho)}.
\end{align}

Set $K_0:=1 \vee c_4$. Then $c_4 t/(2\rho^\beta) \le 1/2$. Using the fact that $\sup_{\eps \in (0,1/2],\, r\ge 2} \eps^{r-1}r^{\beta}<\infty$, we deduce from \eqref{e:NDU-truncated-2} that for any $x,y \in M$ with $d(x,y)\ge k\rho$,
\begin{align*}
&	p^{(\rho)}(t,x,y) \le \frac{c_3\rho^{\beta/\nu-\alpha'}d(x,y)^{\alpha'}}{V(x,\rho)t^{1/\nu}}\left[ \bigg( \frac{ c_4t}{2\rho^\beta}\bigg)^{2d(x,y)/(k\rho)}\right]^{k/20}\nn\\
	&\le\frac{c_3\rho^{\beta/\nu-\alpha'}d(x,y)^{\alpha'}}{V(x,\rho)t^{1/\nu}}  \left[ \frac{ c_5t}{2\rho^\beta} \bigg( \frac{k\rho}{2d(x,y)}\bigg)^\beta \right]^{k/20}=\frac{c_6\rho^{\beta/\nu-\alpha'}t^{k/20-1/\nu}}{V(x,\rho)d(x,y)^{\beta k/20-\alpha'}}.
\end{align*}
 \qed

 \subsection{Diagonal upper heat kernel estimate}

In this subsection, we establish the following  diagonal upper estimate for the heat kernel of $(P_t)_{t\ge 0}$.

\begin{thm}\label{t:NDU}
	Suppose that \VDb,  \Sb \ and \Nashb \ hold, and that one of the following conditions is satisfied:

\smallskip

\setlength{\leftskip}{3mm}

\noindent {\rm (1)} \TJb \  and \IVJb \ hold with $\gamma \ge 0 $ satisfying 
\begin{align}\label{e:NDU-range}
(1-\nu)\gamma< 1+\nu.
\end{align}

\noindent {\rm (2)} \dTJb \ and \RVDb \ hold with $q \in [1,\infty)$ and $\alpha_0'\ge 0$ satisfying 
\begin{align}\label{e:NDU-range-2}
\frac{(1-\nu)}{2\beta}  \left[ \bigg( \frac{2}{q} -1 \bigg) \alpha' - \alpha_0' \right]_+ <1+\nu.
\end{align}

\setlength{\leftskip}{0mm}

\noindent  Then $(P_t)_{t\ge0}$ has a heat kernel $p(t,x,y)$ on $(0,\infty)\times M\times M$ and there exists   $C>0$ such that 
	\begin{align*}
		p(t,x,x) \le  \frac{C}{V(x,t^{1/\beta} \wedge R_0)}\bigg(1+\frac{t}{R_0^\beta}\bigg) \quad \text{for all $t>0$ and $x \in M$}. 
	\end{align*}
\end{thm}

By Remark \ref{r:dTJ}, both  conditions (1) and (2) in the above theorem imply \TJb.

To prove Theorem \ref{t:NDU}, we first establish several  lemmas. The proof of Theorem \ref{t:NDU} will be given at the end of this subsection.

\begin{lem}\label{l:NDU-range}
(i) If \eqref{e:NDU-range} holds, then there exist  $\delta \in [0,1 \wedge \nu)$ and $\lambda \in (0,1)$ such that
\begin{align}\label{e:NDU-range-lambda}
	2 \lambda \delta  - (\gamma-1) (1-\delta)>0.
\end{align}

\noindent (ii) If \eqref{e:NDU-range-2} holds, then there exist  $\delta \in [0,1 \wedge \nu)$ and $\lambda \in (0,1)$ such that
\begin{align}\label{e:NDU-range-lambda-2}
2\lambda \delta  - \frac{ (1-\delta)}{2\beta}  \left[ \bigg( \frac{2}{q}-1 \bigg) \alpha' -\alpha_0' - 2\beta  \right]   >0.
\end{align}
\end{lem}
\pf  (i) If $\nu \ge 1$, then  \eqref{e:NDU-range-lambda} is satisfied with $\lambda =1/2$ and $\delta = 1 - (2 + |\gamma-1|)^{-1}$. If $\nu<1$, then since $2 \nu - (\gamma-1)(1-\nu)>0$ by \eqref{e:NDU-range},  there exist $\delta \in [0,\nu)$ and $\lambda \in (0,1)$ satisfy \eqref{e:NDU-range-lambda}.

\noindent  (ii) If $\nu \ge 1$, then  \eqref{e:NDU-range-lambda-2} is satisfied with $\lambda =1/2$ and $\delta = 1 - (2 + |\frac{\alpha'-\alpha_0'}{2\beta}-1|)^{-1} $. If $\nu<1$, then since $
2 \nu  - \frac{ (1-\nu)}{2\beta}  [ ( \frac{2}{q}-1 ) \alpha' -\alpha_0' - 2\beta  ]      >0$ by \eqref{e:NDU-range-2},  we get the result. \qed

\begin{lem}\label{l:NDU-1}
	Suppose that \Nashb \ holds. Then there exists $C>0$ such that for any ball $B:=B(x_0,R)$, the heat kernel $p^B(t,x,y)$ exists and satisfies that
	\begin{equation*}
	p^B(t,x,y)\le \frac{CR^{\beta/\nu} }{  V(x_0,R) (t\wedge R^\beta)^{1/\nu}}  \bigg( 1\wedge \frac{R_0}{R} \bigg)^{-b'/\nu}  \;\; \text{for all $t>0$ and $x,y \in B$}.
	\end{equation*}
\end{lem}
\pf The result follows from Lemma \ref{l:Nash-truncated-consequence}(i) and Proposition \ref{p:existence-heat-kernel}.\qed

For each $x\in M$ and $\rho>0$, define
\begin{align}\label{e:def-F-rho}
	F_{x,\rho}(t):=\frac{\rho^{\beta/\nu}}{V(x, \rho)(t\wedge \rho^\beta)^{1/\nu}} \bigg(1+ \frac{t}{\rho^\beta}\bigg), \qquad t>0.
\end{align}
 
 For the remainder of this subsection,   let $K$ denote the constant from Proposition \ref{p:NDU-pre-1}.

The following lemma is a crucial   step in the proof of Theorem \ref{t:NDU}.

\begin{lem}\label{l:NDU-2}
Under the setting of Theorem  \ref{t:NDU}, the following holds: Let $\rho \in \Q_+ \cap (0, R_0/K)$, $B:=B(x_0,R)$  and  $\p:(0,\infty)\to (0,\infty)$ be a non-increasing function such that, for some $c_0\ge 1$,
	\begin{align}\label{e:l:NDU-p-scaling}
		\p(t/2) \le c_0\p(t) \quad \text{for all $t>0$}.
	\end{align}
	Then  there exists $C>0$, independent of $\rho,B$ and depending on $\p_0$ only through $c_0$, such that if the heat kernel $p^B(t,x,y)$ of $(P^B_t)_{t\ge 0}$   satisfies 
	\begin{align}\label{e:l:NDU-2-1}
		p^B(t,x,x) \le  AF_{x,\rho}(t) + \p(t) \quad \text{for all $t >0$ and $x \in B$}
	\end{align}
	with a constant $A\ge 1$, then for all $t>0$, $x \in B$ and $y\in B(x,\rho)$,
	\begin{align*}
&\int_0^t \int_B p^{(\rho)}(s,x,z) \int_B  J^{(\rho)}_2(z,dw) p^B(t-s,w,y)\mu(dz)\, ds\le CF_{x,\rho}(t) +  \frac12\left( AF_{x,\rho}(t)   + \p(t)   \right).
	\end{align*}
\end{lem}
\pf  Let $K_0$ be the constant in Proposition \ref{p:NDU-truncated-off-1}, $k:=20(1  + 1/\nu  + \alpha'/\beta)$ and  $\eps \in (0,\eta_0 \wedge (1/2))$ be a constant to be chosen later, where $\eta_0$ is the constant in Proposition \ref{p:TEb}.    Fix  $t>0$, $x \in B$ and $y\in B(x,\rho)$. We have
\begin{align*}
	\int_0^t \int_B p^{(\rho)}(s,x,z) \int_B  J^{(\rho)}_2(z,dw) p^B(t-s,w,y)\mu(dz)\, ds =:I_1+I_2+I_3 + I_4,
\end{align*}
where
\begin{align*}
	I_1&:=	\int_{\eps ( t \wedge (\rho^\beta/K_0))}^{t} \int_B p^{(\rho)}(s,x,z) \int_B  J^{(\rho)}_2(z,dw) p^B(t-s,w,y)\mu(dz)\,ds,\\
	I_2&:=	\int_{0}^{\eps ( t \wedge (\rho^\beta/K_0))} \int_{B\setminus B(x,k\rho)} p^{(\rho)}(s,x,z) \int_B  J^{(\rho)}_2(z,dw) p^B(t-s,w,y)\mu(dz)\,ds,\\
		I_3&:=	\int_{0}^{\eps ( t \wedge (\rho^\beta/K_0))} \int_{B \cap B(x,k\rho)} p^{(\rho)}(s,x,z) \int_{B\cap B(y,4k\rho)}  J^{(\rho)}_2(z,dw) p^B(t-s,w,y)\mu(dz)\,ds,\\
			I_4&:=	\int_{0}^{\eps ( t \wedge (\rho^\beta/K_0))} \int_{B \cap B(x,k\rho)} p^{(\rho)}(s,x,z) \int_{B\setminus B(y,4k\rho)}  J^{(\rho)}_2(z,dw) p^B(t-s,w,y)\mu(dz)\,ds.
\end{align*}

For $I_1$, using Theorem \ref{t:NDU-truncated-on} and symmetry in the first inequality below and  \TJb \ in the second, we obtain
 \begin{align}\label{e:estimate-I1}
I_1	 &\le \frac{c_1\rho^{\beta/\nu}}{V(x,\rho)}	\int_{\eps ( t \wedge (\rho^\beta/K_0))}^{t} s^{-1/\nu}\int_B \int_B  J^{(\rho)}_2(w,dz) \, p^B(t-s,y,w)\mu(dw) \, ds\nn\\
	 &\le \frac{c_2\eps^{-1/\nu} \rho^{\beta/\nu-\beta}}{V(x,\rho)(t\wedge (\rho^\beta/K_0))^{1/\nu} } 	\int_{\eps ( t \wedge (\rho^\beta/K_0))}^{t} \int_Bp^B(t-s,y,w)\mu(dw) \, ds\nn\\
	 &\le \frac{c_2\eps^{-1/\nu} \rho^{\beta/\nu-\beta}t}{V(x,\rho)(t\wedge (\rho^\beta/K_0))^{1/\nu} } \le  c_2K_0^{1/\nu}\eps^{-1/\nu} F_{x,\rho}(t).
\end{align}

For $I_2$, using Proposition \ref{p:NDU-truncated-off-1}  in the first inequality below, $\beta k/20-\alpha'>0$ and symmetry in the second,   \TJb \ in the third, and  $k/20-1/\nu >1$ and $\eps<1\le K_0$ in the fourth, we get
\begin{align}\label{e:estimate-I2}
	I_2	&\le \frac{c_3\rho^{\beta/\nu-\alpha'}}{V(x,\rho)} \int_{0}^{\eps ( t \wedge (\rho^\beta/K_0))}   \int_{B(x,k\rho)^c} \frac{s^{k/20-1/\nu} }{d(x,z)^{\beta k/20 - \alpha'}}\int_B  J^{(\rho)}_2(z,dw)  p^B(t-s,w,y)\mu(dz)\,ds \nn\\
	&\le \frac{c_3\rho^{\beta/\nu-\alpha'}}{(k\rho)^{\beta k/20-\alpha'}V(x,\rho)} \int_{0}^{\eps ( t \wedge (\rho^\beta/K_0))} s^{k/20-1/\nu} \int_B   \int_{B(x,k\rho)^c}  J^{(\rho)}_2(w,dz)  p^B(t-s,y,w)\mu(dw)\,ds \nn\\
			&\le \frac{c_4\rho^{\beta/\nu-\beta-\beta k/20}}{V(x,\rho) } \int_{0}^{\eps ( t \wedge (\rho^\beta/K_0))} s^{k/20-1/\nu}  \int_B p^B(t-s,y,w)\mu(dw)\, ds\nn\\
			& \le \frac{c_5 \rho^{\beta/\nu-\beta-\beta k/20} (t\wedge \rho^\beta)^{k/20 +1-1/\nu}}{V(x,\rho)} \le \frac{c_5\rho^{\beta/\nu} }{V(x,\rho)(t\wedge \rho^\beta)^{1/\nu}} \le c_5 F_{x,\rho}(t).
\end{align}

By using the semigroup property, the  Cauchy-Schwarz inequality, \eqref{e:l:NDU-2-1}, the monotonicity of $\p$ and \eqref{e:l:NDU-p-scaling},  we observe that for all $s \in (0,t/2)$ and $w \in B$,
\begin{align}\label{e:l:NDU-2}
	p^B(t-s,w,y) & \le p^B(t-s,w,w)^{1/2}p^B(t-s,y,y)^{1/2} \nn\\
	&\le  \left( AF_{w,\rho}(t/2) + \p(t/2)  \right)^{1/2} \left( AF_{y,\rho}(t/2) + \p(t/2)  \right)^{1/2}\nn\\
		&\le  c_0 2^{1/\nu}  \left( AF_{w,\rho}(t) + \p(t)  \right)^{1/2} \left( AF_{y,\rho}(t) + \p(t)  \right)^{1/2}.
\end{align}
Further, by \eqref{e:VD2}, there exists $c_6\ge 1$ such that for all $w \in B(y,4k\rho)$,
\begin{align}\label{e:l:NDU-3}
	F_{w,\rho}(t) \le c_6 F_{x,\rho}(t).
\end{align}
Using \eqref{e:l:NDU-2} and \eqref{e:l:NDU-3} in the first inequality below, and  \TJb \ in the second, we obtain
\begin{align}\label{e:estimate-I3}
I_3	&\le c_0 2^{1/\nu} \left( c_6AF_{x,\rho}(t) + \p(t)  \right) \int_{0}^{\eps \rho^\beta/K_0} \int_{B(x,k\rho)} p^{(\rho)}(s,x,z) \int_{B(y,4k\rho)}  J^{(\rho)}_2(z,dw)    \mu(dz)\,ds \nn\\
	&\le \frac{c_7 ( AF_{x,\rho}(t) + \p(t)  )}{\rho^\beta}   \int_{0}^{\eps \rho^\beta/K_0} ds\le \frac{c_7\eps ( AF_{x,\rho}(t) + \p(t)  )}{K_0}.
\end{align}

 Define for  $m\ge 2$ and $s\in (0,\rho^\beta)$, 
\begin{align*}
f_{m}(s)&:=	\int_{B\cap B(x,k\rho)} p^{(\rho)}(s,x,z)  \int_{B \cap (B(y,2^{m+1} k\rho)\setminus B(y,2^m k\rho))}  J^{(\rho)}_2(z,dw) p^B(t-s,w,y)\mu(dz).
\end{align*}
we consider the following two cases  separately:

\smallskip

Case 1: \TJb  \ and \IVJb \ hold with  $\gamma$ satisfying \eqref{e:NDU-range}. By Lemma \ref{l:NDU-range}(i), there exist constants $\delta \in [0,1\wedge \nu)$ and $\lambda \in (0,1)$ satisfying \eqref{e:NDU-range-lambda}. For all  $m\ge 2$ and $s\in (0, \eps(t \wedge (\rho^\beta/K_0)))$, using symmetry of $J$, $d(x,y)<\rho<k\rho$ and Theorem \ref{t:NDU-truncated-on}  in the first inequality below,   \TJb \ and  Proposition \ref{p:esssup1} in the second,   Proposition \ref{p:TEb}   in the third and $\rho<R_0$  in the fourth, we obtain 
\begin{align}\label{e:estimate-fm-1}
	f_{m}(s)&\le \frac{c_{8}\rho^{\beta/\nu}}{V(x,\rho)s^{1/\nu} }  \int_{ B \cap (B(y,2^{m+1} k\rho)\setminus B(y,2^m k\rho))}  J\big(w,B(w,(2^m-2)k\rho)^c\big) \, p^B(t-s,y,w) \mu(dw) \nn\\
	&\le \frac{c_{9}\rho^{\beta/\nu-\beta}}{2^{\beta m}V(x,\rho)s^{1/\nu} }  \esssup_{B(y,2^{m-2}k\rho)}P_{t-s} \1_{B(y,2^mk\rho)^c}  \nn \\
	&\le \frac{c_{10}\rho^{\beta/\nu-2\beta}t}{ 2^{2\beta m}   V(x,\rho)s^{1/\nu} }  \bigg(1 + \frac{2^m\rho}{R_0}\bigg)^{2\beta (1-\lambda)} \le \frac{c_{11}\rho^{\beta/\nu-2\beta}t}{   2^{2\beta \lambda m}V(x,\rho)s^{1/\nu} }.
\end{align}
On the other hand,   for all $m\ge 2$ and $s\in (0,t/2)$, using \eqref{e:l:NDU-2} and \eqref{e:l:NDU-3} in the first inequality below, and Lemma \ref{l:IVJ} (with $a=2$, $r=\rho^\beta$ and $R=(2^m-2)^\beta k^\beta \rho^\beta$)  and \TJb \ in the third, we get
\begin{align}\label{e:estimate-fm-2}
f_m(s)&\le c_0  2^{1/\nu} \left( c_6AF_{x,\rho}(t) + \p(t)  \right)^{1/2}\int_{B\cap B(x,k\rho)} p^{(\rho)}(s,x,z)\nn\\
&\quad  \times \bigg[ \int_{B \cap (B(y,2^{m+1} k\rho)\setminus B(y,2^m k\rho))} \left(   A F_{w,\rho}(t)  + \p(t) \right)^{1/2} J^{(\rho)}_2(z,dw) \bigg] \,\mu(dz)\nn\\
&\le  c_0  2^{1/\nu}\left( c_6AF_{x,\rho}(t) + \p(t)  \right)^{1/2}\int_{B\cap B(x,k\rho)} p^{(\rho)}(s,x,z)\nn\\
&\quad  \times \bigg[ \int_{B \cap (B(z,2(2^m-2) k\rho)\setminus B(z,(2^m-2) k\rho))} \left(   (A F_{w,\rho}(t))^{1/2}  + \p(t)^{1/2} \right) J^{(\rho)}_2(z,dw) \bigg] \,\mu(dz)\nn\\
&\le \frac{c_{12}  (2^mk)^{\gamma \beta} \left( c_6 AF_{x,\rho}(t) + \p(t)  \right)^{1/2} ((AF_{x,\rho}(t))^{1/2}  + \p(t)^{1/2}) }{2^{m\beta} k^\beta \rho^{ \beta}}\int_{M} p^{(\rho)}(s,x,z) \mu(dz) \nn\\
&\le \frac{c_{13}  2^{(\gamma-1) \beta m} \left( AF_{x,\rho}(t) + \p(t)  \right) }{\rho^{ \beta}} .
\end{align}
Combining \eqref{e:estimate-fm-1} and \eqref{e:estimate-fm-2}, we obtain for all  $m\ge 2$ and $s\in (0, \eps(t \wedge (\rho^\beta/K_0)))$,
\begin{align*}
	f_{m}(s) &\le  \bigg[\frac{c_{11}\rho^{\beta/\nu-2\beta}t}{   2^{2\beta \lambda m}V(x,\rho)s^{1/\nu} } \bigg]^{\delta}  \left[\frac{c_{13}  2^{(\gamma-1) \beta m} \left( AF_{x,\rho}(t) + \p(t)  \right) }{\rho^{ \beta}} \right] ^{1-\delta}\nn\\
	&=  \frac{c_{14}\rho^{ \beta \delta/\nu - (1+\delta)\beta}  t^{ \delta}(    AF_{x,\rho}(t)  + \p(t)  )^{1-\delta}}{2^{(2 \lambda \delta  - (\gamma-1) (1-\delta))\beta m}V(x,\rho)^{\delta}s^{\delta/\nu} } .
\end{align*}
By \eqref{e:NDU-range-lambda}, since $\delta<\nu$, it follows that
\begin{align}\label{e:estimate-I4-32}
	I_4&=  \int_0^{\eps(t \wedge (\rho^\beta/K_0))}\sum_{m\ge 2} f_{m}(s) \,ds \le   \frac{c_{15}\rho^{ \beta \delta/\nu - (1+\delta)\beta}  t^{ \delta}(    AF_{x,\rho}(t)  + \p(t)  )^{1-\delta}}{V(x,\rho)^{\delta} }  \int_0^{\eps(t \wedge \rho^\beta)} \frac{ds}{s^{\delta/\nu}}\nn\\
	&= \frac{c_{16} \eps^{1-\delta/\nu} (t \wedge \rho^\beta) \rho^{\beta\delta/\nu - \beta} t^{ \delta}   (    AF_{x,\rho}(t) + \p(t)  )^{1-\delta}}{V(x,\rho)^{\delta}  (t \wedge \rho^\beta)^{\delta/\nu} \rho^{\delta \beta}}   \nn\\
&\le c_{16} \eps^{1-\delta/\nu}\bigg[ \frac{\rho^{\beta/\nu}}{V(x,\rho) (t\wedge \rho^\beta)^{1/\nu}} \bigg( \frac{t}{\rho^\beta}\bigg)\bigg]^\delta (    AF_{x,\rho}(t) + \p(t)  )^{1-\delta}\nn\\
&\le c_{16}\eps^{1-\delta/\nu}(AF_{x,\rho}(t) + \p(t)).
\end{align}
By combining \eqref{e:estimate-I1}, \eqref{e:estimate-I2}, \eqref{e:estimate-I3} and  \eqref{e:estimate-I4-32},  and  choosing  $\eps$ sufficiently small, we arrive at the desired result in this case.

\smallskip

Case 2: \dTJb \ and \RVDb \ hold with $q$ and  $\alpha_0'$ satisfying \eqref{e:NDU-range-2}. By Lemma \ref{l:NDU-range}(ii), there exist constants $\delta \in [0,1\wedge \nu)$ and $\lambda \in (0,1)$ satisfying \eqref{e:NDU-range-lambda-2}. Note that  \eqref{e:estimate-fm-1} is still valid for all  $m\ge 2$ and $s\in (0, \eps(t \wedge (\rho^\beta/K_0)))$.  By  \eqref{e:VD2} and  \eqref{e:RVD2}, we have for all $m\ge 2$, $w\in B(y,2^{m+1} k\rho)$ and $s>0$,
\begin{align}\label{e:estimate-I4-3-F}
	F_{w,\rho}(s)  \le \frac{c_{17} 2^{\alpha' m}\rho^{\beta/\nu}}{V(w, 2^{m+2}k\rho)(s\wedge \rho^\beta)^{1/\nu}} \bigg(1+ \frac{s}{\rho^\beta}\bigg) \le  c_{18}2^{(\alpha'-\alpha_0')m} F_{x,\rho}(s).
\end{align}
For all $z\in B\cap B(x,k\rho)$, $m\ge 2$ and $s\in (0,t/2)$, using  symmetry and H\"older's inequality in the first inequality below, \dTJb \ and  the semigroup property in the second,  \eqref{e:l:NDU-2-1}, \eqref{e:l:NDU-2},  the monotonicity of $\p$ and the fact that  $t<2t-2s<2t$ in the third, \RVDb,  \eqref{e:l:NDU-3}  and \eqref{e:estimate-I4-3-F} in the fourth, and  the fact that $ V(x,\rho)^{-1} \le AF_{x,\rho}(t)$ in the sixth, we obtain
\begin{align*}
&	 \int_{B \cap (B(y,2^{m+1} k\rho)\setminus B(y,2^m k\rho))}  p^B(t-s,w,y) J(z,w) \mu(dw) \nn\\
&\le  \bigg[ \int_{ B(z,(2^m-2)k\rho)^c} J(z,w)^{q}\, \mu(dw) \bigg]^{1/q} \bigg[  \int_{ B} p^B(t-s,y,w)p^B(t-s,w,y) \mu(dw) \bigg]^{1-1/q} \\
&\qquad \times \sup_{w\in B \cap B(y,2^{m+1} k\rho)}p^B(t-s,w,y)^{2/q-1} \nn\\
&\le  \frac{c_{19}  p^B(2t-2s, y,y)^{1-1/q}  \sup_{w\in B \cap (B(x,2^{m+1} k\rho)\setminus B(x,2^m k\rho))}p^B(t-s,w,y)^{2/q-1}  }{V(z, (2^m-2)k\rho)^{1-1/q} \,((2^m-2)k\rho)^\beta} \nn\\
&\le \frac{c_{20}  (AF_{y,\rho}(t) + \p(t))^{1/2}  \sup_{w\in B \cap B(y,2^{m+1} k\rho)}   (AF_{w,\rho}(t) + \p(t))^{1/q-1/2}  }{V(z, (2^m-2)k\rho)^{1-1/q}\, ((2^m-2)k\rho)^\beta} \nn\\ 
&\le \frac{c_{21} ( 2^{(\alpha'-\alpha_0')m} AF_{x,\rho}(t) + \p(t)  )^{1/q-1/2} \left(AF_{x,\rho}(t) + \p(t)  \right)^{1/2}}{2^{((1-1/q) \alpha_0' + \beta)m}V(z, 2k\rho)^{1-1/q} \rho^\beta} \nn\\
&\le \frac{c_{22} 2^{(\alpha'-\alpha_0')(1/q-1/2) m}  \left( AF_{x,\rho}(t) + \p(t)  \right)^{1/q}  }{2^{((1-1/q) \alpha_0' + \beta) m} V(x,\rho)^{1-1/q}  \rho^\beta }\nn\\
& \le \frac{c_{22} 2^{((1/q-1/2)\alpha' - \alpha_0'/2 - \beta) m}  \left( AF_{x,\rho}(t) + \p(t)  \right)  }{  \rho^\beta }.
\end{align*} Since $ 	\int_{B\cap B(x,k\rho)} p^{(\rho)}(s,x,z) \mu(dz) \le 1$, it follows that
for all $m\ge 2$ and $s\in (0,t/2)$, 
\begin{align*}
	f_{m}(s)&\le  \frac{c_{22} 2^{((1/q-1/2)\alpha' - \alpha_0'/2 - \beta) m}  \left( AF_{x,\rho}(t) + \p(t)  \right)  }{  \rho^\beta }.
\end{align*}
Combining this with \eqref{e:estimate-fm-1}, we obtain for all  $m\ge 2$ and $s\in (0, \eps(t \wedge (\rho^\beta/K_0)))$,
\begin{align*}
	f_{m}(s) &\le  \bigg[\frac{c_{11}\rho^{\beta/\nu-2\beta}t}{   2^{2\beta \lambda m}V(x,\rho)s^{1/\nu} } \bigg]^{\delta}  \left[ \frac{c_{22} 2^{((1/q-1/2)\alpha' - \alpha_0'/2 - \beta) m}  \left( AF_{x,\rho}(t) + \p(t)  \right)  }{  \rho^\beta } \right] ^{1-\delta}\nn\\
	&=  \frac{c_{23}\rho^{ \beta \delta/\nu - (1+\delta)\beta}  t^{ \delta}(    AF_{x,\rho}(t)  + \p(t)  )^{1-\delta}}{2^{a\beta m}V(x,\rho)^{\delta}s^{\delta/\nu} } ,
\end{align*}
where $a:=2\lambda \delta  - (1-\delta)((2/q-1)\alpha' -\alpha_0' - 2\beta  )/(2\beta)$.
Since $a>0$ by \eqref{e:NDU-range-lambda-2}  and $\delta<\nu$, repeating the arguments for \eqref{e:estimate-I4-32},  we arrive at
\begin{align*}
	I_4&=  \int_0^{\eps(t \wedge (\rho^\beta/K_0))}\sum_{m\ge 2} f_{m}(s)\,ds \le c_{24}\eps^{1-\delta/\nu}(AF_{x,\rho}(t) + \p(t)).
\end{align*}By choosing $\eps$ sufficiently small, the proof is complete.  \qed

\begin{lem}\label{l:NDU-pre}
	Under the setting of Theorem  \ref{t:NDU}, there exists $C>0$  such that for any $\rho \in \Q_+ \cap (0, R_0/K)$ and $B:=B(x_0,R)$, the heat kernel $p^B(t,x,y)$   satisfies the following estimate:
	\begin{align*}
		p^B(t,x,x) \le  CF_{x,\rho}(t) \quad \text{for all $t >0$ and $x \in B$}.
	\end{align*}
\end{lem}
\pf By Lemma \ref{l:NDU-1},  the heat kernel $p^B(t,x,y)$  exists and
	\begin{align}\label{e:l:NDU-pre-1} 
	p^B(t,x,x)\le \p_0(t)  \quad \text{for all $t>0$ and $x \in B$},
\end{align}
where
\begin{align*}
	\p_0(t):=\frac{c_1R^{\beta/\nu} }{  V(x_0,R) (t\wedge R^\beta)^{1/\nu}}  \bigg( 1\wedge \frac{R_0}{R} \bigg)^{-b'/\nu}.
\end{align*}
   Applying \eqref{e:comparison-Meyer-upper} and using  Proposition  \ref{p:esssup1},  we obtain for all $t>0$ and  $x\in B$,
 \begin{align}\label{e:l:NDU-pre-2}
	p^{B}(t,x,x) &\le \limsup_{\eps\to0} \sup_{x',y' \in B(x,\eps)} p^{(\rho)}(t,x',y')\nn\\
	&\quad +\limsup_{\eps\to0} \sup_{x',y' \in B(x,\eps)} \int_0^t \int_B   p^{(\rho)}(s,x',z)  \int_B J_2^{(\rho)}(z,dw)  p^{B}(t-s,w,y')  \mu(dz) ds\nn\\
	&=:I_1+I_2.
 \end{align} 
 By Theorem \ref{t:NDU-truncated-on} and \eqref{e:VD2}, there exists $c_2\ge 1$ such that
 \begin{align}\label{e:l:NDU-pre-3} 
 I_1\le c_2F_{x,\rho}(t).
 \end{align}
 Note that $\p_0$ is non-increasing, $\p_0(t/2)\le 2^{1/\nu}\p_0(t)$ for all $t>0$ and,  by \eqref{e:l:NDU-pre-1}, \eqref{e:l:NDU-2-1} is satisfied with $A=1$ and $\p=\p_0$. 
 Thus, applying Lemma \ref{l:NDU-2}, we obtain
	\begin{align}\label{e:l:NDU-pre-4} 
I_2&\le c_3 F_{x,\rho}(t) + 2^{-1}\left( F_{x,\rho}(t) + \p_0(t)\right).
\end{align}
Combining \eqref{e:l:NDU-pre-2}, \eqref{e:l:NDU-pre-3} and \eqref{e:l:NDU-pre-4}, we obtain for all $t>0$ and $x\in B$,
\begin{align*}
	p^B(t,x,x)\le (c_2+c_3 +2^{-1})F_{x,\rho}(t)  + 2^{-1}\p_0(t).
\end{align*} 
This means that \eqref{e:l:NDU-2-1} is satisfied with $A=c_2+c_3 +2^{-1}$ and $\p=2^{-1}\p_0$.  Applying Lemma \ref{l:NDU-2} again, we obtain
\begin{align*} 
	I_2&\le c_3 F_{x,\rho}(t) + 2^{-1} \left(  (c_2+c_3 +2^{-1})F_{x,\rho}(t) ++2^{-1}\p_0(t)\right).
\end{align*}
 Combining this with  \eqref{e:l:NDU-pre-2} and \eqref{e:l:NDU-pre-3}, we get that for all $t>0$ and $x\in B$,
 \begin{align*}
 	p^B(t,x,y) \le \left(c_2 + c_3+ 2^{-1}(c_2+c_3) + 2^{-2}\right)F_{x,\rho}(t)  + 2^{-2}\p_0(t).
 \end{align*}
 By iterating this procedure, we  deduce  that  for all $t>0$ and $x\in B$,
 \begin{align*}
p^B(t,x,x)\le ( c_2 + c_3) F_{x,\rho}(t)\sum_{n=0}^\infty 2^{-n}  	  = 2(c_2+c_3)F_{x,\rho}(t).
 \end{align*}
The proof is complete.  \qed

Now, we give the proof of Theorem \ref{t:NDU}.

\medskip

\noindent \textbf{Proof of Theorem \ref{t:NDU}.} Let $\rho \in \Q_+\cap (0,R_0/K)$  and $\{B(z_i,\rho)\}_{i\ge 1}$ be an open covering of $M$. For all $i\ge 1$, $t>0$ and  any $f\in L^2(M)$ with $\lVert f \rVert_2=1$, using H\"older's inequality in the first inequality below,  symmetry and the semigroup property in the third equality, Lemma \ref{l:NDU-pre} in the second inequality and \eqref{e:VD2} in the third inequality,  we obtain  
\begin{align*}
	\lVert P_tf \rVert_{L^\infty(B(z_i,\rho))} &= \lim_{R\to \infty}	\lVert P^{B(z_i,R)}_tf \rVert_{L^\infty(B(z_i,\rho))} = \lim_{R\to \infty}	\sup_{x \in B(z_i,\rho)} \bigg| \int_M p^{B(z_i,R)}(t,x,y)f(y)\mu(dy)\bigg|\nn\\
& \le  \lim_{R\to \infty}	 \sup_{x \in B(z_i,\rho)} \bigg( \int_M p^{B(z_i,R)}(t,x,y)^2\mu(dy)\bigg)^{1/2}  \nn\\
&=  \lim_{R\to \infty}	 \sup_{x \in B(z_i,\rho)} p^{B(z_i,R)}(2t,x,x)^{1/2}\le 	c_1 \sup_{x \in B(z_i,\rho)} F_{x,\rho}(2t)^{1/2}\le  	c_2 F_{z_i,\rho}(t)^{1/2},
\end{align*}
where $F_{x,\rho}$ is the function defined as \eqref{e:def-F-rho}. Thus, $\lVert P_t \rVert_{L^2(M)\to L^\infty(B(z_i,\rho))} \le c_2F_{z_i,\rho}(t)^{1/2}$ for all $i\ge 1$ and $t>0$. By Proposition \ref{p:existence-heat-kernel}, we deduce that $(P_t)_{t\ge0}$ has a heat kernel $p(t,x,y)$ on $(0,\infty)\times M\times M$.

 Fix $x_0\in M$. By 
\cite[Theorem 2.12(b) and (c)]{GT12}, there exists a properly exceptional set $\sN\subset M$ such that  for each  $t>0$ and $x,y \in M\setminus \sN$,   $p^{B(x_0,n)}(t,x,y)$ increases as $n\to \infty$ and converges to a  measurable function $q(t,x,y)$. By the monotone convergence theorem, it follows that for any $t>0$, non-negative $f\in L^2(M)$ and $\mu$-a.e. $x\in M$,
\begin{align*}
	P_t f(x) = \lim_{n\to \infty} P^{B(x_0,n)} f(x) = \lim_{n\to \infty} \int_M p^{B(x_0,n)}(t,x,y)f(y)\mu(dy)  = \ \int_M q(t,x,y)f(y)\mu(dy).
\end{align*}
Thus, for each fixed $t>0$, we have $p(t,x,y)=q(t,x,y)$ for $\mu$-a.e. $x,y \in M$.

Let $t>0$ and $x\in M$. Set $r_t:=t^{1/\beta}\wedge (R_0/K)$ and pick $\rho_t \in \Q_+\cap (r_t/2,r_t)$. By Lemma  \ref{l:NDU-pre} and \eqref{e:VD2}, we have 
\begin{align*}
	\esssup_{z,w \in B(x,\rho_t)}q(t,z,w) &\le \limsup_{n\to \infty} 
	\esssup_{z,w \in B(x,\rho_t)}p^{B(x_0,n)}(t,z,w) \le c_3 F_{x, \rho_t}(t) \nn\\
	& \le  \frac{c_4r_t^{\beta/\nu}}{V(x, r_t)(t\wedge r_t^\beta)^{1/\nu}} \bigg(1+ \frac{t}{r_t^\beta}\bigg)\le  \frac{c_5}{V(x, t^{1/\beta} \wedge R_0)} \bigg(1+ \frac{t}{R_0^\beta}\bigg).
\end{align*}
By Proposition \ref{p:esssup1}, it follows that
\begin{align*}
	p(t,x,x)\le 	\esssup_{z,w \in B(x,\rho_t)}p(t,z,w)  = 	\esssup_{z,w \in B(x,\rho_t)}q(t,z,w) \le  \frac{c_5}{V(x, t^{1/\beta} \wedge R_0)} \bigg(1+ \frac{t}{R_0^\beta}\bigg).
\end{align*}
The proof is complete.\qed

\section{Upper estimates for the semigroup}\label{ss:7}

In this section, we return to the case of a general scale function $\phi$. By applying a metric transformation, we obtain the following results based on those established  in Section \ref{ss:special}.

\begin{prop}\label{p:conservative-TE}
Suppose that \S \ holds. Then the following hold:

\noindent (i) $(\sE,\sF)$  is conservative. 

\noindent (ii) If \TJ \ also holds, then \TE \ holds.
\end{prop}

\begin{thm}\label{t:NDU-main}
	Suppose that \VD, \RVD,  \S \ and \Nash \ hold with $\alpha_0\ge 0$. Assume further that one of the following conditions holds:

	\medskip

	\setlength{\leftskip}{3mm}
	
	\noindent {\rm (1)} \TJ \  and \IVJ \ hold with $\gamma \ge 0 $ satisfying \eqref{e:main1-assumption}.

	\noindent {\rm (2)} \dTJ \ and \RVD \ hold with $q \in [1,\infty)$ and $\alpha_0\ge 0$ satisfying  \eqref{e:main2-assumption}.
	
	\setlength{\leftskip}{0mm}
	
	\medskip 
	
	\noindent  Then $(P_t)_{t\ge0}$ has a heat kernel $p(t,x,y)$ on $(0,\infty)\times M\times M$ and there exists   $C>0$ such that 
	\begin{align*}
		p(t,x,x) \le  \frac{C}{V(x,\phi^{-1}(x,t \wedge T_0))}\bigg(1+\frac{t}{T_0}\bigg) \quad \text{for all $t>0$ and $x \in M$}. 
	\end{align*}
	Consequently, \DUE  \ holds.
\end{thm}

\subsection{Change of metric}

The following result on a change of metric, motivated by  \cite{Ki12} and \cite{BM18}, was established in \cite[Proposition 6.1 and (6.2)]{GHH24+}.
\begin{prop}\label{p:change-of-metric}
Let $\phi$ be a scale function on $(M,d)$.	There exist constants $\beta>0$ and $C\ge 1$, and a metric $d_*$ on $M$ such that
	\begin{align*}
	C^{-1}	d_*(x,y)^\beta \le  \phi(x,d(x,y)) \le C d_*(x,y)^\beta \quad \text{for all $x,y\in M$}.
	\end{align*}
\end{prop}

Throughout the remainder of this section, we fix a scale function $\phi$ on $(M,d)$, and let $\beta$ and $d_*$ denote the  constant and the metric from Proposition \ref{p:change-of-metric}. 

For $x\in M$ and $r>0$, let $B_*(x,r):=\left\{ y \in M: d_*(x,y)<r\right\}$ denote the metric ball with respect to  $d_*$,
and let $V_*(x,r):=\mu(B_*(x,r))$.  According to \cite[Proposition 6.2(i)]{GHH24+}, there exist constants $C'\ge C^2>1$ such that
\begin{equation}\label{e:ball-compare}	B_*(x, C'^{-1}r^{1/\beta}) \subset B(x, \phi^{-1}(x, C^{-1}r)) \subset B_*(x,r^{1/\beta}) \quad \text{for all $x\in M$ and $r>0$}.
\end{equation}

\begin{lem}\label{l:volume-compare}
	Suppose that $(M,d,\mu)$ satisfies \VD. The there exists $C>1$ such that
	\begin{align*}
	C^{-1} V_*(x, r^{1/\beta})\le	V(x, \phi^{-1}(x,r)) \le C V_*(x, r^{1/\beta}) \quad \text{for all $x\in M$ and $r>0$}.
	\end{align*}
\end{lem}
 \pf Using   \eqref{e:ball-compare}, \VD \ and \eqref{e:phi-scale}, we obtain for all $x\in M$ and $r>0$,
 \begin{align*}
 	 V_*(x, r^{1/\beta}) \ge V(x,\phi^{-1}(x,c_1^{-1}r)) \ge c_2	V(x,\phi^{-1}(x,r))
 \end{align*}
 and $	 V_*(x, r^{1/\beta})  \le  V(x,\phi^{-1}(x,c_1r)) \le  c_3	V(x,\phi^{-1}(x,r))$. \qed

\begin{prop}\label{p:change-VD}
 (i) If $(M,d,\mu)$ satisfies \VD, then $(M,d_*,\mu)$ satisfies  {\rm VD}$(\alpha\beta/\beta_1)$.

\noindent (ii) If  $(M,d,\mu)$ satisfies \RVD, then $(M,d_*,\mu)$ satisfies  {\rm RVD}$(\alpha_0\beta/\beta_2)$.
\end{prop}
\pf (i) Using Lemma \ref{l:volume-compare}, \VD \ and  \eqref{e:phi-scale}, we obtain for all $x \in M$ and $R\ge r>0$,
\begin{align*}
	\frac{V_*(x,R)}{V_*(x,r)} \le 
	\frac{c_1V(x,\phi^{-1}(x,R^\beta))}{V(x,\phi^{-1}(x, r^\beta))} \le c_2 \bigg( 
	\frac{\phi^{-1}(x,R^\beta)}{\phi^{-1}(x, r^\beta)}\bigg)^{\alpha} \le c_3 \bigg( \frac{R}{r}\bigg)^{\alpha\beta/\beta_1}.
\end{align*}

\noindent (ii)  Fix $x_0\in M$ and set $R_1:=\sup_{z \in M}\phi(z,	\text{diam$(M)$})$.  Let diam$_*(M)$ denote the diameter of $M$ under $d_*$. By \eqref{e:ball-compare}, we have
\begin{align*}
	M = \overline{ B(x_0, 	\text{diam$(M)$})} \subset \overline{ B_*(x_0,  c_1 \phi(x_0,	\text{diam$(M)$})^{1/\beta})} \subset \overline{ B_*(x_0,  c_1R_1^{1/\beta})}.
\end{align*}
Thus, diam$_*(M) \le c_1 R_1^{1/\beta}$. Note that if $\diam(M)<\infty$, by   \eqref{e:phi-comp}, there exists $c_2>0$ such that $\phi(x,\diam(M)) \ge c_2 R_1$ for all $x \in M$. Thus, by \eqref{e:phi-scale}, for all $x \in M$,
\begin{align*}
	\phi^{-1}(x, c_1^\beta R_1) \le c_3 \phi^{-1}(x, c_2 R_1) = c_3 \,\diam(M).
\end{align*}
 Using   Lemma \ref{l:volume-compare},  \eqref{e:RVD2} and  \eqref{e:phi-scale}, we conclude that  for all $x \in M$ and $0<r\le R<  c_1R_1^{1/\beta}$, 
\begin{align*}
	\frac{V_*(x,R)}{V_*(x,r)} \ge 
	\frac{c_2V(x,\phi^{-1}(x,R^\beta))}{V(x,\phi^{-1}(x, r^\beta))} \ge c_3 \bigg( 
	\frac{\phi^{-1}(x,R^\beta)}{\phi^{-1}(x,r^\beta)}\bigg)^{\alpha_0} \ge c_4 \bigg( \frac{R}{r}\bigg)^{\alpha_0\beta/\beta_2}.
\end{align*}
 \qed

 \begin{prop}\label{p:change-of-metric-J}
 	Suppose that $(M,d,\mu)$ satisfies \VD.

 	\noindent (i)  If  $(\sE,\sF)$ on $(M,d,\mu)$ satisfies \TJ, then $(\sE,\sF)$ on  $(M,d_*,\mu)$ satisfies \TJb.

 	\noindent (ii) If  $(\sE,\sF)$ on $(M,d,\mu)$ satisfies \dTJ, then  $(\sE,\sF)$ on  $(M,d_*,\mu)$  satisfies \dTJb.

 	\noindent (i)   If  $(\sE,\sF)$ on $(M,d,\mu)$ satisfies \IVJ, then   $(\sE,\sF)$ on  $(M,d_*,\mu)$ satisfies \IVJb.
 \end{prop}
 \pf (i)-(ii) See \cite[Proposition 7.4(iii)]{GHH24+}.
 
 \noindent (iii)  For all $x\in M$ and $0<r\le R$, using \eqref{e:ball-compare} and Lemma \ref{l:volume-compare} in the first inequality below,  Lemma \ref{l:IVJ} in the second, and Lemma \ref{l:volume-compare} in the third, we obtain
 \begin{align*}
 	\int_{B_*(x,2R^{1/\beta})\setminus B_*(x,R^{1/\beta})} \frac{ J(x,dy)}{\sqrt{V_*(y,r^{1/\beta})}}&\le  c_1	\int_{B(x,c_2\phi^{-1}(x,R))\setminus B(x,c_3 \phi^{-1}(x,R))} \frac{ J(x,dy)}{\sqrt{V(y, \phi^{-1}(y, r))}} \\
 &\le \frac{c_2}{R\sqrt{V(x,\phi^{-1}(x,r))}} \bigg( \frac{R}{r}\bigg)^{\gamma}\le \frac{c_3}{R\sqrt{V_*(x,r^{1/\beta})}} \bigg( \frac{R}{r}\bigg)^{\gamma}.
 \end{align*}

 \qed

\begin{prop}\label{p:change-of-metric-S}
		Suppose that $(\sE,\sF)$ on $(M,d,\mu)$ satisfies \S \ with the localizing constant $T_0$. Then   $(\sE,\sF)$  on  $(M,d_*,\mu)$  satisfies \Sb \ with the localizing constant $R_0=T_0^{1/\beta}$.
\end{prop}
\pf The result follows from the proof of \cite[Proposition 7.4(ii)]{GHH24+} where   the case    $R_0= \diam_*(M)$ was considered. We provide the proof for completeness.

Let $B_*:=B_*(x_0,r)$ with $r\in (0,T_0^{1/\beta})$. By \eqref{e:ball-compare}, there exists $c_1\in (0,1)$ such that $B:=B(x_0, \phi^{-1}(x_0,c_1r^\beta))\subset B_*$.  Moreover, by \eqref{e:phi-scale} and \eqref{e:ball-compare},  we get $4^{-1}B \supset  c_2 B_*$ for some $c_2\in (0,1/4]$. By \S, there exist  $\eps_0,a_0\in (0,1)$ such that for all $t \in (0, a_0c_1^\beta r^\beta]$,	\begin{align*}
\essinf_{c_2B_*}	P^{B_*}_t \1_{B_*} \ge \essinf_{	4^{-1} B }	P^{	 B}_t \1_{	 B} \ge \eps_0.
\end{align*}
By Lemma \ref{l:S-self-improvement}, this implies the desired result. \qed

\begin{prop}\label{p:change-of-metric-N}
		Suppose that $(M,d,\mu)$ satisfies \VD \ and  $(\sE,\sF)$ on $(M,d,\mu)$ satisfies \Nash \ with the localizing constant $T_0$. Then   $(\sE,\sF)$ on  $(M,d_*,\mu)$  satisfies \Nashb \ with $b'=b\beta$ and the localizing constant $R_0=T_0^{1/\beta}$.
\end{prop}
\pf Let $B_*:=B_*(x_0,r)$  and $f\in \sF^{B_*}$. By \eqref{e:ball-compare}, there exists $c_1>1$ such that $B:=B(x_0, \phi^{-1}(x_0, c_1r^\beta))\supset B_*$.   Using \Nash \ for $f\in \sF^B$ in the first inequality below, \VD \ and \eqref{e:phi-scale} in the second and Lemma \ref{l:volume-compare} in the third,  we obtain
\begin{align*}
		\lVert f \rVert_{2}^{2+2\nu} &\le   \frac{c_2c_1^\beta r^\beta}{   V(x_0,\phi^{-1}(x_0, c_1r^\beta))^\nu}  \bigg( 1\wedge \frac{T_0^{1/\beta}}{c_1 r} \bigg)^{-b \beta} \left[ \sE(f,f) +  \frac{\lVert f \rVert_2^2}{c_1r^\beta}\right] \lVert f \rVert_1^{2\nu}\nn\\
		&\le   \frac{c_3 r^\beta}{   V(x_0,\phi^{-1}(x_0, r^{\beta}))^\nu}  \bigg( 1\wedge \frac{T_0^{1/\beta}}{r} \bigg)^{-b \beta} \left[ \sE(f,f) +  \frac{\lVert f \rVert_2^2}{r^\beta}\right] \lVert f \rVert_1^{2\nu}\nn\\
			&\le   \frac{c_4 r^\beta}{   V_*(x_0,r))^\nu}  \bigg( 1\wedge \frac{T_0^{1/\beta}}{r} \bigg)^{-b \beta} \left[ \sE(f,f) +  \frac{\lVert f \rVert_2^2}{r^\beta}\right] \lVert f \rVert_1^{2\nu}.
\end{align*} \qed 

\subsection{Proofs of Proposition  \ref{p:conservative-TE} and Theorem \ref{t:NDU-main}} Let $\beta$ and $d_*$ be the  constant and the metric from Proposition \ref{p:change-of-metric}.

\medskip

\noindent \textbf{Proof of Proposition \ref{p:conservative-TE}.} By Proposition  \ref{p:change-of-metric-S},  as a form on $(M,d_*,\mu)$, $(\sE,\sF)$ satisfies  \Sb \ with the localizing constant $R_0=T_0^{1/\beta}$.

\noindent (i) The result follows from Lemma \ref{l:conservative}.

\noindent (ii) By Proposition \ref{p:change-of-metric-J}(i), the form   $(\sE,\sF)$  on $(M,d_*,\mu)$ satisfies \TJb. Applying Proposition \ref{p:TEb}, we obtain for any ball $B_*:=B_*(x_0,R)$,
	\begin{align}\label{e:TE-main-1}
	\esssup_{y \in 4^{-1}B_*} P_{t} \1_{B_*^c}(y)\le   \frac{c_1t}{R^\beta \wedge T_0} \quad \text{for all $t>0$}.
\end{align}
 By \eqref{e:ball-compare},  there exists $c_2>1$ such that for any $x_0 \in M$ and $r>0$,
\begin{align}\label{e:TE-main-2}
B_*(x_0,c_2^{-1}\phi (x_0,r)^{1/\beta}) \subset 	B\subset B_*(x_0,c_2\phi (x_0,r)^{1/\beta}).
\end{align}
Set $\eta:=1/(4c_2^2)$. By  \eqref{e:TE-main-1} and \eqref{e:TE-main-2},  for any ball $B:=B(x_0,r)$  and  $t>0$,
\begin{align}\label{e:TE-main-3}
	\esssup_{y \in \eta B} P_{t} \1_{B^c}(y) &\le \esssup_{y \in B_*(x_0,4^{-1}c_2^{-1}\phi(x_0,r)^{1/\beta})} P_{t} \1_{B_*(x_0,c_2^{-1}\phi(x_0,r)^{1/\beta})^c}(y)\le \frac{c_1c_2^{\beta} t}{\phi(x_0,r) \wedge T_0 }.
\end{align}

Now, let us prove \eqref{e:TE}. Let $B:=B(x_0,r)$  and $\{B(z_i,\eta r/4 )\}_{i\ge 1}$ be a covering of $4^{-1}B$  where $z_i \in \overline{4^{-1}B}$ for all $i \ge 1$.  Using \eqref{e:TE-main-3} and \eqref{e:phi-scale}, we obtain for all $t>0$,
\begin{align*}
		\esssup_{y \in 4^{-1} B} P_{t} \1_{B^c}(y) &\le \sup_{i\ge 1} 	\esssup_{y \in  B(z_i,\eta r/4)} P_{t} \1_{B^c}(y) \le \sup_{i\ge 1} 	\esssup_{y \in  B(z_i,\eta r/4)} P_{t} \1_{B(z_i,r/4)^c}(y)\\
		&\le \frac{c_1c_2^{\beta} t}{\phi(x_0,r/4) \wedge T_0 } \le \frac{c_3 t}{\phi(x_0,r) \wedge T_0 }.
\end{align*} 
The proof is complete.\qed

\noindent \textbf{Proof of Theorem \ref{t:NDU-main}.} By Propositions \ref{p:change-VD}, \ref{p:change-of-metric-S} and \ref{p:change-of-metric-N}, $(M,d_*,\mu)$ satisfies \VDb \ and \RVDb \ with $\alpha'=\alpha\beta/\beta_1$ and $\alpha_0'=\alpha_0\beta/\beta_2$, and as a form on $(M,d_*,\mu)$, $(\sE,\sF)$ satisfies  \Sb \ and \Nashb \  with  $b'=b\beta$ and the localizing constant $R_0=T_0^{1/\beta}$. Moreover, by Proposition \ref{p:change-of-metric-J},  one of the additional assumptions in Theorem \ref{t:NDU} is satisfied. Therefore, by Theorem \ref{t:NDU},  the semigroup $(P_t)_{t\ge0}$ has a heat kernel $p(t,x,y)$ on $(0,\infty)\times M\times M$ and there exists   $c_1>0$ such that 
\begin{align*}
	p(t,x,x) \le  \frac{c_1}{V_*(x,(t\wedge T_0)^{1/\beta})}\bigg(1+\frac{t}{T_0}\bigg) \quad \text{for all $t>0$ and $x \in M$}. 
\end{align*}
By Lemma \ref{l:volume-compare}, this yields the desired result. \qed

\section{Proofs of main results}\label{ss:proofs}

\subsection{Proofs of Theorem \ref{t:main1}, Corollary \ref{c:main1} and  Theorem \ref{t:main2}}

\

Under \VD \ and \TJ, we have the following implications:
\begin{align}\label{e:implications}
	\text{\wFK \ $+$ \CS} \;\,&\text{$\;\implies\;$ \GFK \ $+$ \Gl   \quad  (by Lemma \ref{l:GFK} and  Proposition \ref{p:G})},\nn\\
	\text{\GFK} \;\,&\text{$\ifandonlyif\,$ \Nash     \qquad \qquad \qquad\,   (by Lemma \ref{l:Nash})},\nn\\
	\text{\Gl} \;\,&\text{$\;\implies\;$ \S    \qquad \qquad \qquad \qquad \  (by Lemma \ref{l:G->S})}.
\end{align}

\noindent \textbf{Proof of Theorem \ref{t:main1}.} (i)  follows from Lemmas \ref{l:DUE-GFK}, \ref{l:TE->S} and \ref{l:S->CS}, and Proposition \ref{p:conservative-TE}. Moreover,  (ii) is a consequence of  \eqref{e:implications} and  Theorem \ref{t:NDU-main}. By Lemma \ref{l:wFK}, the additional assertion holds. The proof is complete. \qed

\noindent \textbf{Proof of Corollary \ref{c:main1}.} 
 (c) clearly implies (a), and (b) implies (c) by Theorem \ref{t:main1}(i).  To establish the implication (a) $\Rightarrow$ (b), we observe that if either  hypothesis (1) or (3) holds, then \IVJ  \ holds with $\gamma=0$, and that  if  (2) holds, then \IVJ \ holds with $\gamma\le 1$ by Lemma \ref{p:IVJ-general}. Thus, in all cases,  \IVJ \ holds with $\gamma\le 1$,  ensuring that \eqref{e:main1-assumption}  is satisfied for any $\nu>0$. Consequently, by Theorem \ref{t:main1}(ii) and the fact that \TJ \ trivially holds if $(\sE,\sF)$ is strongly local,  we conclude that (a) implies (b).   This completes the proof. \qed

\noindent \textbf{Proof of Theorem \ref{t:main2}.} Recall that   \TJ \ holds under \VD \ and \dTJ. Thus, \eqref{e:main2-result} follows from \eqref{e:implications} and Theorem \ref{t:NDU-main}. Since \eqref{e:main2-assumption} is true for any $\nu>0$ if $q=2$, \eqref{e:main2-result-2} follows. \qed

\subsection{Proof of Theorem \ref{t:main1-counterexample}}

In this subsection, we prove  Theorem \ref{t:main1-counterexample} by explicitly constructing a counterexample.

 For a fixed constant $\xi \in (0,1)$,  define $I_0:=[0,1]$ and construct $I_i$ for $i\ge 1$ by removing a centrally situated open subinterval of relative length $\xi$ from each interval in $I_{i-1}$. The $\xi$-Cantor set $\sC_{\xi}$ is then defined as $\sC_{\xi}:=\cap_{i\ge 0}I_i$. Note that the case $\xi=1/3$ corresponds to the standard Cantor set. For each $i \ge 0$, let $\mu_i$ be the Lebesgue measure restricted to $I_i$,  normalized such that $\mu_i(I_i)=1$. Then $\mu_i$ converges weakly to a measure $\mu$ supported  on $\sC_\xi$. Moreover,  $(\sC_\xi, |\cdot|, \mu)$ is a metric measure space satisfying
 \begin{align}\label{e:alpha-set}
 C^{-1}r^{\alpha_\xi} \le \mu(	B(x,r)) \le C r^{\alpha_\xi} \quad \text{for all $x\in \sC_\xi$ and $0<r\le 1$},
 \end{align}
 for some $C>1$, where
 \begin{align}\label{e:alpha-xi}
 	\alpha_\xi:=\frac{\log 2}{\log(2/(1-\xi))}.
 \end{align}
 For $n\ge 1$, denote by $\sC_\xi^n = \sC_\xi \times \cdots \times \sC_\xi$  the $n$-fold  product of $\sC_\xi$ with itself, $d^n$ the sup metric on $\sC_{\xi}^n$ defined as $d^n((x_1,\cdots, x_n), (y_1,\cdots,y_n))= \max_{1\le i\le n} |x_i-y_i|$, and $\mu^n$ the product measure on $\sC_{\xi}^n$. Note that \diam$(C^n_\xi)=1$ and by \eqref{e:alpha-set}, 
  \begin{align}\label{e:alpha-n-set}
 	C^{-1}r^{n\alpha_\xi} \le \mu^n(	B(x,r)) \le C r^{n\alpha_\xi} \quad \text{for all $x\in \sC_\xi^n$ and $0<r\le 1$},
 \end{align}
 for some $C>1$. In particular, $(\sC_\xi^n, d^n, \mu^n)$ satisfies VD($n\alpha_\xi$) and RVD($n\alpha_\xi$).

 Let $n\ge 1$ and $0<\beta_1\le \beta_2<2$ be fixed constants, and let $\bbeta_n:\sC^n_\xi \to [\beta_1,\beta_2]$ be a function satisfying $|\bbeta_n(x)-\bbeta_n(y)| \le d^n(x,y)$ for all $x,y \in \sC^n_\xi$. Define
\begin{align}\label{e:def-phi-bbeta}
	\phi_{\bbeta_n}(x,r)= \begin{cases}
		r^{\bbeta_n(x)}, &  0<r\le 1;\\
		r^{\beta_1}, & r>1.
	\end{cases} 
\end{align}
Then $\phi_{\bbeta_n}$ is a scale function on $\sC^n_\xi$ satisfying  \eqref{e:phi-scale} with $\beta_1$ and $\beta_2$ (see \cite[p.36-37]{CKW-elp}). Associated with  $\phi_{\bbeta_n}$, we define a transition jump kernel $J_{\bbeta_n}(x,dy)$ on $\sC^n_\xi$ as
\begin{align*}
	J_{\bbeta_n}(x,dy)  = 	J_{\bbeta_n}((x_1,\cdots, x_n),d(y_1,\cdots, y_n))= \sum_{i=1}^n |x_i - y_i|^{-\alpha_\xi - \bbeta_n(x) \wedge \bbeta_n(y)} \mu(dy_i) \prod_{j\neq i} \delta_{x_j} (dy_j),
\end{align*}
where $\delta_{x_j}$ denotes the Dirac measure at $x_j$. Denote by Lip$(\sC^n_\xi)$ the family of all Lipschitz continuous functions on $\sC^n_\xi$. Define a bilinear form $(\sE^{\bbeta_n},\sF^{\bbeta_n})$ on $L^2(\sC^n_\xi)$ by
 \begin{align}\label{e:def-sE-bbeta}
 	\sE^{\bbeta_n} (f,g) &= \int_{\sC^n_\xi \times \sC^n_\xi} (f(x)-f(y))(g(x)-g(y)) J_{\bbeta_n}(x,dy) \, \mu^n(dx),\nn\\
 	\sF^{\bbeta_n}&= \text{$(\sE^{\bbeta_n}_1)^{1/2}$-closure of Lip$(\sC^n_\xi)$}.
 \end{align}
 Since $\beta_2<2$, it is straightforward to verify that   $\sE^{\bbeta_n}(f,f)<\infty$ for all $f\in \text{Lip$(\sC^n_\xi)$}$ and that $(\sE^{\bbeta_n},\sF^{\bbeta_n})$ constitutes a regular Dirichlet form on $L^2(\sC^n_\xi)$. 

 \begin{lem}\label{l:bbeta-form}
 	 $(\sE^{\bbeta_n},\sF^{\bbeta_n})$ satisfies {\rm TJ}$(\phi_{\bbeta_n})$ and {\rm CS}$(\phi_{\bbeta_n})$ with $T_0=1$.
 \end{lem}
 \pf Let $x=(x_1,\cdots,x_n) \in \sC^n_\xi$.  Since \diam$(\sC^n_\xi)=1$, we have $J_{\bbeta_n}(x,B(x,r)^c)=0$ for all $r>1$. Further, for all  $r\in (0,1]$,
 \begin{align}\label{e:bbeta-TJ}
 	J_{\bbeta_n}(x,B(x,r)^c) \le \sum_{i=1}^n \int_{y_i \in \sC_\xi,\, |x_i-y_i|\ge r} | x_i-y_i|^{-\alpha_\xi- \bbeta_n(x)} \mu(dy_i).
 \end{align}
 For each $1\le i\le n$, using \eqref{e:alpha-set} and $\bbeta_n(x) \ge \beta_1$, we get
 \begin{align*}
 	\int_{y_i \in \sC_\xi,\, |x_i-y_i|\ge r} | x_i-y_i|^{-\alpha_\xi- \bbeta_n(x)} \mu(dy_i) &\le \sum_{n\ge 0, \, 2^n r\le 1}  (2^nr)^{-\alpha_\xi- \bbeta_n(x)}  \int_{y_i \in \sC_\xi,\, 2^n r \le |x_i-y_i|<2^{n+1}r }\mu(dy_i)\\
 &\le \frac{c_1}{r^{\bbeta_n(x)}}\sum_{n\ge 0} 2^{-n \bbeta_n(x)}  \le  \frac{c_1}{r^{\bbeta_n(x)}}\sum_{n\ge 0} 2^{-n \beta_1} = \frac{c_2}{\phi_{\bbeta_n}(x,r)}.
 \end{align*}
 Combining this with \eqref{e:bbeta-TJ}, we deduce that   {\rm TJ}$(\phi_{\bbeta_n})$ holds. Since $\beta_2<2$,  by Proposition \ref{l:beta<2}, {\rm CS}$(\phi_{\bbeta_n})$ also holds. The proof is complete.\qed
 
\begin{lem}\label{l:bbeta-constant}
Suppose that $\bbeta_n(x)=\beta_1$ for all $x\in \sC^n_\xi$. Then the following hold.

\noindent (i) The heat kernel $p_{\bbeta_n}(t,x,y)$ for   $(\sE^{\bbeta_n}, \sF^{\bbeta_n})$ exists. Moreover,  for any $T>0$,  there exists a constant $C=C(T)>0$ such that
\begin{equation}\label{e:axis-NDU}
	p_{\bbeta_n}(t,x,y) \le \frac{C}{\phi_{\bbeta_n}^{-1}(x_0,t)^{n\alpha_\xi}} \quad \text{for all $ t\in (0,T]$ and $\mu^n$-a.e. $x,y \in  \sC^n_\xi$}.
\end{equation}

\noindent (ii) For any ball $B\subset \sC^n_\xi$, the heat kernel $p_{\bbeta_n}^B(t,x,y)$ for the part  of $(\sE^{\bbeta_n}, \sF^{\bbeta_n})$ on $B$ exists. Moreover, there exist constants  $C>0$ and $\eta\in (0,1)$ such that for any ball $B:=B(x_0,r)$ with $r\in (0,1]$,
\begin{equation}\label{e:axis-NDL}
	p_{\bbeta_n}^B(t,x,y) \ge \frac{C}{\phi_{\bbeta_n}^{-1}(x_0,t)^{n\alpha_\xi}} \;\; \text{for all $ t\in (0, \eta\phi_{\bbeta_n}(x_0,r)]$ and $\mu^n$-a.e. $x,y \in B(x_0, \eta \phi_{\bbeta_n}^{-1}(x_0,t))$}.
\end{equation}

\noindent (iii)  $(\sE^{\bbeta_n},\sF^{\bbeta_n})$ satisfies {\rm FK}$_{\nu}(\phi_{\bbeta_n})$ with $\nu=\beta_1/n\alpha_\xi$ and $T_0=1$.

\noindent (iv) (Poincar\'e inequality) There exist constants $C>0$ and $\delta \in (0,1]$ such that for any ball $B:=B(x_0,r)$ with $r\in (0,1]$ and  any bounded $f \in \sF^{\bbeta_n}$,
\begin{align}\label{e:PI}
	\inf_{a\in \R}\int_{\delta B} \left|f- a\right|^2 d\mu^n \le C \phi_{\bbeta_n}(x_0,r)\int_{B\times B} (f(x)-f(y))^2 J_{\bbeta_n}(x,dy)\, \mu^n(dx).
\end{align} 

\end{lem}
\pf By assumption,  $\phi_{\bbeta_n}(x,r)=r^{\beta_1}$ for all $x\in \sC^n_\xi$ and $r>0$.

\noindent (i)-(ii) For $n=1$,   the results follow from  \eqref{e:alpha-set} and \cite[Propositions 2.2(i) and  2.3]{CKK09}. For $n\ge 2$, observe that 
$$p_{\bbeta_n}(t,(x_1,\cdots,x_n),(y_1,\cdots,y_n))= \prod_{i=1}^n p_{\bbeta_1}(t,x_i,y_i) \quad \text{for $\mu^n$-a.e. $x,y\in \sC^n_\xi$}$$
and for any ball $B:=B(x_0,r)\subset \sC^n_\xi$ with $r\in (0,1]$, 
$$p^B_{\bbeta_n}(t,(x_1,\cdots,x_n),(y_1,\cdots,y_n))= \prod_{i=1}^n p^{B((x_0)_i, r)}_{\bbeta_1}(t,x_i,y_i) \quad \text{for $\mu^n$-a.e. $x,y\in B$}. $$
Using these, since \eqref{e:axis-NDU} and \eqref{e:axis-NDL} hold for $n=1$, we obtain that   \eqref{e:axis-NDU} and \eqref{e:axis-NDL} hold for all $n\ge 1$.

\noindent (iii) By \eqref{e:alpha-set}, and Lemmas \ref{l:DUE-GFK} and \ref{l:wFK}, the result follows from (i). 

\noindent (iv) By (ii), the result follows from  \cite[Proposition 3.5(i)]{CKW-jems}. Although \cite{CKW-jems} assumes  the underlying space to be unbounded,  the proof  also holds  for  bounded spaces (see \cite[Remark 1.19]{CKW-jems}). \qed

\begin{lem}\label{l:bbeta-general}  
	The assertions (ii)-(iv) in Lemma \ref{l:bbeta-constant} hold without  the  assumption that $\bbeta_n(x)=\beta_1$  for all $x\in \sC^n_\xi$ .
\end{lem}
\pf  Fix $x_0 \in \sC^n_\xi$ and $r\in (0,1]$. Let $\beta^*:=\inf_{z \in B(x_0,r)} \bbeta_n(z)$ and define $\bbeta^*_n(z):=\beta^*$ for all $z \in \sC^n_\xi$. By \eqref{e:phi-comp}, there exists $c_1>0$ independent of $x_0$ and $r$ such that
\begin{align}\label{e:comp-bbeta-general}  
 \phi_{\bbeta_n}(x_0,r)   \ge c_1 \sup_{z\in B(x_0,r)} \phi_{\bbeta_n}(z,r) =  c_1 r^{\beta^*}.
\end{align}

  Let $D$ be a non-empty open subset of $B(x_0,r)$ and $f\in \sF^{\bbeta_n,D}$ with $\lVert f \rVert_2=1$, where $\sF^{\bbeta_n,D}$ is the $(\sE^{\bbeta_n}_1)^{1/2}$-closure of $\sF^{\bbeta_n}\cap C_c(D)$.  
Using \diam$(\sC_\xi)=1$,  Lemma \ref{l:bbeta-constant}(iii) and  \eqref{e:comp-bbeta-general}, we get
\begin{align*}
	\sE^{\bbeta_n}(f,f) \ge\sE^{\bbeta_n^*}(f,f) \ge  \frac{c_2}{r^{\beta^*}} \left( \frac{\mu^n(B(x_0,r))}{\mu^n(D)}\right)^{\beta^*/n\alpha_\xi} \ge  \frac{c_3}{\phi_{\bbeta_n}(x_0,r)} \left( \frac{\mu^n(B(x_0,r))}{\mu^n(D)}\right)^{\beta_1/n\alpha_\xi}.
\end{align*}
Hence, $(\sE^{\bbeta_n},\sF^{\bbeta_n})$ satisfies {\rm FK}$_{\beta_1/n\alpha_\xi}(\phi_{\bbeta_n})$.

For \eqref{e:PI}, using  \diam$(\sC_\xi)=1$, \eqref{e:comp-bbeta-general}  and  Lemma \ref{l:bbeta-constant}(iv) (with $\beta_1$ replaced by $\beta^*$), we see that for any  $f \in $ Lip$(\sC^n_\xi)$,
\begin{align*}
	  &\phi_{\bbeta_n}(x_0,r)\int_{B\times B} (f(x)-f(y))^2 J_{\bbeta_n}(x,dy)\, \mu^n(dx) \\
	  &\ge  c_1r^{\beta^*}\int_{B\times B} (f(x)-f(y))^2 J_{\bbeta^*_n}(x,dy)\, \mu^n(dx)\ge c_4 \inf_{a\in \R}\int_{\delta B} \left|f- a\right|^2 d\mu^n .
\end{align*}
Thus, since Lip$(\sC^n_\xi)$ is dense in $\sF^{\beta_n}$ with $\sE^{\bbeta_n}_1$-norm, \eqref{e:PI} holds.

 By \eqref{e:alpha-set}, Lemma \ref{l:bbeta-form} and \eqref{e:PI}, the conditions (VD), (RVD), (TJ), (PI) and (ABB) in \cite{GHH23} are satisfied. Therefore, by \cite[Theorem 2.10]{GHH23},  assertion (ii) of Lemma \ref{l:bbeta-constant} holds (see condition (LLE) therein). 
The proof is complete. \qed 

Now we present the proof of Theorem \ref{t:main1-counterexample}.

\medskip

\noindent \textbf{Proof of Theorem \ref{t:main1-counterexample}.}  Let $\eps>0$ be fixed. Choose $\xi \in (0,1)$ such that $\alpha_\xi\le \eps/2$, where $\alpha_\xi$ is defined in \eqref{e:alpha-xi}, and   $n\ge 1$ such that
\begin{align}\label{e:def-eps}
\frac{2(1+\eps)}{n \alpha_\xi} < \frac12  \quad \text{ and } \quad \bigg(1 + \frac{\eps}{2} \bigg) \bigg(	1 - \frac{2(1+\eps)}{n \alpha_\xi} \bigg) >1.
\end{align}
Denote ${\bf 0}=(0,\cdots,0) \in \sC_\xi^n$ and  ${\bf e}_1=(1,0,\cdots,0) \in \sC_\xi^n$.  Set 
\begin{align*}
	\beta_1:=1 \quad \text{and} \quad \beta_2:= \bigg( 1- \frac{2(1+\eps)}{n \alpha_\xi}\bigg)^{-1} \in (1,2),
\end{align*}
and let $\bbeta_n:\sC^n_\xi \to [\beta_1,\beta_2]$ be such that $|\bbeta_n(x)-\bbeta_n(y)| \le d^n(x,y)$ for all $x,y \in \sC^n_\xi$, $\bbeta_n(x)=\beta_2$ for all  $x \in B({\bf 0}, 1/4)$ and  $\bbeta_n(x)=\beta_1$ for all  $x \in B({\bf e}_1, 1/4)$.  Define $\phi_{\bbeta_n}$ as \eqref{e:def-phi-bbeta}. Observe that
\begin{align}\label{e:def-eps-2}
\frac{n\alpha_\xi}{2} - \frac{n\alpha_\xi}{2\beta_2}= \frac{n\alpha_\xi}{2} \left[ 1 -\bigg( 1- \frac{2(1+\eps)}{n \alpha_\xi}\bigg) \right]  = 1+\eps.
\end{align}
Since $(\sC_\xi^n, d^n, \mu^n)$ satisfies VD($n\alpha_\xi$) and RVD($n\alpha_\xi$), by Lemma  \ref{l:bbeta-form} and \ref{l:bbeta-general}, and Proposition \ref{p:IVJ-general} with \eqref{e:def-eps-2}, the regular Dirichlet  form $(\sE^{\bbeta_n},\sF^{\bbeta_n})$ on $L^2(\sC^n_\xi)$ defined as \eqref{e:def-sE-bbeta} satisfies  {\rm TJ}$(\phi_{\bbeta_n})$, {\rm WFK}$_{\nu}(\phi_{\bbeta_n})$, {\rm CS}$(\phi_{\bbeta_n})$ and {\rm IJ}$_{2,\gamma}(\phi_{\bbeta_n})$ with $\nu=1/n\alpha_\xi$ and $\gamma=1+\eps$.  Note that $(1-\nu)\gamma < 1+ \nu + \eps$.
 
 We  show that  $(\sE^{\bbeta_n},\sF^{\bbeta_n})$ does not  satisfy  DUE$(\phi_{\bbeta_n})$. 
 Let $(\sE^{\bbeta_n,(1/4)}, \sF^{\bbeta_n})$ be the  $1/4$-truncated form of $(\sE^{\bbeta_n}, \sF^{\bbeta_n})$. By Lemmas \ref{l:bbeta-form} and \ref{l:truncated-general}, $(\sE^{\bbeta_n,(1/4)}, \sF^{\bbeta_n})$ is a regular Dirichlet form on $L^2(\sC^n_\xi)$. Note that $(\sE^{\bbeta_n}, \sF^{\bbeta_n})$ satisfies    {\rm Nash}$_{\nu}(\phi_{\bbeta_n})$ with $\nu=1/n\alpha_\xi$   by Lemmas \ref{l:GFK} and \ref{l:Nash}. Thus, by Lemma \ref{l:Nash-truncated-consequence}(ii) and  Proposition \ref{p:existence-heat-kernel}, for any open subset $D$ of $\sC^n_\xi$,  the heat kernel $p^{(1/4),D}_{\bbeta_n}(t,x,y)$ for the part of $(\sE^{\bbeta_n,(1/4)}, \sF^{\bbeta_n})$ on $D$ exists. In particular, the heat kernel $p^{(1/4)}_{\bbeta_n}(t,x,y)$ for $(\sE^{\bbeta_n, (1/4)}, \sF^{\bbeta_n})$ exists.
    By Lemma \ref{l:bbeta-general}, there exist constants $\eta\in (0,1)$ and $c_1>0$ such that for all $t\in (0,\eta]$,
\begin{align}
	\essinf_{x,z\in B({\bf 0}, \eta t^{1/\beta_2})} p^{B({\bf 0}, (t/\eta)^{1/\beta_2})}_{\bbeta_n} (t,x,z) &\ge c_1t^{-n\alpha_\xi/\beta_2},\label{e:counter-NDL-1}\\
	\essinf_{w,y\in B({\bf e}_1, \eta t)} p^{B({\bf e}_1, t/\eta)}_{\bbeta_n} (t,w,y) &\ge c_1t^{-n\alpha_\xi}.\label{e:counter-NDL-2}
\end{align} For all $t\in (0,8^{-\beta_2}\eta]$ and $z \in B({\bf 0},\eta t^{1/\beta_2})$, we have $\overline{B({\bf 0}, (t/\eta)^{1/\beta_2})} \subset B(z, 1/4)$. Thus, for all $s>0$, $w \in \sC^n_\xi \setminus B(z,1/4)$ and $y \in \sC^n_\xi$, it holds that $p^{B({\bf 0}, (t/\eta)^{1/\beta_2})}_{\bbeta_n} (s,w,y) =0$.  Consequently, by \eqref{e:comparison-Meyer-upper} and \eqref{e:counter-NDL-1}, we obtain for all $t \in (0, 8^{-\beta_2}\eta]$,
\begin{equation}\label{e:counter-NDL-3}
		\essinf_{x,z\in B({\bf 0}, \eta t^{1/\beta_2})} p^{(1/4), B({\bf 0}, (t/\eta)^{1/\beta_2})}_{\bbeta_n} (t,x,z) \ge 	\essinf_{x,z\in B({\bf 0}, \eta t^{1/\beta_2})} p^{B({\bf 0}, (t/\eta)^{1/\beta_2})}_{\bbeta_n} (t,x,z)\ge c_1t^{-n\alpha_\xi/\beta_2} .
\end{equation}
Observe that for all $r\in (0,1/4]$,
\begin{align}\label{e:counter-NDL-4}
	&\int_{B({\bf 0},\eta r)} \int_{B({\bf e}_1,\eta  r)} \1_{\{d^n(z,w) \ge 1/4\}}J_{\bbeta_n}(z,dw)\, \mu^n(dz) \nn\\
	&= 	\int_{B({\bf 0},\eta r)} \int_{w_1 \in \sC_\xi,\, |w_1 -1| < \eta r } \frac{\mu(dw_1)}{|z_1- w_1|^{\alpha_\xi +1}}\,\mu^n(dz) \nn\\
	&\ge 	\int_{B({\bf 0},\eta r)} \int_{w_1 \in \sC_\xi,\, |w_1 -1| < \eta r }  \mu(dw_1)\, \mu^n(dz) \ge c_2  (\eta r)^{(n+1)\alpha_\xi},
\end{align}
where we used \eqref{e:alpha-set}  and \eqref{e:alpha-n-set} in the last inequality. For all $t\in (0,\eta]$, $\mu^n$-a.e. $x \in B({\bf 0}, \eta t/8)$ and $\mu^n$-a.e. $  y\in B({\bf e}_1, \eta t/8)$, using \eqref{e:comparison-Meyer-lower} in the first inequality below, \eqref{e:counter-NDL-2} and \eqref{e:counter-NDL-3}  in the second, and \eqref{e:counter-NDL-3}  and  \eqref{e:counter-NDL-4}  in the third, we obtain 
\begin{align}\label{e:counter-NDL-5}
	&	p_{\bbeta_n}(t,x,y)\ge   \int_{0}^{(t/8)^{\beta_2}} \int_{B({\bf 0}, \eta s^{1/\beta_2})} p_{\bbeta_n}^{(1/4)}(s, x, z) \nn\\
		&\qquad \qquad \qquad\qquad \qquad  \times \int_{B({\bf e}_1, \eta s^{1/\beta_2})} \1_{\{d^n(z,w) \ge 1/4\}}J_{\bbeta_n}(z,dw) p_{\bbeta_n}(t-s,w,y) \mu^n(dz) ds\nn\\
		&\ge  c_1^2  \int_{0}^{(t/4)^{\beta_2}}   s^{-n\alpha_\xi/\beta_2}  (7t/8)^{-n\alpha_\xi} \int_{B({\bf 0}, \eta s^{1/\beta_2})}  \int_{B({\bf e}_1, \eta s^{1/\beta_2})} \1_{\{d^n(z,w) \ge 1/4\}}J_{\bbeta_n}(z,dw)  \mu^n(dz) ds\nn\\
		&\ge c_3 t^{-n\alpha_\xi}  \int_{0}^{(t/4)^{\beta_2}}   s^{\alpha_\xi/\beta_2}   ds = c_4 t^{-(n-1)\alpha_\xi + \beta_2}.
\end{align}
 
If  $(\sE^{\bbeta_n},\sF^{\bbeta_n})$ also satisfies  DUE$(\phi_{\bbeta_n})$, then by \eqref{e:DUE-2}, we get for all $t\in (0,8^{-\beta_2}\eta]$,
\begin{align}\label{e:main1-contradiction}
	\sup_{x\in B({\bf 0}, \eta t/8), \, y\in B({\bf e}_1, \eta t/8)}	p_{\bbeta_n}(t,x,y)\le c_5 t^{-(1+1/\beta_2)n\alpha_\xi/2}.
\end{align}
Recall  $\alpha_\xi \le \eps/2$. Since $\beta_2<1+\eps/2$ by the second inequality in  \eqref{e:def-eps}, using \eqref{e:def-eps-2}, we get
\begin{align*}
	(n-1)\alpha_\xi - \beta_2 -\frac{n\alpha_\xi}{2}\bigg( 1 + \frac{1}{\beta_2}\bigg)  = \frac{n\alpha_\xi}{2} - \frac{n\alpha_\xi}{2\beta_2} - \alpha_\xi - \beta_2 = 1+\eps - \alpha_\xi - \beta_2 >0,
\end{align*}
implying that $\lim_{t\to 0} t^{-(n-1)\alpha_\xi + \beta_2}/  t^{-(1+1/\beta_2)n\alpha_\xi/2} = \infty$. Thus, by \eqref{e:counter-NDL-5}, \eqref{e:main1-contradiction} leads to a contradiction. Consequently,  DUE$(\phi_{\bbeta_n})$ fails for $(\sE^{\bbeta_n},\sF^{\bbeta_n})$. The proof is complete. \qed

\noindent \textbf{Proof of Theorem \ref{t:main1-special-counterexample}.} By applying the change of metric from Proposition \ref{p:change-of-metric} to the example constructed in the proof of Theorem \ref{t:main1-counterexample}, and using Lemma \ref{l:volume-compare} and  Propositions \ref{p:change-VD}, \ref{p:change-of-metric-J}, \ref{p:change-of-metric-S} and \ref{p:change-of-metric-N}, we obtain the desired result. The details are omitted. \qed

\section{Applications}\label{ss:9}

\subsection{Non-local operators with variable orders}
In this subsection, we apply our main theorems to symmetric non-local operators on $\R^d$, whose jump kernel is bounded below by that of a stable-like operator of variable order. This example is motivated by \cite[Example 2.3]{BKK10} and \cite[Section 6]{CKW-elp}. Unlike  these references, we also consider the case where the jump density does not exist, and we establish a sharp on-diagonal upper heat kernel estimate.

	Let $d\ge 1$, $0<\beta_1\le \beta_2<2$ and  $\bbeta:\R^d \to [\beta_1,\beta_2]$ be a function satisfying
\begin{align*}
	|\bbeta(x) - \bbeta(y)| \le \frac{c_0}{\log (2/|x-y|)} \quad \text{for $|x-y|<1$}
\end{align*}
with some $c_0>0$. Define
\begin{align*}
	\phi_{\bbeta}(x,r)= \begin{cases}
		r^{\bbeta(x)}, &  0<r\le 1;\\
		r^{\beta_1}, & r>1.
	\end{cases} 
\end{align*}
Then $\phi_{\bbeta}$ is a scale function on $\R^d$ satisfying  \eqref{e:phi-scale} with $\beta_1$ and $\beta_2$ (see \cite[p.36-37]{CKW-elp}). 

Let $\sE$ be a bilinear form defined by
\begin{align*}
	 \sE(f,g) = \int_{\R^d\times \R^d} (f(x)-f(y))(g(x)-g(y))  J(dx,dy),
\end{align*}
where $ J(dx,dy)$ is a symmetric Radon measure on $\R^d\times \R^d\setminus \diag$.  Assume that $ J(dx,dy) =  J(x,dy)dx$ satisfies TJ$(\phi_\bbeta)$. Since $\beta_2<2$, by an  argument similar to that in \eqref{e:beta<2}, we have $ \sE(f,f)<\infty$ for all $f \in $ Lip$_c(\R^d)$, where Lip$_c(\R^d)$ denotes the family of all Lipschitz  functions on $\R^d$ with compact supports.  Let $ \sF$ be the $ \sE_1$-closure of Lip$_c(\R^d)$. Then $( \sE,  \sF)$ is a regular Dirichlet form on $L^2(\R^d)$. We consider the following assumptions on $J$:

\medskip

\noindent {\bf (A)} Suppose that there exist constants $R_0\in (0,\infty]$, $C>1$ and a Radon measure $j$ on $\R^d\setminus \{0\}$ such that the following properties  hold:

\vspace{-0.2in}

\begin{align*}
	(\rm A1)& \;	J(x,dy)  \ge   (|x-y| \wedge 1)^{\beta_1- \bbeta(x)} j(x-dy) \;\; \text{ for all $x\in \R^d$}.\\
	(\rm A2)& \; C^{-1} r^{-\beta_1}\le  \sup_{|z|\le 1/r}\int_{\R^d\setminus \{0\}} (1- \cos \la z,y \ra) j(dy) \le Cr^{-\beta_1} \;\;\text{ for all $r \in (0,R_0)$}. \hspace{0.6in}
\end{align*}

 A typical example of $j(dy)$ satisfying (A2) is given by $j(dy)= c|y|^{-d-\beta_1}dy$ for some $c>0$. In this case, (A1) takes the form that, for all $x\in \R^d$, $J(x,dy)=J(x,y)dy$,  where $J(x,y)$ is a measurable function on $\R^d\times \R^d \setminus \diag$ satisfying
 \begin{align*}
 	J(x,y) \ge  \frac{c}{|x-y|^{d} \phi_\bbeta(x,|x-y|)}=c \begin{cases}
 		|x-y|^{-d-\bbeta(x)}  &\mbox{ for $|x-y|< 1$},\\
 		|x-y|^{-d-\beta_1} &\mbox{ for $|x-y|\ge 1$}.
 	\end{cases}
 \end{align*} 
 
\begin{thm}\label{t:perturb}
	Suppose that  $( \sE, \sF)$ is defined as above, with the jump kernel $J(dx,dy)=J(x,dy)dx$ satisfying {\rm TJ}$(\phi_\bbeta)$ and {\bf (A)}.

 \noindent (i) If \begin{align}\label{e:TJ-perturb-example}
 (d-\beta_1)\bigg[ \frac{1}{\beta_1} - \frac{1}{\beta_2}\bigg]  < 2+ \frac{2\beta_1}{d},
	\end{align}
then $( \sE, \sF)$ satisfies {\rm DUE}$(\phi_\bbeta)$ with the localizing constant $R_0^{\beta_1}$.

\noindent (ii) Suppose additionally that $J(dx,dy)=J(x,y)dxdy$ for a measurable function $J$ on $\R^d\times \R^d\setminus \diag$ such that 
\begin{align}\label{e:TJ-theta-example}
\bigg( 	\int_{\R^d} J(x,y)^{q} dy \bigg)^{1/q} \le \frac{C}{r^{(1-1/q)d} \phi_\bbeta(x,r)}  \quad \text{for all $x\in \R^d$ and $r>0$},
\end{align}
with some constants $C>0$ and  $q \in [1,\infty)$ satisfying
\begin{align}\label{e:perturb-2}
 (d-\beta_1)\bigg[ \bigg(\frac{2}{q}-1\bigg) \frac{1}{\beta_1} - \frac{1}{\beta_2}\bigg]  < 2+ \frac{2\beta_1}{d}.
\end{align}
Then $( \sE, \sF)$ satisfies  {\rm DUE}$(\phi_\bbeta)$ with the localizing constant $R_0^{\beta_1}$.
\end{thm}

\begin{remark}
	By Remark  \ref{r:dTJ}, if there exists a constant $C>0$ such that
	\begin{align}\label{e:A-upper}
		J(x,y) \le \frac{C}{|x-y|^d \phi_\bbeta(x,|x-y|)} \quad \text{for all $x,y \in \R^d$},
	\end{align}
	then TJ$(\phi_\bbeta)$ and \eqref{e:TJ-theta-example} hold with $q=2$.   \eqref{e:perturb-2} clearly holds for $q=2$. Thus, when the jump density $J(x,y)$ exists and satisfies both {\bf (A)} and \eqref{e:A-upper},   by Theorem \ref{t:perturb}, $( \sE, \sF)$ satisfies  {\rm DUE}$(\phi_\bbeta)$. 
\end{remark}

To prove Theorem \ref{t:perturb}, we first establish the following lemma.

\begin{lem}\label{l:perturb}
If  \TJ \ and {\bf (A)} hold, then  {\rm FK}$_\nu(\phi_\bbeta)$ holds with $\nu=\beta_1/d$.
\end{lem}
\pf Define $j_s(dx)=2^{-1}j(dx) + 2^{-1}j(-dx)$ and $j_s(dx,dy)= j_s(x-dy)dx$. Then $j_s(dx,dy)$ ia a symmetric Radon measure on $\R^d\times \R^d\setminus \diag$. Further, by the symmetry of $J(dx,dy)$, (A1) and the fact that  $\bbeta(x) \ge \beta_1$ for all $x\in \R^d$, we have
\begin{align}\label{e:perturb-below}
	j_s(dx,dy)&\le  (|x-y| \wedge 1)^{\bbeta(x)-\beta_1}   J(dx,dy) \le J(dx,dy)\quad \text{ in $\R^d\times \R^d \setminus \diag$}.
\end{align}
Define a bilinear form $(\wt \sE,\wt \sF)$ on $L^2(\R^d)$ as follows:
\begin{align*}
	\wt \sE(f,g):=\int_{\R^d\times \R^d} (f(x)-f(y))(g(x)-g(y))  j_s(dx,dy)
\end{align*}
and $\wt \sF$ is the $\wt \sE_1$-closure of Lip$_c(\R^d)$. By \eqref{e:perturb-below},  $(\wt \sE,\wt \sF)$ is a well-defined, regular Dirichlet form on $L^2(\R^d)$ satisfying TJ$(\phi_\bbeta)$. The generator $\wt \sL$ of $(\wt \sE,\wt\sF)$ is  a L\'evy operator associated with a symmetric L\'evy process $Y=(Y_t)_{t\ge 0}$ on $\R^d$. Moreover, the characteristic exponent $\Psi$ of $Y$, defined by $\Psi(x):=-\log \E e^{i \la x, Y_1 \ra}$, is given by $	\Psi(x) = \int_{\R^d} \left( 1 - e^{i\la x, z \ra} - i \la x,z \ra \1_{\{|z|<1\}} \right) j_s(dz).$ By (A2), it follows from \cite[(1), and Theorems 3.1 and 3.12]{GS20} that $(\wt \sE,\wt \sF)$ satisfies DUE$(r^{\beta_1})$ with the localizing constant $R_0^{\beta_1}$.  Consequently, by Theorem \ref{t:main1}(i) and Lemma \ref{l:wFK}, $(\wt \sE, \wt \sF)$ satisfies {\rm FK}$_\nu(\phi_\bbeta)$ for $\nu = \beta_1/d$, with some $\delta_1\in (0,1)$ and the localizing constant $R_0^{\beta_1}$.

Let  $B(x_0,r)\subset \R^d$ with $r\in (0, \delta_1^{1/\beta_1} R_0 )$ and $D$ be a non-empty open subset of $B(x_0,r)$. Let $f \in \sF^D$ with $\lVert f \rVert_2=1$. Since $(\wt \sE, \wt \sF)$ satisfies {\rm FK}$_\nu(\phi_\bbeta)$ for $\nu = \beta_1/d$, we have
\begin{align}\label{e:perturb-FK}
	\wt \sE(f,f) \ge  c_1\mu(D)^{-\beta_1/d}.
\end{align}
If $r\ge1$, then by \eqref{e:perturb-below} and \eqref{e:perturb-FK}, we obtain $\sE(f,f) \ge \wt \sE(f,f)\ge  c_1\mu(D)^{-\beta_1/d}.$ Further, if $r<1$, then  using \eqref{e:perturb-below}, the fact that $D\subset B(x_0,r)$, \eqref{e:phi-comp} for $\phi_\bbeta$, and \eqref{e:perturb-FK}, we obtain
\begin{align*}
	\sE(f,f) \ge \inf_{x \in B(x_0,r)}(2r)^{\beta_1- \bbeta(x)} \wt \sE(f,f) \ge c_2  r^{\beta_1- \bbeta(x_0)} \mu(D)^{-\beta_1/d} = \frac{c_3}{r^{\bbeta(x_0)}} \bigg( \frac{|B(x_0,r)|}{\mu(D)}\bigg)^{\beta_1/d}.
\end{align*}
By \eqref{e:def-lambda1}, we arrive at the desired result. \qed

\medskip

\noindent\textbf{Proof of Theorem \ref{t:perturb}.} Set $\gamma:=d/(2\beta_1)- d/(2\beta_2)$ and $\nu:=\beta_1/d$.  A direct computation show that  \eqref{e:main1-assumption} holds for these parameters if and only if 
\eqref{e:TJ-perturb-example} is satisfied. Similarly, 
\eqref{e:main2-assumption} holds for $\alpha=\alpha_0=d$ and $\nu$ if and only if \eqref{e:perturb-2} is satisfied.
  By Proposition \ref{p:IVJ-general}, and   Lemmas \ref{l:beta<2} and  \ref{l:perturb}, the form $(\sE,\sF)$ satisfies IJ$_{2,\gamma}(\phi_\bbeta)$, CS$(\phi_\bbeta)$ and FK$_\nu(\phi_\bbeta)$.  Consequently, the desired results follow immediately from  Theorems \ref{t:main1} and \ref{t:main2}. \qed

\subsection{Singular jump kernel} 

Let  $T_0\in (0,\infty]$ and $n\ge 2$ be fixed,  and let $\phi: (0,\infty)\to (0,\infty)$  be a strictly increasing,  continuous function satisfying the following scaling condition:
\begin{align*}
	C^{-1}\left( \frac{R}{r} \right)^{\beta_1} \le 		\frac{\phi(R)}{\phi(r)} \le C\left( \frac{R}{r} \right)^{\beta_2} \quad \text{for all} \;\, 0<r\le R,
\end{align*}
for some  $C\ge 1$ and $0<\beta_1\le \beta_2$. Note that  for any metric  space $(M,d)$, the function $\phi$ induces  a natural scale function on $M$, defined by $\phi(x,r)=\phi(r)$. 

For $1\le i\le n$, let $(M^{(i)}, d^{(i)},\mu^{(i)})$ be a metric measure space in which all metric balls are precompact. We assume that each $(M^{(i)}, d^{(i)},\mu^{(i)})$ satisfies VD$( \alpha^{(i)})$ for some $\alpha^{(i)}>0$. Let $(M,d,\mu)$ be the Cartesian product of $(M^{(i)}, d^{(i)},\mu^{(i)})$ for $1\le i\le n$, equipped with the sup metric and the product measure. For each $i$, let $B^{(i)}(x_i,r)$ denote the metric ball in $M^{(i)}$ centered at $x_i \in M^{(i)}$ of radius $r$ and $V^{(i)}(x_i,r):=\mu^{(i)}(B^{(i)}(x_i,r))$. Then for all $x=(x_1,\cdots, x_n)\in M$ and $r>0$,
\begin{align}\label{e:product-ball}
	B(x,r) = \prod_{i=1}^n B^{(i)}(x_i,r) \quad \text{and} \quad 	V(x,r) = \prod_{i=1}^n V^{(i)}(x_i,r).
\end{align}
Thus, $(M,d,\mu)$ is a metric measure space satisfying \VD  \ with $\alpha=	 \sum_{i=1}^n \alpha^{(i)}$.

 For each $i$, let  $(\sE^{(i)},\sF^{(i)})$ be a  regular Dirichlet form on $L^2(M^{(i)}; \mu^{(i)})$.  We  assume the following hypothesis for $(\sE^{(i)},\sF^{(i)})$:  Let $T_0\in (0,\infty]$ be a fixed constant.

\medskip

\noindent {\bf (B)} For all $1\le i\le n$, the form $(\sE^{(i)},\sF^{(i)})$ satisfies the following conditions:

\vspace{-0.2in}

\begin{align}
	(\rm B1)& \; \text{CS($\phi$)  and  FK$_{\nu_i}(\phi)$ hold for some $\nu_i>0$,  with the localizing constant $T_0$.} \nn\\
	(\rm B2)& \; \text{TJ($\phi$) and  IJ$_{2,\gamma_i}(\phi)$ hold for some $\gamma_i\ge 0$ satisfying} \hspace{2.1in}\nn\\
	&  \qquad \qquad \qquad	 (1-\overline \nu_i) \gamma_i < 1+ \overline  \nu_i, \quad \text{where} \;\; \overline \nu_i: = \nu_i \wedge \frac{\beta_1}{\sum_{i=1}^n \alpha^{(i)}} . \label{e:B2} 
\end{align}

Let  $\sE^{(i),(L)}$ denote the strongly local part  of $(\sE^{(i)},\sF^{(i)})$, and   $J^{(i)}(dx_i,dy_i)=J^{(i)}(x_i,dy_i)\mu^{(i)}(dx_i)$ denote
the  jump measure of $(\sE^{(i)},\sF^{(i)})$.
Define a symmetric  bilinear  form $ 	\wh \sE^{(L)}$ on $L^2(M,\mu)$ as
\begin{align*}
	\wh \sE^{(L)}(f,f):=\sum_{i=1}^n \int_{\wt M^{(i)}} \sE^{(i),(L)}(f(x_1, \dots, x_{i-1}, \cdot, x_{i+1}, \dots, x_n),f(x_1, \dots, x_{i-1}, \cdot, x_{i+1}, \dots, x_n))&\\
	\mu^{(i)}(dx_1) \cdots \mu^{(i-1)}(dx_{i-1})\, \mu^{(i+1)} (dx_{i+1}) \cdots \mu^{(n)}(dx_n),&
\end{align*}
where $\wt M^{(i)}:= M^{(1)} \times \dots M^{(i-1)} \times M^{(i+1)} \times \cdots M^{(n)}$, and a symmetric Radon measure $\wh J$ on $M\times M\setminus \diag$ as
\begin{align*}
	\wh J(dx,dy)&=	\wh J(x,dy)\mu(dx):= \sum_{i=1}^n  J^{(i)}(x_i, dy_i) \prod_{j\neq i} \delta_{x_j} (dy_j) \mu(dx).
\end{align*}

Let $(\sE,\sF)$ be a regular Dirichlet form on $L^2(M,\mu)$ with the  Beurling-Deny decomposition:
\begin{align*}
	\sE(f,f) = \sE^{(L)}(f,f) + \int_{M\times M} (f(x)-f(y))^2J(dx,dy),
\end{align*}
where $\sE^{(L)}$ is the strongly local part of $(\sE,\sF)$ and $J(dx,dy)$ is a symmetric Radon measure on $M\times M\setminus \diag$. Consider the following hypothesis for $(\sE,\sF)$:
 
\medskip

\noindent {\bf (C)} There exists a constant $C>1$ such that the following properties hold:

\vspace{-0.2in}

\begin{align*}
	(\rm C1)& \;	C^{-1} \wh \sE^{(L)}(f,f)\le \sE^{(L)}(f,f) \le C\wh \sE^{(L)}(f,f) \quad \text{for all $f  \in L^2(M,\mu)$},\\
	(\rm C2)& \; \text{The kernel $J$ is of the form $J(dx,dy)=J(x,dy) \mu(dx)$, where $J(x,dy)$ satisfies}\hspace{0.2in}\nn\\
	&  \qquad \qquad \quad	 C^{-1}\wh J(x,dy)\le   J(x,dy) \le 	C^{-1}\wh J(x,dy) \;\; \text{ for all $x \in M$.} \hspace{0.75in}
\end{align*}

Typical examples of $( \sE, \sF)$ satisfying {\bf (C)} include the  form associated with the \textit{cylindrical fractional Laplacian} 
$\sum_{i=1}^n -(-\partial_{ii})^{\beta/2}$ on $\R^n$ for some $\beta \in (0,2)$, and those involving the Laplacian   $\Delta + \sum_{i=1}^n -(-\partial_{ii})^{\beta/2}$. Specifically, for $\sum_{i=1}^n -(-\partial_{ii})^{\beta/2}$, {\bf (B)} is satisfied with $T_0=\infty$, $\phi(r)=r^\beta$, $\nu_1 = \cdots= \nu_n=\beta$ and $\gamma_1 = \cdots= \gamma_n=0$, and  for   $\Delta + \sum_{i=1}^n -(-\partial_{ii})^{\beta/2}$,   {\bf (B)} is satisfied with $T_0=1$,  $\phi(r)=r^2$, $\nu_1 = \cdots= \nu_n=2$ and $\gamma_1 = \cdots= \gamma_n=0$.  The cylindrical fractional Laplacian is a fundamental example of a  non-local operator  that does not possess a jump density. For a basic understanding to this operator in Euclidean spaces, we refer to \cite{BC06, BC10, Xu13, KKK22}, and for general metric measure spaces, we refer to \cite{IRS13, St05}.

As an application of our main result, we obtain the following theorem. 
\begin{thm}\label{t:singular}
Let $(M,d,\mu)$ be a product space and $(\sE,\sF)$ be a regular Dirichlet form on $L^2(M,\mu)$ defined as above. If {\bf (B)} and {\bf (C)} hold, then $(\sE,\sF)$ satisfies  \DUE.
\end{thm}
\pf By {\bf (B)} and Theorem \ref{t:main1}, for each $i$, the semigroup $(P^{(i)}_t)_{t\ge 0}$ associated with  $\sE^{(i)}$ is conservative, satisfies \TE \ and has a heat kernel $p^{(i)}(t,x_i,y_i)$ on $(0,\infty) \times M^{(i)} \times M^{(i)}$ satisfying \DUE. 
  Define a Dirichlet form $(\wh\sE, \sF)$ on $L^2(M,\mu)$ by
\begin{align*}
	\wh \sE(f,f)= \sE^{(L)}(f,f) + \int_{M\times M} (f(x)-f(y))^2 \wh J(dx,dy),
\end{align*}
and let $(P_t)_{t\ge 0}$ be the  semigroup associated with  $(\wh \sE, \sF)$. Observe that
\begin{align*}
	P_t f = (P^{(1)}_t f_1, \cdots, P^{(n)}_t f_n) \quad \text{
		for all $t>0$ and  $f = (f_1, \cdots, f_n) \in L^2(M,\mu)$.}  
\end{align*}
Hence, $(P_t)_{t\ge 0}$ is conservative and has a heat kernel $p(t,x,y)$ on $(0,\infty)\times M\times M$ given by
\begin{align}\label{e:product-heatkernel}
	p(t,(x_1,\cdots,x_n),(y_1,\cdots,y_n))= \prod_{i=1}^n p^{(i)}(t,x_i,y_i) .
\end{align}
By  \eqref{e:product-heatkernel} and \eqref{e:product-ball},  we see that $(\wh\sE, \sF)$ satisfies \DUE. Further, for any  $B(x,r) = \prod_{i=1}^n B^{(i)}(x_i,r)$  and $t>0$, since all $(P^{(i)}_t)_{t\ge 0}$ satisfy \TE, we obtain
\begin{align*}
		\esssup_{ B(x,r/4)}	P_t \1_{B^c} &= 	\esssup_{ B(x,r/4)}\left(1 - 	P_t \1_{B}\right) = 	\esssup_{ B(x,r/4)}\bigg(1 - \prod_{i=1}^n	P^{(i)}_t \1_{B^{(i)}(x_i,r)}\bigg) \\
		&\le	\esssup_{ B(x,r/4)}\bigg(1 -  \left(1 - \frac{c_1t}{\phi(r) \wedge T_0} \right)_+^n \bigg) \le \frac{nc_1t}{\phi(r) \wedge T_0}.
\end{align*}
We used the inequality $(1-r)^n \ge 1-nr$ for all $r\in [0,1]$ in the last inequality. Thus, $(P_t)_{t\ge 0}$ satisfies  \TE. Recall that  $(M,d,\mu)$ satisfies \VD \ with $\alpha:=	 \sum_{i=1}^n \alpha^{(i)}$. Applying Theorem \ref{t:main1}(i), it follows that $(\wh \sE,\sF)$ satisfies \wFK \ with $\nu=\beta_1/\alpha$ and  \CS.

We claim that $(\wh \sE,\sF)$ satisfies \IVJ \ with $\gamma:=\max_{1\le i\le n} \gamma_i$, where $\gamma_i$, $1\le i\le n$, are the constants in {\bf (B)}. Indeed, for all $x=(x_1,\cdots, x_n)\in M$ and $0<r\le R$, using \eqref{e:product-ball} in the equality below, and  IJ$_{2,\gamma_i}(\phi)$ for $(\sE^{(i)},\sF^{(i)})$ in the first inequality,  we get
\begin{align*}
	&\int_{B(x,2\phi^{-1}(R))\setminus B(x,\phi^{-1}(R))} \frac{ \wh J(x,dy)}{\sqrt{V(y,\phi^{-1}(r))}}\\
	&=\sum_{i=1}^n \frac{1}{\prod_{j\neq i} \sqrt{ V^{(j)}(x_j,\phi^{-1}(r))}}\int_{B(x_i,2\phi^{-1}(R))\setminus B(x_i,\phi^{-1}(R))} \frac{ J^{(i)}(x_i,dy_i)}{\sqrt{V^{(i)}(y_i,\phi^{-1}(r)) }} \\
		&\le  \frac{c_2}{R\sqrt{V(x,\phi^{-1}(x,r))}}\sum_{i=1}^n \bigg( \frac{R}{r}\bigg)^{\gamma_i} \le  \frac{nc_2}{R\sqrt{V(x,\phi^{-1}(x,r))}} \bigg( \frac{R}{r}\bigg)^{\gamma} .
\end{align*}

By {\bf (C)}, we deduce that $(\sE,\sF)$ satisfies  \wFK \ with $\nu=\beta_1/\alpha$,  \CS, and  \IVJ \ with $\gamma$. From \eqref{e:B2}, we have $(1-\nu)\gamma <1+\nu$. Consequently, applying Theorem \ref{t:main1}, we  obtain the desired conclusion. \qed

\appendix

\section{Appendix}

Throughout the appendix,  let  $(\X,d_\X,m)$ be a metric measure space,  and $(\sG,\sH)$ be a regular Dirichlet form on $L^2(\X,m)$. For $p \in [1,\infty]$ and a Borel set $U\subset \X$, we write $L^p(U):=L^p(U,m)$.

\subsection{Pointwise existence and regularity of heat kernels}

Let   $(Q_t)_{t\ge 0}$ be the semigroup associated with $(\sG,\sH)$. The following was established in \cite{GHH21}.

\begin{prop}\label{p:existence-heat-kernel}   Let $p\in [1,2]$. Suppose that there exist  $T\in (0,\infty]$, a countable  open covering $\sU$ of  $\X$, and  $\Phi:\sU\times (0,T)\to (0,\infty)$ such that for any $U \in \sU$ and $f \in L^p(\X) \cap L^2(\X)$, 
	\begin{align*}
		\lVert Q_t f \rVert_{L^\infty(U)} \le \Phi(U,t) \lVert f \rVert_{L^p(\X)} \quad \text{ for all $t\in (0,T)$}.
	\end{align*}
	Then $(Q_t)_{t\ge 0}$ has a pointwise defined heat kernel $q(t,x,y)$ on $(0,\infty) \times \X\times \X$. Moreover,  for each $t\in (0,T)$, $U \in \sU$ and $x \in U$,
	\begin{align*}
	\lVert  q(t,x,\cdot) \rVert_{L^{p/(p-1)}(\X)}   \le \Phi(U,t) \quad &\mbox{ if $p>1$},\\
		\sup_{y\in \X} q(t,x,y)  \le \Phi(U,t) \quad &\mbox{ if $p=1$}.
	\end{align*}
\end{prop}
\pf The result follows from \cite[Theorem 2.1 and Corollary 4.2]{GHH21} for $p=1$ and \cite[Theorem 2.2]{GHH21} for $p\in (1,2]$.\qed

The following regular property of heat kernels was established in \cite{GT12}. Although \cite{GT12} treats only strongly local Dirichlet forms, the same proof applies to general Dirichlet forms as well.

\begin{prop}\label{p:esssup1} 
	Let  $U,V\subset \X$ be  open sets.  Suppose that $(Q_t)_{t\ge 0}$ has a pointwise defined heat kernel $q(t,x,y)$. Then for  any $t>0$ and Borel $f \in L^1(\X)$, we have
	\begin{align*}
		\esssup_{x\in U} \int_\X q(t,x,y)f(y)m(dy) &= 	\sup_{x\in U}\int_\X q(t,x,y)f(y)m(dy)
	\end{align*}and
	\begin{align*}
		\esssup_{x\in U,\, y\in V} q(t,x,y) &= 	\sup_{x\in U,\, y\in V} q(t,x,y).
	\end{align*}
\end{prop}
\pf The results follow from \cite[Corollaries 2.7 and 2.9]{GT12} and \eqref{e:heatkernel-domain}. \qed

\subsection{Some results for truncated Dirichlet forms}\label{ss:truncation}

In this subsection, we assume that $(\sG,\sH)$ has no killing part. 
By the Beurling-Deny decomposition, we have \begin{align*}
	\sG(f,g)=\sG^{(L)}(f,g) + \sG^{(J)}(f,g), \quad f,g \in \sH,
\end{align*}
where $\sG^{(L)}$ is the strongly local part of $(\sG,\sH)$ and 
\begin{align*}
	\sG^{(J)}(f,g)&= \int_{\X \times \X} (f(x)-f(y))(g(x)-g(y))J(dx,dy),
\end{align*}
for  a symmetric Radon measure $J$ on $\X\times \X \setminus \diag$.  We further assume that there is a transition jump kernel $J(x,dy)$ on $\X \times \sB(\X)$, where $\sB(\X)$ denotes the $\sigma$-algebra of Borel sets of $\X$, such that
\begin{align*}
	J(dx,dy) = J(x,dy) m(dx) \quad \text{in $\X \times \X$}.
\end{align*} 
For $\rho>0$ and  $x\in \X$, let
\begin{align*}
	J_1^{(\rho)}(x,dy) :=  \1_{\{d_\X(x,y) < \rho\}} J(x,dy) \quad \text{ and } \quad J_2^{(\rho)}(x,dy) = J(x,dy)- J_1^{(\rho)}(x,dy).
\end{align*} 
Then we define a bilinear form $(\sG^{(\rho)},\sH)$ on $L^2(\X)$ by
\begin{align*}
	\sG^{(\rho)}(f,g) = \sG^{(L)}(f,g) + \int_{\X\times \X} (f(x)-f(y))(g(x)-g(y))  J_1^{(\rho)}(x,dy)m(dx).
\end{align*}
 Note that $(\sG^{(\rho)},\sH)$ is   a Dirichlet form on $L^2(\X)$. Whenever $(\sG^{(\rho)},\sH)$ is  regular, for each $f\in \sH$, we  use a $\sG^{(\rho)}$-quasi-continuous  version  of it.

\begin{lem}\label{l:truncated-general}
	Let $\rho>0$.	Suppose that
	\begin{align}\label{e:truncated-general}
		\sJ(\rho):=	\esssup_{x\in \X} J_2^{(\rho)}(x,\X)<\infty.
	\end{align}
	Then we have
	\begin{align*}
		\sG(f,f) - \sG^{(\rho)}(f,f)\le  4 \sJ(\rho)\lVert f \rVert_{L^2(\X)}^2 \quad \text{for all $f \in \sH$}.
	\end{align*}
	Consequently, the norms   $\sG_1^{1/2}$ and $(\sG^{(\rho)}_1)^{1/2}$ are equivalent on $\sH$, and  $(\sG^{(\rho)}, \sH)$ is regular.
\end{lem}
\pf  For any $f \in \sH$, using the symmetry of $J_2^{(\rho)}$ and \eqref{e:truncated-general}, we obtain
\begin{align*}
	\sG(f,f) - \sG^{(\rho)}(f,f) &\le 4 \int_\X f(x)^2  J^{(\rho)}_2(x,\X) \,m(dx) \le   4\sJ(\rho)\lVert f \rVert_{L^2(\X)}^2.
\end{align*} 
It follows that $\sG_1^{(\rho)}(f,f)\le \sG_1(f,f) \le  (1+4\sJ(\rho)) \sG_1^{(\rho)}(f,f)$ for all $f\in \sH$. Thus,  $\sG_1^{1/2}$ and $(\sG^{(\rho)}_1)^{1/2}$ are equivalent on $\sH$. Since $(\sG,\sH)$  is  regular, we deduce  that  $(\sG^{(\rho)}, \sH)$ is also regular.  \qed

For  an open subset $D$ of $\X$ and $\rho>0$,    let $(Q^D_t)_{t\ge 0}$, $(Q^{(\rho)}_t)_{t\ge 0}$  and  $(Q^{(\rho),D}_t)_{t\ge 0}$ be the semigroups associated with   $(\sG, \sH^D)$, $(\sG^{(\rho)},\sH)$ and $(\sG^{(\rho)},\sH^{(\rho),D})$, respectively, where $
\sH^{D} $ is the $\sG_1^{1/2}$-closure of $ \sH\cap C_c(D)$ and 
$\sH^{(\rho),D} $ is the $(\sG^{(\rho)}_1)^{1/2}$-closure of $ \sH\cap C_c(D)$. If \eqref{e:truncated-general} holds, then  by Lemma \ref{l:truncated-general}, we have $\sH^D = \sH^{(\rho),D}$ for any $\rho>0$ and open set $D\subset \X$.

The following two results were   proved in \cite[Proposition 4.6]{GHL14} and \cite[Lemma 3.1(b) and (3.5)]{BGK09} respectively under the existence of the density  for $K(x,dy)$.  The proofs remain valid under the current setting.

\begin{prop}\label{p:comparison-tail}
	Let $0<\rho<\rho'$.	Suppose that  \eqref{e:truncated-general} holds. Then for any  open set $D\subset \X$  and   $0\le f \in L^\infty(D)$, we have that, for all $t>0$ and $m$-a.e. $x \in D$,
	\begin{align}\label{e:comparison-tail-1}
		\big|	Q^{(\rho'),D}_t f(x) -Q^{(\rho),D}_t f(x) \big|\le  2t \lVert f \rVert_{L^\infty(D)}	\esssup_{z\in \X}  \left( J_2^{(\rho)}(z,\X) - J_2^{(\rho')}(z,\X) \right)
	\end{align}
	and
	\begin{align}\label{e:comparison-tail-2}
		\big|	Q^D_t f(x) -Q^{(\rho),D}_t f(x) \big|\le  2t \lVert f \rVert_{L^\infty(D)} \sJ(\rho).
	\end{align}
\end{prop}

\begin{prop}\label{p:comparison-Meyer}
	Let $\rho>0$ and $D$ be an open subset of $\X$. 	Suppose that  \eqref{e:truncated-general} holds, and the semigroups $(Q^D_t)_{t\ge 0}$ and $(Q^{(\rho),D}_t)_{t\ge 0}$ possess heat kernels $q^D(t,x,y)$ and $q^{(\rho),D}(t,x,y)$ on $(0,\infty) \times D\times D$, respectively. Then for all $t>0$ and $m$-a.e. $x,y \in D$,
	\begin{equation}\label{e:comparison-Meyer-upper}
	q^D(t,x,y)\le q^{(\rho),D}(t,x,y) +  \int_0^t \int_D q^{(\rho), D}(s, x, z) \int_D J_2^{(\rho)}(z, dw)  q^D(t-s,w,y)\, m(dz) ds
	\end{equation}
	and
		\begin{equation}\label{e:comparison-Meyer-lower}
		q^D(t,x,y) \ge   \int_0^t \int_D q^{(\rho), D}(s, x, z) \int_D J_2^{(\rho)}(z, dw)  q^D(t-s,w,y)\, m(dz) ds.
	\end{equation}
\end{prop}

For an open set $U\subset \X$ and $\rho>0$, the $\rho$\textit{-neighborhood} of $U$ is defined by
\begin{align*}
	U_\rho:=\big\{ x \in \X: \text{dist}(x,U)<\rho\big\}.
\end{align*}

\begin{prop}\label{p:comparison-entrance}
	Let  $\rho>0$ and $U,D$ be open subsets of $\X$ such that $U\subset D$ and $U_{\rho+\delta}$ is precompact for some $\delta>0$.	Suppose that $(\sG^{(\rho)},\sH)$ is a regular Dirichlet form on $L^2(\X)$. Then for any  subset $K$ of $U$ which is closed relative to the topology on $D$, we have for all $0\le f \in  L^\infty(D)$, $t>0$ and $m$-a.e. $x \in  U_\rho$, 
	\begin{align}\label{e:comparison-entrance}
		Q^{(\rho),D}_t f(x) \le Q^{(\rho),U}_t f(x) + \left( 1 - Q^{(\rho),U}_t \1_U(x) \right) \sup_{0<s \le t} \lVert 	Q^{(\rho),D}_s f\rVert_{L^\infty(U_{\rho+\delta} \setminus K)}.
	\end{align}
\end{prop}
A result similar to  Proposition \ref{p:comparison-entrance} was proved in  \cite[Corollary 4.8]{GHL10} under the assumption that   $K$ is compact. In contrast, here we require only that  $K$ be closed in the subspace topology on $D$.  Moreover,   by  using $U_{\rho+\delta}$ instead of $U_\rho$     in the right-hand side of \eqref{e:comparison-entrance}, we  do not need  $U_\rho\subset D$, unlike  in \cite{GHL10}. These modifications are essential for our application.  The proof of  Proposition \ref{p:comparison-entrance} follows the approach of \cite{GHL10} but involves non-trivial adaptations to account for the weaker topological condition on  $K$.

Before proving Proposition \ref{p:comparison-entrance}, we briefly recall the notion of  weak derivatives and the  parabolic maximum principle for weak subsolutions.

	Let $U$ be an open subset of $\X$ and $I$ be an open interval. A function $u:I\to L^2(U)$ is said to be \textit{weakly differentiable} at $t\in I$ if there is a function $w \in L^2(U)$ such that
	\begin{align*}
		\lim_{\eps\to 0} \frac1{\eps} \la u(t+\eps)-u(t), \vp \ra = \la w, \vp\ra \quad \text{for all $\vp \in L^2(U)$}.
	\end{align*}
	The function $w$ is called the \textit{weak derivative} of $u$ at $t$ and we denote $w=\partial_t u$.

\begin{prop}\label{p:maximum}
{\rm{  \!\!\! (\!\!\cite[Proposition 5.2]{GHL09}).}} 	Let $(\sE,\sF)$ be a regular Dirichlet form on $L^2(\X)$, $T \in (0,\infty]$ and $U$ be an open subset of $\X$. Suppose that $u:(0,T) \times \X\to \R$ satisfies the following properties:

	\smallskip

	$\bullet$ $u_+(t,\cdot)$ converges to $0$ in $L^2(U)$ as $t \to 0$.

	$\bullet$ $t\mapsto u(t)|_{U}$ is weakly differentiable in $L^2(U)$ at any $t \in (0,T)$;

	$\bullet$ $u_+(t,\cdot) \in \sF^U$ for all $t \in (0,T)$;

	$\bullet$ For any  $0\le \psi \in \sF^{U}$, we have $\la u_t, \psi\ra + \sE(u,\psi)\le 0$ for all $t \in (0,T)$,

	\smallskip
	
\noindent where  $ \sF^U$ denotes  the closure of  $\sF \cap C_c(U)$ in $\sE_1$-norm.   Then 
\begin{align*}
	u(t,x)\le 0 \quad \text{for all $t\in (0,T)$ and $m$-a.e. $x \in \X$.} 
\end{align*}
\end{prop}

\noindent \textbf{Proof of Proposition \ref{p:comparison-entrance}.} 
Let $0\le f\in  L^\infty(D)$. Fix $\eps\in (0,\delta)$ and let $\vp \in \sH$ be a cutoff function for  $U_{\rho+\eps} \Subset  U_{\rho+\delta}$. Such a function $\vp$ exists since $U_{\rho+\eps}$ is precompact and $(\sG^{(\rho)},\sH)$ is regular. Define 
$$u(s,x):= Q^{(\rho),D}_s f(x) - Q^{(\rho),U}_s f(x) \quad \text{and} \quad  w(s,x):=u(s,x) \vp(x).$$
Note that for all $s>0$, $ \lVert u(s,\cdot) \rVert_{L^\infty(\X)} \le \lVert Q^{(\rho),D}_s f\rVert_{L^\infty(\X)} \le \lVert f \rVert_{L^\infty(\X)}$ and $u(s,\cdot) \in \sH^{(\rho),D}$.
Thus, for any $s>0$,  by \cite[Theorem 1.4.2(ii)]{FOT},  $w(s,\cdot ) \in \sH$. Besides, by \cite[Lemma 1.3.4(i)]{FOT}, we have
\begin{align}\label{e:comparison-1}
		\left\langle  \partial_s u(s,\cdot), \psi \right\rangle = - \sG^{(\rho)}(u(s,\cdot),\psi) \quad \text{for all $s>0$ and $\psi \in \sH^{(\rho),U}$}.
\end{align}
 Observe that  $\sG^{(\rho)}(\xi,\psi)=0$ for any $\xi,\psi \in \sH$ with dist(supp$[\xi]$, supp$[\psi])>\rho$.  Thus, since supp$[(\vp-1)u(s,\cdot)]\subset$ supp$[\vp-1] \subset U_{\rho+\eps}^c$ and dist$(U_{\rho+\eps}^c, U) \ge \rho + \eps>\rho$, 
\begin{align}\label{e:comparison-2}
\sG^{(\rho)}( (\vp-1)u(s,\cdot), \psi )=0 \quad \text{for all $s>0$ and $\psi \in \sH^{(\rho),U}$}.
\end{align}
Using \eqref{e:comparison-1} and \eqref{e:comparison-2}, we obtain for any $s>0$ and $ \psi \in \sH^{(\rho),U}$,
\begin{equation}\label{e:comparison-3}
	\left\langle  \partial_s w(s,\cdot), \psi \right\rangle = 	\left\langle  \partial_s u(s,\cdot), \psi \right\rangle  = -\sG^{(\rho)}\big(w(s,\cdot) -(\vp-1)u(s,\cdot),\psi\big)= -\sG^{(\rho)}(w(s,\cdot),\psi).
\end{equation}

   Fix $t>0$ and a cutoff function $\xi \in \sH$ for $U_{\rho+\delta} \Subset \X$. Define
\begin{align}\label{e:comparison-m}
	a=a(t)&:=\sup_{0<s\le t} \lVert Q_s^{(\rho),D}f\rVert_{L^\infty(U_{\rho+\delta}\setminus K)},\\
	 g(s,x)&:=\big( 1- Q^{(\rho),U}_s \1_U(x)\big)\xi(x) \quad \text{ and } \quad v:= w -ag. \nn
\end{align}
As $s\to 0$, by the strong continuity of the semigroups $(Q^{(\rho),D}_s )_{s\ge 0}$ and  $(Q^{(\rho),U}_s )_{s\ge 0}$,      $w(s,\cdot)$ converges to $0$ in $L^2(U)$.  Thus, $v_+(s,\cdot)$  converges to $0$ in $L^2(U)$ as $s\to 0$. We claim that
\begin{align}\label{e:comparison-requirement}
v_+(s,\cdot) \in \sH^{(\rho),U} \quad \text{for all $0<s\le t$}.
\end{align}
To obtain \eqref{e:comparison-requirement},   by \cite[Theorem 4.4.3]{FOT},  it suffices to show that for any $s\in (0,t]$, 
\begin{align}\label{e:comparison-requirement-claim}
	w(s,x)\le ag(s,x) \quad \text{for $\sG^{(\rho)}$-q.e. $x \in U^c$}.
\end{align}
Let $s \in (0,t]$. Since $u(s,\cdot) \in \sH^{(\rho),D}$, by \cite[Theorem 4.4.3]{FOT}, we have
\begin{align}\label{e:comparison-requirement-claim-1}
	w(s,x)\le u(s,x)=0\le a g(s,x) \quad \text{for $\sG^{(\rho)}$-q.e.  $x \in D^c$.} 
\end{align}
Besides,  since $\vp=0$ in $U_{\rho+\delta}^c$,    $w(s,x) = 0 \le a \xi(x)$   for $m$-a.e. $x\in D\setminus U_{\rho+\delta}$. Furthermore, by  \eqref{e:comparison-m}, 
 $	 w(s,x) \le Q_s^{(\rho),D}f(x) \le a = a \xi (x)$ for $m$-a.e. $x\in U_{\rho+\delta}\setminus K$. Hence,  $w(s,x) \le a\xi(x)$ for $m$-a.e. $x\in D\setminus K$. 
   Note that $	 w(s,\cdot) -   a\xi$ is an element of $\sH$, so it is $\sG^{(\rho)}$-quasi-continuous by our standing assumption. Since $D\setminus K$ is open, by \cite[Lemma 2.1.4]{FOT}, we obtain  $w(s,x) \le a\xi(x)$ for $\sG^{(\rho)}$-q.e. $x\in D\setminus K$. 
  Since $Q^{(\rho),U}_s\1_U \in \sH^{(\rho),U}$ so that $Q^{(\rho),U}_s\1_U(x) = 0$ for $\sG^{(\rho)}$-q.e. $x \in U^c$, it follows that
  \begin{align*}
  	w(s,x) \le a\xi(x) = a g(s,x) \quad \text{for $\sG^{(\rho)}$-q.e. $x \in D\setminus U$}.
  \end{align*}
  Combining this with \eqref{e:comparison-requirement-claim-1}, we obtain \eqref{e:comparison-requirement-claim}. This  finishes the proof of  \eqref{e:comparison-requirement}.
  
By \cite[Lemma 1.3.4(i)]{FOT},  $s\mapsto g(s)|_U$ is weakly differentiable in $L^2(U)$ at all $s>0$ and, for all $s>0$ and $0 \le \psi \in \sH^{(\rho), U}$,\begin{align*}
	\sG^{(\rho)} (g(s,\cdot), \psi ) & =   \lim_{h\to 0} \frac{1}{h} \la 1-  Q_{h}^{(\rho)}  \1_\X  + Q_h^{(\rho)} Q_s^{(\rho),U}\1_U   - Q^{(\rho),U}_s \1_U, \psi \ra  \\
	&\ge  \lim_{h\to 0} \frac{1}{h} \la   Q_{s+h}^{(\rho),U}\1_U   - Q^{(\rho),U}_s \1_U, \psi \ra  =\la \partial_s Q_s^{(\rho),U}\1_U,\psi \ra = - \la \partial_s  g(s,\cdot), \psi\ra.
\end{align*}
Combining this with \eqref{e:comparison-3}, we get that for all $s>0$ and  $0 \le \psi \in \sH^{(\rho),U}$,
	\begin{align*}
		\la \partial_s  v(s,\cdot), \psi \ra + \sG^{(\rho)}( v(s,\cdot),\psi) = -a \big( 	\la \partial_s  g(s,\cdot), \psi \ra + \sG^{(\rho)}( g(s,\cdot),\psi)\big)  \le 0.
	\end{align*}

Now, by applying Proposition \ref{p:maximum}, we conclude that for $m$-a.e. $x \in U_\rho$,
\begin{align*}
 u(t,x) - a( 1- Q^{(\rho),U}_t \1_U(x)) =  v(t,x) \le 0.
\end{align*} 
The proof is complete. \qed

\begin{cor}\label{c:comparison-entrance-2}
	Let  $\rho>0$ and $U,V,W,D$ be open subsets of $\X$ such that $U\subset V$, $V_{2\rho} \Subset W$, $W \subset D$ and $D_{\rho+\delta}$ is precompact for some $\delta\in (0,\rho/4)$.	Suppose that $(\sG^{(\rho)},\sH)$ is a regular Dirichlet form on $L^2(\X)$.
	Then for any  $0\le f \in L^\infty(D)$ with {\rm supp}$[f]\subset \overline U$, for all $t>0$ and $m$-a.e. $x \in D$,
	\begin{align}\label{e:comparison-entrance-2-claim}
		Q^{(\rho),D}_t f(x) \le Q^{(\rho),W}_t f(x) + \left( 1- Q^{(\rho),W}_t \1_W(x)\right) 	\sup_{0<s\le t} \lVert Q^{(\rho),D}_{s} f\rVert_{L^\infty(V_{2\rho}\setminus \overline V)}.
	\end{align}
\end{cor}
\pf Let $t>0$ and $0\le f \in L^\infty(D)$ with supp$[f]\subset \overline U$.  
Set $E:=D\setminus \overline{V_{\frac32\rho}}$ and $K:=D\setminus V_{2\rho}$. Then $K$ is closed relative to the topology on $D$ and $K\subset E$.

By Proposition \ref{p:comparison-entrance}, we obtain   for  $m$-a.e. $x \in W_{\rho}$,
\begin{align}\label{e:comparison-entrance-2-0}
	Q^{(\rho),D}_t f(x) \le Q^{(\rho),W}_t f(x) + \left( 1 - Q^{(\rho),W}_t \1_W(x) \right) \sup_{0<s \le t} \lVert 	Q^{(\rho),D}_{s} f\rVert_{L^\infty(W_{\rho+\delta} \setminus \overline{V_{2\rho}})}.
\end{align}
Since supp$[f]\subset \overline U \subset \overline V$ and $\overline E \cap \overline V= \emptyset$, we have
\begin{align}\label{e:comparison-entrance-2-1}
	Q^{(\rho),E}_t f= 0 \quad \text{in $D$}.
\end{align}
Applying Proposition \ref{p:comparison-entrance} and using \eqref{e:comparison-entrance-2-1}, we obtain for all $s\in (0,t]$, 
\begin{align}\label{e:comparison-entrance-2-2}
	&\lVert 	Q^{(\rho),D}_{s} f\rVert_{L^\infty(W_{\rho+\delta} \setminus \overline{V_{2\rho}})} = \lVert 	Q^{(\rho),D}_{s} f\rVert_{L^\infty((W_{\rho+\delta}\cap D) \setminus \overline{V_{2\rho}})} \le	\lVert 	Q^{(\rho),D}_{s} f\rVert_{L^\infty(D \setminus \overline{V_{2\rho}})} \nn\\
	&\le \lVert Q^{(\rho),E}_t f \rVert _{L^\infty(D \setminus \overline{V_{2\rho}})} + 	\sup_{0<s'\le s} \lVert Q^{(\rho),D}_{s'} f\rVert_{L^\infty(E_{\rho+\delta}\setminus K)}= \sup_{0<s'\le t} \lVert Q^{(\rho),D}_{s'} f\rVert_{L^\infty(E_{\rho+\delta}\setminus K)}.
\end{align}  Note that   $D\cap E_{\rho+\delta} \subset D\setminus \overline V$ and $(D\setminus \overline V) \setminus K = V_{2\rho} \setminus \overline V$. Thus, we get
\begin{equation}\label{e:comparison-entrance-2-3}
	\sup_{0<s'\le t} \lVert Q^{(\rho),D}_{s'} f\rVert_{L^\infty(E_{\rho+\delta}\setminus K)} \le 	\sup_{0<s'\le t} \lVert Q^{(\rho),D}_{s'} f\rVert_{L^\infty((D\setminus \overline V)\setminus K)} =	\sup_{0<s'\le t} \lVert Q^{(\rho),D}_{s'} f\rVert_{L^\infty(V_{2\rho}\setminus \overline V)}.
\end{equation}
Combining  \eqref{e:comparison-entrance-2-0}, \eqref{e:comparison-entrance-2-2} and \eqref{e:comparison-entrance-2-3}, we conclude that \eqref{e:comparison-entrance-2-claim} holds for $m$-a.e. $x \in W_\rho$.

Besides, using Proposition \ref{p:comparison-entrance},   \eqref{e:comparison-entrance-2-1} and \eqref{e:comparison-entrance-2-3}, we get that for $m$-a.e. $x \in D\setminus W_\rho$, 
\begin{align*}
	Q^{(\rho),D}_t f(x) &\le Q^{(\rho),E}_t f(x) + \sup_{0<s \le t} \lVert 	Q^{(\rho),D}_{s} f\rVert_{L^\infty(E_{\rho+\delta} \setminus K)}\le	\sup_{0<s\le t} \lVert Q^{(\rho),D}_{s} f\rVert_{L^\infty(V_{2\rho}\setminus \overline V)}.
\end{align*} 
Since $Q^{(\rho),W}_t \1_W = 0$ in $W_\rho^c$, this completes the proof. \qed

\end{document}